\documentclass[10pt]{article}

\usepackage{mathtools}
\usepackage{fancyhdr}
\usepackage{amsmath,amsfonts}
\usepackage{graphicx}
\usepackage{color}
\usepackage[footnotesize]{caption}
\usepackage{epstopdf}
\usepackage[Lenny]{fncychap}
\usepackage{algorithm}
\usepackage{algorithmic}
\usepackage{tocvsec2}

\usepackage[cal=boondoxo]{mathalfa}
\newcommand{\iu}{\mathbf{i}}
\renewcommand{\vec}[1]{\mathbf{#1}}

\usepackage{tikz}
\usetikzlibrary{positioning}
\usetikzlibrary{patterns}
\usepackage[european,nooldvoltagedirection]{circuitikz}
\usepackage{pgfplots}
\pgfplotsset{compat=1.14}
\usepgfplotslibrary{fillbetween}

\newtheorem{lemma}{Lemma}
\newtheorem{proposition}{Proposition}
\newtheorem{defn}{Definition}

\usepackage[normalem]{ulem}

\makeatletter
\renewcommand*\env@matrix[1][*\c@MaxMatrixCols c]{%
  \hskip -\arraycolsep
  \let\@ifnextchar\new@ifnextchar
  \array{#1}}
\makeatother

\makeatletter

\newcommand{\Rmnum}[1]{\expandafter\@slowromancap\romannumeral #1@}
\makeatother

\newtheorem{remark}{Remark}

\title{A ghost-point smoothing strategy for geometric multigrid on curved boundaries}
\author{
Armando Coco\thanks{Dipartimento di Matematica e Informatica -- Università di Catania -- email:armando.coco@unict.it}
\and
Mariarosa Mazza\thanks{Dipartimento di Scienza e Alta Tecnologia -- Università dell'Insubria -- email:mariarosa.mazza@uninsubria.it}
\and
Matteo Semplice\thanks{Dipartimento di Scienza e Alta Tecnologia -- Università dell'Insubria -- email:matteo.semplice@uninsubria.it}
}

\begin{document}

\maketitle

\begin{abstract}
We present a Boundary Local Fourier Analysis (BLFA) to optimize the relaxation parameters of boundary conditions in a multigrid framework.
The method is implemented to solve elliptic equations on curved domains embedded in a uniform Cartesian mesh,
although it is designed to be extended for general PDEs in curved domains, wherever a multigrid technique can be implemented.
The boundary is implicitly defined by a level-set function and a ghost-point technique is employed to treat the boundary conditions.
Existing strategies in literature adopt a constant relaxation parameter on the whole boundary. In this paper, the relaxation parameters are optimized in terms of the distance between ghost points and boundary, with the goal of smoothing the residual along the tangential direction.
Theoretical results are confirmed by several numerical tests in 1D, 2D and 3D, showing that the convergence factor associated with the smoothing on internal equations is not degraded by boundary effects.
\end{abstract}


\section{Introduction}
Multigrid (MG) method is among the most efficient iterative solver for discrete Partial Differential Equations (PDEs). Initially designed for elliptic equations, MG methods have been widely adopted to solve countless problems thanks to its optimal computation cost that scales linearly with respect to the number of computational nodes for sparse matrices, outperforming many other numerical solvers. Although initially MG approaches were proposed for simple domains such as rectangles in 2D and cuboids in 3D, the application to more complex-shaped domains gained an increasing interest in the scientific community in the last years in several contexts~\cite{gibou2018review}. 

For example, in the context of energy storage problems one of the recent challenges in computer simulations is the representation of 3D microstructure effects of lithium-ion batteries. It is common practice to represent the complex geometry of this class of problems by adopting simplified approaches, such as union of spheres~\cite{wiedemann2013effects}. In this way, the complexity of the geometry is more manageable, though it comes at the expense of accuracy.
Representing complex-shaped objects is a bottleneck to countless applications in scientific computing, ranging from micro to macro scale problems, such as fluid dynamics past complex capillary structure in cardiovascular systems or hydrothermal circulations in subsurface reservoirs.  

In the context of MG for complex-shaped domains, it is well known that boundary effects can propagate from the boundary to the overall domain and then degrade the efficiency of the entire MG unless a proper treatment is adopted. To this purpose, it is common practice in literature to add extra relaxations on the boundary conditions in order to improve the accuracy near the boundary and maintain the overall performance acceptable~\cite{Brandt1984}. This approach has been adopted, for example, in~\cite{coco2020multigrid, coco2018second, CocoRusso:Elliptic}. The rationale is that the extra computational effort required to perform the local boundary relaxations is negligible with respect to the overall internal relaxations and vanishes when $h \rightarrow 0$, where $h$ is the spatial step~\cite{Trottemberg:MG}. This rationale is well justified for computational grids that have an almost uniform resolution in the whole domain. However, several real-world problems require a more refined mesh close to the boundary to capture the curvature of complex-shaped domains and to maintain the overall numerical error almost uniform everywhere. In these cases, the computational cost of extra relaxations along the boundary might dominate the overall cost, since the majority of the mesh is concentrated in the surrounding area of the boundary, resulting in a sub-optimal solver.

A way to overcome this difficulty is to adopt Algebraic Multigrid (AMG) methods~\cite{ruge1987algebraic}, in which the MG operators are defined merely from the linear system. We refer the reader to~\cite{xu2017algebraic} for a complete review on the a 
general theory of AMG methods. The two main approaches used to setup the hierarchy of prolongation matrices in a completely algebraic setting are classical coarsening
and coarsening by aggregation; see e.g.~\cite{briggs2000multigrid,notay1}. The main advantage of AMG methods is that they only requires matrices and can
thus be easily integrated into existing
code (see \cite{d2021amg} for a recent implementation of parallel coarsening algorithms based on aggregation of unknowns). On the other hand, in AMG methods one cannot easily incorporate a priori knowledge about the model
problem, and more efficient solvers could be designed if geometrical information is included into the definition of the operators. In the context of MG methods for curved boundaries, this was highlighted in~\cite{adams2005comparison}, where a comparison between geometric and algebraic approaches is performed.

Complex-shaped domains are usually modelled by fitted-boundary methods such as the Finite Element Method (FEM)~\cite{lozovskiy2018quasi, shaidurov2013multigrid, de2020multigrid, danilov2017finite}, that has been by far the dominant method to represent the geometry of domains, thanks to the use of isoparametric elements that allow a greater flexibility in the domain discretization. However, there are critical limitations in FEM when applied to highly complex-shaped domains: (i) the mesh generation process may become cumbersome when generating elements to conform highly varying curvatures of the boundary and (ii) the design and implementation of parallel solvers in FEM require impractical implementation efforts to partition the mesh for a balanced distribution of the computational load between cores. For this reason, unfitted-boundary methods such as those based on Finite Difference Methods (FDM) and level-set approaches~\cite{gibou2018review} have become more popular in recent years, since they allow a natural design of parallel solvers in the context of High Performance Computing and at the same time they do not require any mesh generation effort, since the domain is immersed into a fixed grid.

FDM in non-rectangular domains on unfitted grid is usually approached by the Ghost-Point method.
The idea appeared for the first time in the Immersed Boundary Method~\cite{Peskin:IBM} to model blood flows in the heart, subsequently extended by LeVeque and Li in~\cite{LeVequeLi:IIM} to higher accuracy.
More recent development can be found in~\cite{Fedkiw:GFM, Gibou:Ghost, Gibou:fourth_order, Gibou:fluid_solid}.

To capture the multi-scale nature of the problem, a very effective and popular approach is the Adaptive Mesh Refinement (AMR), where a uniform Cartesian grid is refined according to a recursive criterion and the entire grid structure is effectively stored in a quad-tree in 2D and oct-tree in 3D. MG methods for AMR (and its parallel implementation) have been elected by the community in scientific computing to be very efficient solvers for PDEs in highly complex-shaped geometries~\cite{dubey2014survey, teunissen2018afivo,runnels2020massively:arxiv}. However, the treatment of boundary conditions in MG methods has been so far undervalued, leading to suboptimal solvers. 

In this paper, we improve the technique proposed in~\cite{CocoRusso:Elliptic} and propose an efficient and accurate technique to embed boundary conditions in multigrid techniques on curved domains without the necessity of extra-relaxations on the boundary equations. The study is proposed for uniform Cartesian grids, providing a building-block to future implementations of MG for AMR, where the main benefit in terms of computational cost will be observed. 

Following the approach proposed by one of the authors in~\cite{CocoRusso:Elliptic} in the framework of finite-difference unfitted-boundary methods, a ghost-point technique will be adopted to enforce the boundary conditions, together with a relaxation technique for boundary equations which is based on adopting a fictitious time approach. For example, a Dirichlet boundary condition $u=g$ is relaxed by the fictitious time evolution $\displaystyle \frac{\partial u}{\partial \tau} = g-u$. In~\cite{CocoRusso:Elliptic}, the fictitious time step was determined by a CFL stability condition to guarantee the convergence of the smoother. However, after a few iterations the boundary relaxations spoiled the smoothness of the solution and degraded the overall performance, needing some extra-relaxations along the boundary to maintain an appropriate overall performance. To overcome this issue, an improved technique will be proposed in this paper by adopting a {\it local} fictitious time step that varies along the boundary for each ghost point and that depends on geometrical properties like the distance of the ghost point from the actual boundary.
The local fictitious time step will be determined by optimising the boundary relaxation parameters through a Boundary Local Fourier Analysis (BLFA) in order to maximise the smoothing property along the tangential direction.
In this way, the overall performance is not degraded and then the MG will maintain an efficiency  similar to that of classical MG on rectangular domains.
The approach will be described in 1D, 2D and 3D. In the latter case, the implementation will be performed in {\tt C++} using {\tt PETSc} \cite{petsc-user-ref,petsc-efficient} to allow for a more natural parallel implementation. Our MG Poisson solver will be compared with the AGgregation-based algebraic MultiGrid (AGMG) of \cite{notay-sw, notay1}.

The paper is organized as follows. In \S\ref{sec:1D} we pose the problem in the 1D setting. Such case, that can indeed be easily managed by eliminated boundary conditions, is intended as starting point for the presentation and discussion of our approach in 2D and 3D settings given in \S\ref{sec:2D} and \S\ref{sect:3D}, respectively. We end with \S\ref{sec:concl} where we draw conclusions.

\section{One-dimensional case} \label{sec:1D}

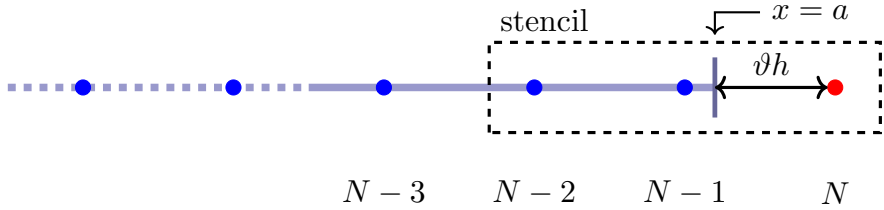
\begin{figure}
\begin{center}
\begin{tikzpicture}[scale=2, every node/.style={scale=1.2}]
\definecolor{colordomain}{rgb}{0.6,0.6,0.8}
\definecolor{colorborder}{rgb}{0.4,0.4,0.6}
\draw[line width=1mm,colordomain] (-0.5,0) -- (2.2,0);
\draw[line width=1mm,colordomain,dashed] (-2.5,0) -- (-0.5,0);
\draw[line width=1.8,color=colorborder] (2.2,-0.2) -- (2.2,0.2);
\foreach \x in {-2,-1,0,1,2}
\draw[blue, fill=blue] (\x,0) circle (0.05);
\draw[red, fill=red] (3,0) circle (0.05);
\definecolor{colordashedbox}{rgb}{0.5569,0.2941,0.1137}
\draw[very thick,dashed,color=black] (0.7,-0.3) rectangle (3.3,0.3);
\node[anchor=south west] at (0.7,0.3) {stencil};
\draw[very thick,<->,color=black] (2.2,0.) -- 
     node[pos=0.5,above] {$\vartheta h$} 
     (3-0.05,0.);
\draw (3,-0.7) node {$N$};
\foreach \x in {1,2,3}
\draw (3-\x,-0.7) node {$N-\x$};
\draw[<-, thick] (2.2,0.35) -- (2.2,0.5) -- (2.5,0.5)  node[right] {$x=a$};	
\end{tikzpicture}
\end{center}
\caption{Stencil for the boundary conditions in 1D. The ghost point is in red, internal grid points are in blue.}
\label{fig:1Ddomain}
\end{figure}

Consider a one-dimensional uniform grid $\mathcal{G}_h = \left\{ x_0, x_1, \ldots, x_N \right\}$ with spatial step $h = x_{i}-x_{i-1}$, for $i=1,\ldots,N$.
The elliptic boundary-value problem is
\begin{equation}
\begin{split}\label{mainProblem}
- \Delta u &= f \text{ on } \Omega = (x_0,a) \\
u(x_0)&=g_0, \quad Bu(a) = g_a,
\end{split}
\end{equation}
where $Bu(a)=g_a$ is a generic boundary condition on $x=a$, with $x_{N-1}<a \leq x_N$ (see Fig.~\ref{fig:1Ddomain}). We call $\vartheta = (x_N-a)/h$. Then $0 \leq \vartheta < 1$.
In this paper we will discuss Dirichlet ($Bu(a)=u(a)$) and Neumann ($Bu(a)=u'(a)$) boundary conditions. The extension to more general boundary conditions, such as Robin boundary conditions, is straightforward.
The elliptic equation $- \Delta u = f$ is discretized by central differences on $x_i \in \Omega = (x_0,a)$
\begin{equation}\label{Disc1Dx0}
\frac{2u_1- u_{2}}{h^2} = f_1+\frac{g_0}{h^2} 
\end{equation}
\[
\frac{-u_{i-1}+2u_i- u_{i+1}}{h^2} = f_i \text{ for } i=2,\ldots,N-1,
\]
while the boundary condition is discretized by $Bu(a) = Bq(a)$, where $q(x)$ is the polynomial of degree $s-1$ that interpolates $u$ on nodes $u_N, u_{N-1}, \ldots, u_{N-s+1}$, where $s$ is the stencil size. 

In general, the discretization of the boundary condition for $s\leq3$ is
\begin{equation}\label{BC1D}
c_{-2} u_{N-2} + c_{-1} u_{N-1} + c_0 u_{N} = g_a,
\end{equation}
where the coefficients $c_{-k}$ are the values (respectively derivatives) of the Lagrange polynomials associated to the nodes $-s+1,\ldots,0$. In particular:
\begin{itemize}
\item[] $(c_{-2},c_{-1},c_0) = (0, \vartheta, 1-\vartheta)$ for Dirichlet b.c. and $s=2$,
\item[] $(c_{-2},c_{-1},c_0) = \left( \vartheta (\vartheta-1)/2, \; \vartheta (2-\vartheta), \;(1-\vartheta) (2-\vartheta)/2 \right)$ for Dirichlet b.c.~and $s=3$,
\item[]$(c_{-2},c_{-1},c_0) = \frac{1}{h} \left( 0 -1, 1 \right)$ for Neumann b.c.~and $s=2$,
\item[] $(c_{-2},c_{-1},c_0) = \frac{1}{h} \left( (0,\; -1,\; 1) - \left( \vartheta-\frac{1}{2} \right) (1,\; -2,\; 1) ) \right)$ for Neumann b.c.~and $s=3$.
\end{itemize}

The grid point $x_N$ is called a ghost point and $u_N$ is the ghost value. Although the boundary condition can be easily eliminated in 1D by solving~\eqref{BC1D} for $u_N$ and substituting into the internal equations (as it has been done for the boundary condition on $x_0$, see Eq.~\ref{Disc1Dx0}), we maintain a non eliminated approach since we aim at the extension 
to higher dimensions, where the eliminated boundary condition approach is not possible (or extremely impractical). We note that a possible approach based on eliminated boundary conditions in higher dimensions for Dirichlet boundary conditions is proposed in~\cite{Gibou:Ghost}, although the extension to Neumann boundary conditions is not straightforward.

The discretized problem is then a linear system $M_h \textbf{u}_h = \textbf{f}_h$ where $M_h \in \mathbb{R}^{N \times N}$:
\[
M_h = 
\begin{pmatrix}
2/h^2 & -1/h^2 & \ldots & & & & 0 \\
\\
-1/h^2 & 2/h^2 & -1/h^2 & \ldots & & & 0 \\
\\
0 & -1/h^2 & 2/h^2 & -1/h^2 & \ldots & & 0 \\
\vdots \\
\vdots \\
\vdots \\
0 & \ldots & & & -1/h^2 & 2/h^2 & -1/h^2  \\
\\
0 & \ldots & & 0 & c_{-2} & c_{-1} & c_0\\
\end{pmatrix}
\]

\subsection{Multigrid method}\label{MG1D}


We quickly recall the main idea of the multigrid method, referring the reader to more comprehensive literature for more details (see, for example,~\cite{Trottemberg:MG}).
First, a relaxation scheme is applied to the linear system $M_h \textbf{u}_h = \textbf{f}_h$ for a few iterations, say $\nu_1$ \emph{pre-smoothing} iterations, obtaining an approximated solution $\vec{u}^*$.
Commonly, the relaxation scheme is of Richardson type, namely:
\begin{equation}\label{RichRel}
\textbf{u}_h^{(m+1)} = \textbf{u}_h^{(m)} + P^{-1} (\textbf{f}_h-M_h \textbf{u}_h^{(m)})
\end{equation}
for a suitable preconditioner matrix $P$.
The residual $\vec{r}_h=\vec{f}_h-M_h \vec{u}^*_h$ is evaluated on the grid $\mathcal{G}_h$ and transferred to a coarser grid $\mathcal{G}_{2h}$ through a restriction operator $\vec{r}_{2h} = \mathcal{I}^h_{2h} \vec{r}_h$, defined below. Then, the residual equation $M_{2h} \vec{e}_{2h} = \vec{r}_{2h}$ is solved exactly on the coarse grid. The algebraic error $\vec{e}_{2h}$ is transferred back to the fine grid $\mathcal{G}_h$ through an interpolation operator $\vec{e}_{h} = \mathcal{I}^{2h}_{h} \vec{e}_{2h}$. The approximation $\vec{u}^*_h$ is then updated by $\vec{u}^*_h \coloneqq \vec{u}^*_h + \vec{e}_h$ and another few, say $\nu_2$, \emph{post-smoothing} iterations of \eqref{RichRel} are applied to the fine grid. 
The approximation $\vec{u}^*$ represents the output of a single two-grid iteration (cycle). The cycle is iterated until the residual on the fine grid falls under a suitable tolerance. 

The above approach is also called Two Grid Correction Scheme (TGCS), since only two grids are involved in the process, $\mathcal{G}_h$ and $\mathcal{G}_{2h}$. More grids can be adopted if the residual equation $M_{2h} \vec{e}_{2h} = \vec{r}_{2h}$ is not solved exactly, but is instead approximated by $\gamma$ cycles of the TGCS on grids $\mathcal{G}_{2h}$ and $\mathcal{G}_{4h}$, and so on recursively, until a very coarse grid $\mathcal{G}_{H}$ is reached and where the residual equation can be solved exactly at a low computational cost. Depending on the value of $\gamma$, we have the V-cycle ($\gamma=1$) or W-cycle ($\gamma=2$). It is uncommon to use higher values of $\gamma$ since the it would not be convenient in terms computational effort~\cite{Trottemberg:MG}.

We now describe the relaxation scheme and the transfer operators that we adopt for the multigrid strategy in 1D.

\subsubsection*{Relaxation scheme}
In order to have an effective multigrid strategy, the relaxation scheme \eqref{RichRel} should be able to quickly reduce the high-oscillatory modes of the error $e_h$ ({\em smoothing property}), regardless of the the reduction of low-oscillatory modes. A relaxation scheme with this property is called a {\it smoother}.
It is well known that the Jacobi scheme (when $P=\text{diag}(M_h)$ is the diagonal part of $M_h$) is not a proper smoother.
Relaxation schemes that have the smoothing property are weighted-Jacobi
$P=\omega \text{diag}(M_h) + (1-\omega) \mathbb{I}$ with $\omega=2/3$ in 1D, where $\mathbb{I}$ is the identity matrix, and the Gauss-Seidel scheme, obtained when the preconditioner matrix $P=L_h$ is the lower triangular part of $M_h$ (including the diagonal).
It can be proved that the convergence of both Jacobi and Gauss-Seidel schemes is not guaranteed for the presence of the last equation in the linear system, especially when $\vartheta > 1/2$. 
For example, when $\vartheta$ is close to $1$ then the diagonal entry $c_0$ of $M_h$ is in magnitude much larger than off-diagonal entries $c_{-1},c_{-2}$ and we observe numerically that the iteration scheme is not convergent.

To maintain the convergence, the ghost value $u_N$ is not updated by the classical Gauss-Seidel scheme. It is instead obtained by a relaxation of the boundary conditions on $x=a$:
\begin{equation}\label{eq:relghostdt}
u_N^{(m+1)} = u_N^{(m)} + \Delta \tau \, (g_a - (c_{-2} \, u^{(m+1)}_{N-2} + c_{-1} \, u^{(m+1)}_{N-1}+ c_0 u^{(m)}_N)).
\end{equation}
for some fictitious positive relaxation parameter $\Delta \tau>0$ to be chosen in order to ensure the convergence of the entire iteration scheme. Observe that \eqref{eq:relghostdt} reverts to the classical Gauss-Seidel scheme if $\Delta \tau= 1/c_0$.

With \eqref{eq:relghostdt}, the preconditioner matrix $P$ is the lower triangular part of $M_h$ except for $P_{NN}$ that is replaced by
\[
P_{NN} = \frac{1}{ \Delta \tau}.
\]
We have used the symbol $\Delta \tau$ because the iteration \eqref{eq:relghostdt} can be also seen as the discretization of the evolutionary differential equation:
\begin{equation}\label{eq:pdefict}
\frac{\partial u}{\partial \tau} = g_a - u \text{ in } x=a
\end{equation}
in a fictitious time $\tau$ (see~\cite{CocoRusso:Elliptic}), where the time derivative in the right hand side is discretized on $x=x_N$ by the explicit Euler scheme and the left hand side on $x=a$. The steady-state solution of \eqref{eq:pdefict} satisfies the boundary condition $u(a)=g_a$.

\subsubsection*{Restriction operator}
Following the approach of~\cite{CocoRusso:Elliptic}, the transfer operators (restriction $\mathcal{I}^h_{2h}$ and interpolation $\mathcal{I}^{2h}_{h}$) are defined as follows. 
The residual $\vec{r}_h$ contains the residual of internal equations for $x_1\leq x_ i \leq x_{N-1}$ and that of the boundary condition for the ghost point $x_N$. These two sub-residuals scale at different powers of $h$: internal equations scale at $h^2$, while the relaxation of the boundary condition scales at $h^0$ or $h$ for Dirichlet and Neumann boundary conditions, respectively. The restriction is then discontinuous across the boundary $x=a$ and for this reason, the ghost value $r_{N+1}$ is not used in the restriction of internal equations. In detail, if $x+h < a $ the restriction is Full-Weighting (FW):
\begin{equation}\label{eq:coarseFW}
\vec{r}_{2h} (x) = \left( \mathcal{I}^h_{2h} \vec{r}_h \right) (x) = \frac{1}{4} \left[ 1 \quad 2 \quad 1 \right] \times
\begin{bmatrix}
\vec{r}_h (x-h) \\
\\
\vec{r}_h (x) \\
\\
\vec{r}_h (x+h) \\
\end{bmatrix}
=
\frac{\vec{r}_h (x-h) + 2 \vec{r}_h (x) + \vec{r}_h (x+h)}{4},
\end{equation}
while it is the average of the internal values if $x<a$ but $x + h \geq a$:
\[
\vec{r}_{2h} (x) = \left( \mathcal{I}^h_{2h} \vec{r}_h \right) (x) = \frac{1}{2} \left[ 1 \quad 1 \quad 0 \right] \times
\begin{bmatrix}
\vec{r}_h (x-h) \\
\\
\vec{r}_h (x) \\
\\
\vec{r}_h (x+h) \\
\end{bmatrix}
=
\frac{\vec{r}_h (x-h) + \vec{r}_h (x)}{2}.
\]
The restriction of the residual of the boundary condition is trivial:
\[
\vec{r}_{2h} (x_N) =  \vec{r}_{h} (x_N).
\]

\subsubsection*{Interpolation operator}
The interpolation operator acts on the error $\vec{e}_{2h}$, which is continuous across the boundary, and therefore no modification from the standard linear interpolation is needed:
\[
\vec{e}_{h} (x) = \left( \mathcal{I}^{2h}_{h} \vec{e}_{2h} \right) (x) = \vec{e}_{2h} (x) \text{ if } x \in \mathcal{G}_{2h}
\]
\[
\vec{e}_{h} (x) = \left( \mathcal{I}^{2h}_{h} \vec{e}_{2h} \right) (x) = \frac{\vec{e}_{2h} (x-h) + \vec{e}_{2h} (x+h)}{2} \text{ if } x \notin \mathcal{G}_{2h}.
\]

\begin{defn}
We call \emph{Boundary Condition Multigrid} \textsc{(BCMG}) method the multigrid technique proposed in this paper.
\end{defn}

\subsection{Numerical results in 1D}
The efficiency of the multigrid is measured through the asymptotic convergence factor
\[
\rho = \lim_{m \rightarrow \infty} \rho^{(m)}
= \lim_{m \rightarrow \infty} \frac{ \left\| \vec{r}^{(m+1)}_h \right\|_p }{ \left\| \vec{r}^{(m)}_h \right\|_p }.
\]

In order to avoid numerical issues due to machine precisions, the convergence factor is measured for the homogeneous problem \eqref{mainProblem} with $f=g_0=g_a=0$ and with an initial guess of the multigrid cycle different from zero, say $u^{(0)}_i$ is randomly uniformly distributed between $-1$ and $1$ for all $i$. The asymptotic convergence factor $\rho$ is then approximated by $\rho =  \frac{1}{5} \sum_{k=11}^{15} \rho^{(m)}$.

\begin{figure}
\begin{center}
\begin{tabular}{cc}
(Dirichlet) & (Neumann) \\
\includegraphics[width=0.45\textwidth]{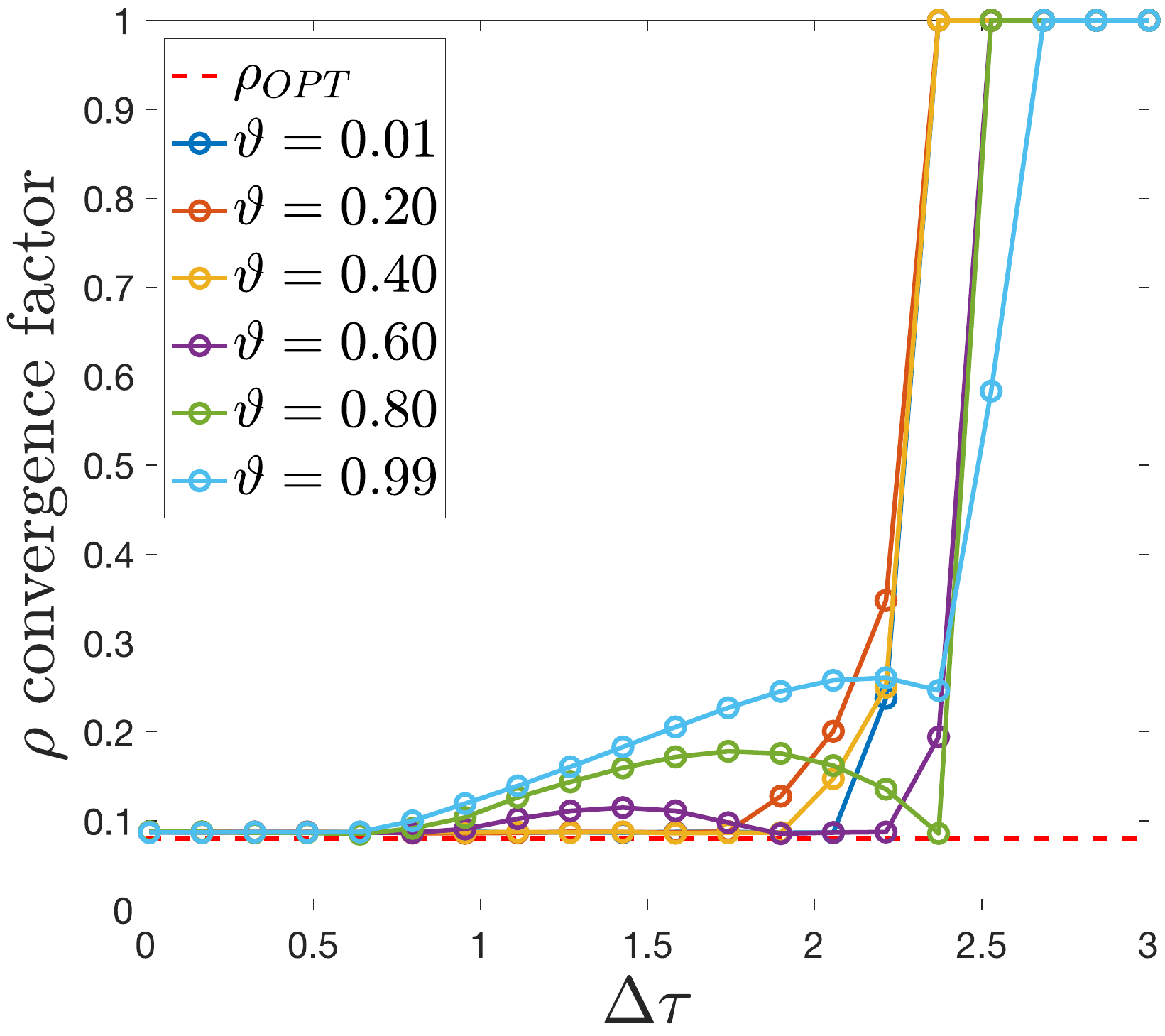}
&
\includegraphics[width=0.45\textwidth]{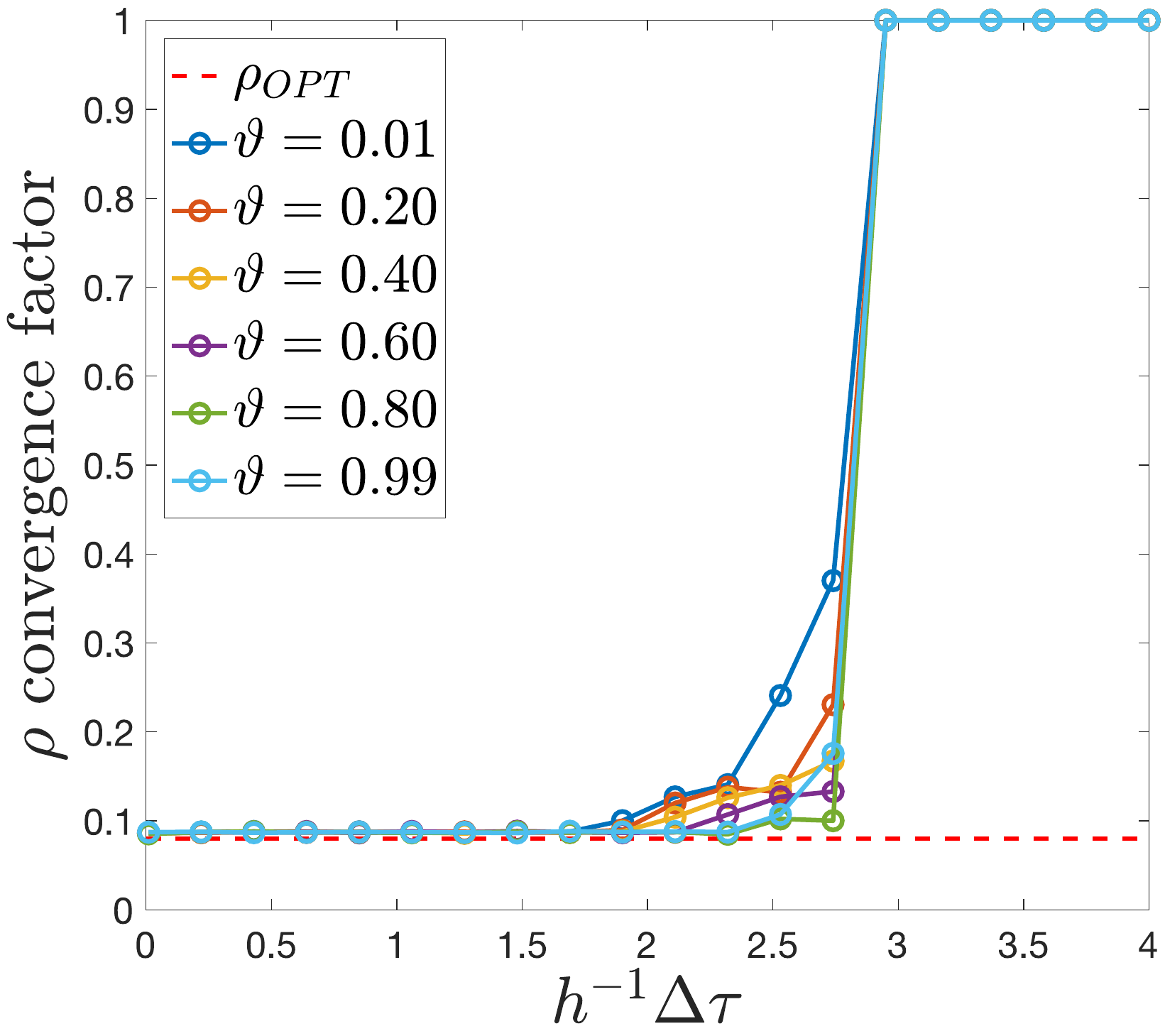}
\end{tabular}
\end{center}
\caption{TGCS asymptotic convergence factor in one space dimension varying the relaxation parameter $\Delta\tau$. Left: Dirichlet boundary condition. Right: Neumann boundary condition.
The red dashed line is the optimal convergence factor without ghost nodes.}
\label{fig:1D_rho}
\end{figure}

In Fig.~\ref{fig:1D_rho} we plot the convergence factors versus $\Delta \tau$ for different values of $\vartheta$ for the TGCS with $N=512$, $\nu_1=2$ and $\nu_2=1$, for Dirichlet (left) and Neumann (right) boundary conditions. We report the results obtained with $s=3$, but very similar results are obtained for $s=2$. The dashed line represents the optimal convergence factor $\rho_{OPT} \approx 0.08$, that is the convergence factor achieved when $\vartheta=0$ and eliminated boundary conditions are implemented~\cite{Trottemberg:MG} (i.e.~the convergence factor of internal relaxations only). If the multigrid does not converge, we cut the convergence factor at $\rho=1$ for graphical purposes. We observe that $\rho$ is close to the optimal one (the convergence factor is not degraded by boundary effects of the relaxation of boundary conditions) or the multigrid does not converge (when $\Delta \tau$ is too large). 
The plot shows  that optimality is reached even for very small values of $\Delta \tau$, indicating that no relaxation is needed for the residual of the boundary condition. We point out however that this is a specific feature of the one-dimensional case, see also \cite{Brandt1984}. In fact in this case, the boundary residual is transferred as it is to the coarser grid and the error equation is solved exactly there. As we will see, this cannot be extended to the case in more space dimensions.

\section{Two-dimensional case} \label{sec:2D}
In higher dimensions, the representation of the residual of boundary conditions on the coarser grid becomes a key point for the efficiency of the multigrid,  since it is not an isolated value as in 1D, but it is instead represented by a set of ghost values just outside the domain and coarsening happens on ghost nodes too. 

For all 2D problems we use a uniform Cartesian grid with $h=\Delta x = \Delta y$, albeit the method can be easily extended to anisotropic grids with $\Delta x \neq \Delta y$. 
Let $\Omega \subseteq \mathcal{D} = [-1,1]^2$ be the computational domain. Namely the grid covers $[-1,1]^2$ with $(N+1)^2$ grid points $X_{ij} = (x_i,y_j)$, where $x_i = -1+i\,h$ and $y_j=-1+j\,h$ for $i,j=0,\ldots,N$ and $h=2/N$.

We consider the elliptic problem
\begin{equation}\label{mainProblem2D}
\begin{split}
- \Delta u &= f \text{ in } \Omega \\
u&=g_\mathcal{B} \text{ on } \Gamma_\mathcal{B} = \partial \mathcal{D} \cap \partial \Omega  \\
Bu &= g_\Gamma \text{ on } \Gamma = \partial \Omega \backslash \Gamma_\mathcal{B}.
\end{split}
\end{equation}

We implement an eliminated boundary conditions approach on $\Gamma_\mathcal{B}$ (similarly to the 1D case \eqref{mainProblem} for $x=x_0$) and a ghost-point technique on $\Gamma$, explained in \S\ref{sect:rectangle}.

\subsection{Set-up of 2D and 3D numerical tests and computation of the convergence factor}\label{sect:setup}
Unless otherwise stated, we perform in all numerical tests the TGCS with $\nu_1=2$, $\nu_2=1$, GS-LEX as a smoother and FW restriction operator.
The convergence factor is measured for the homogeneous problem \eqref{mainProblem2D} with $f=g_\mathcal{B}=g_\Gamma=0$ and with an initial guess of the multigrid cycle different from zero and in particular randomly uniformly distributed $u_{ij} \sim (-1,1)$ for all $i$ and $j$. The asymptotic convergence factor $\rho$ is then approximated by $\rho =  \frac{1}{5} \sum_{k=11}^{15} \rho^{(m)}$.

\subsection{Rectangular domains}\label{sect:rectangle}

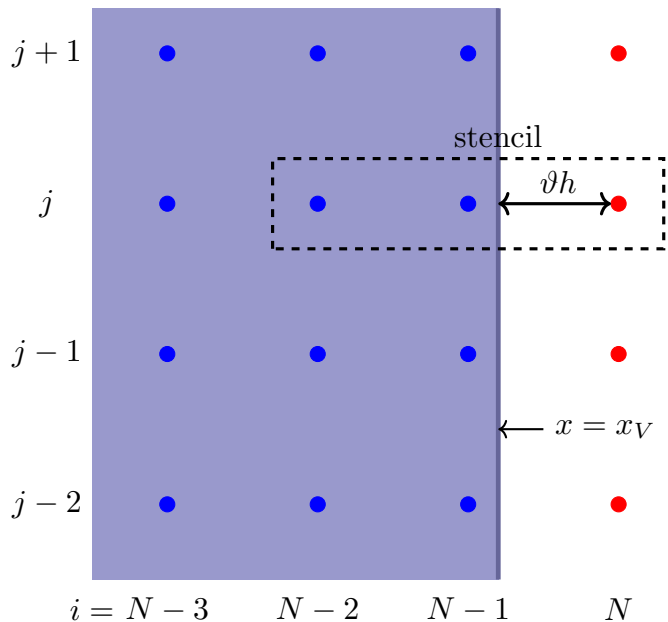
\begin{figure}
\begin{center}
\begin{tikzpicture}[scale=2, every node/.style={scale=1.2}]
\definecolor{colordomain}{rgb}{0.6,0.6,0.8}
\definecolor{colorborder}{rgb}{0.4,0.4,0.6}
\fill[color=colordomain]
(-0.5,-0.5) -- (2.2,-0.5) -- (2.2,3.3) -- (-0.5,3.3) -- cycle;
\draw[line width=1.8,color=colorborder] (2.2,-0.5) -- (2.2,3.3);
\foreach \x in {0,1,2}
\foreach \y in {0,1,2,3}
\draw[blue, fill=blue] (\x,\y) circle (0.05);
\foreach \y in {0,1,2,3}
\draw[red, fill=red] (3,\y) circle (0.05);
\definecolor{colordashedbox}{rgb}{0.5569,0.2941,0.1137}
\draw[very thick,dashed,color=black] (0.7,1.7) rectangle (3.3,2.3);
\draw (2.2,2.45) node {stencil};
\draw[very thick,<->,color=black] (2.2,2) -- (3-0.05,2);
\node[above] at (2.6,2) {$\vartheta h$};
\draw (3,-0.7) node {$N$};
\foreach \x in {1,2,3}
\draw (3-\x,-0.7) node {$N-\x$};
\draw (-0.5,-0.7) node {$i=$};
\foreach \y in {2,1}
\draw (-0.8,2-\y) node {$j-\y$};
\draw (-0.8,2) node {$j$};
\draw (-0.8,3) node {$j+1$};
\draw[<-, thick] (2.2,0.5) -- (2.5,0.5)  node[right] {$x=x_V$};	
\end{tikzpicture}
\end{center}
\caption{Grid setup in the two-dimensional rectangular domain tests.}
\label{fig:2DdomainVLINE}
\end{figure}

Consider the rectangular domain
\[
\Omega = \left\{ (x,y) \colon 0 \leq x \leq x_V, 0\leq y \leq 1 \right\},
\]
where $1-h<x_V \leq 1$.
Setting $\vartheta = (1-x_V)/h$, then $0 < \vartheta \leq 1$.
This is illustrated in Fig.~\ref{fig:2DdomainVLINE}.
For the multigrid method, we use the technique and notations introduced in \S\ref{MG1D} for the 1D case. For the 2D case, we have to define the relaxation operator, and the transfer (restriction and interpolation) operators.

\paragraph*{Relaxation scheme.}
The value on grid points of $\Gamma_\mathcal{B}$ is set from the boundary condition, so $u^{(m)}_{ij} = g_\mathcal{B} (x_i,y_j)$ for any $m$, where $i=0$ or $j \in
\left\{ 0,N \right\}$.

The value on the internal grid points is obtained by the standard Gauss-Seidel lexicographic (GS-LEX) iteration:
\begin{equation}\label{GSLEX_2D}
u_{i,j}^{(m+1)} = \frac{h^2}{4} \, f_{ij} + \frac{u_{i+1,j}^{(m)}+u_{i-1,j}^{(m+1)}+u_{i,j+1}^{(m)}+u_{i,j-1}^{(m+1)}}{4}
,
\quad
\text{for }i,j=1,\ldots,N-1.
\end{equation}

Finally, a fictitious-time relaxation (similar to the 1D case \eqref{eq:relghostdt}) is performed on the boundary $\Gamma$, i.e.  $x=x_V$, for which the ghost points are $X_{Nj}=(x_N,y_j)$ for $j=1,\ldots,N-1$:
\begin{equation}\label{relax2DV}
u_{N,j}^{(m+1)} = u_{N,j}^{(m)} + \Delta \tau \, (g_{j} - (c_{-2} \, u^{(m+1)}_{N-2,j} + c_{-1} \, u^{(m+1)}_{N-1,j}+ c_0 u^{(m)}_{N,j}))
\end{equation}
where $g_{j} = g_\Gamma (x_V,y_j)$.

\paragraph*{Restriction operator: Full-Weighting (FW)}

With the notation $T \, v_{ij}$, where $T$ is a $3 \times 3$ matrix with elements $T_{rs}$ for $r,s,=-1,0,1$, we denote the value:
\[
T v_{ij} = \sum_{r,s=-1}^1 T_{rs} v_{i+r \, j+s}.
\]
Let $u^*_{ij}$ be the approximation after $\nu_1$ relaxations on the fine grid.
The residual after the internal relaxations \eqref{GSLEX_2D}
\[
(\vec{r}_h)_{ij} = f_{ij} - \frac{4 u^*_{ij} - (u_{i+1,j}^{*}+u_{i-1,j}^{*}+u_{i,j+1}^{*}+u_{i,j-1}^{*})}{h^2}
\]
is transferred to the coarser grid by the usual Full-Weighting (FW) 2D nine-point stencil (see e.g.~\cite{Trottemberg:MG} for more details):
\[
(\vec{r}_{2h})_{2i, 2j} = 
\frac{1}{16}
\begin{bmatrix}
1 & 2 & 1 \\
2 & 4 & 2 \\
1 & 2 & 1 
\end{bmatrix}
(\vec{r}_{h})_{2i, 2j}.
\]
The residual after ghost value iterations \eqref{relax2DV}
\[
r_{Nj} = 
g_{j} - (c_{-2} \, u^{*}_{N-2,j} + c_{-1} \, u^{*}_{N-1,j}+ c_0 u^{*}_{N,j})
\]
scales differently than the internal one and is thus transferred via the 1D operator acting only on ghost values:
\[
(\vec{r}_{2h})_{N , 2j} = 
\frac{(\vec{r}_h)_{N,2j-1} + 2 (\vec{r}_h)_{N,2j} + (\vec{r}_h)_{N,2j+1}}{4}.
\]
Note that in the present test case setup, FW on internal equations never mixes internal with ghost residuals.

\paragraph*{Interpolation operator.}
Standard 2D bilinear interpolation operator is applied all over the grid (internal and ghost grid points).

\subsubsection{Dirichlet boundary conditions: Boundary Local Fourier Analysis (BLFA)}\label{2D_dirLFA}
For Dirichlet boundary conditions on $x=x_V$, the coefficients $c_{-2}$, $c_{-1}$ and $c_0$ are defined in \eqref{BC1D}. Assume a stencil size of $s=3$. The problem has $(N-1)^2$ internal points $(x_i,y_j)$ for $i,j=1,\ldots,N-1$ and $N-1$ ghost points $(x_N,y_j)$ for $j=1,\ldots,N-1$, for a total of $(N-1)N$ active points.

The discrete problem can be written as a linear system $M_h \vec{u}_h = f_h$, with:

\begin{equation*}
M_h =
\left(
\begin{array}{ccccccccc | ccc}
L & & D & & & &  & &  & & & \\
D & & L & & D & &  & &  & & & \\
& & D & & L & & D & &  & & & \\
& &  & & \ddots & & \ddots & & \ddots & & & \\
& &  & &  & & D & & L & & D &    \\
\hline
\\
& &  & &  & & c_{-2} {I}_{N-1} & & c_{-1} {I}_{N-1} & & c_0 {I}_{N-1} &     \\
\end{array}
\right)
\end{equation*}
\begin{equation*}
=
\begin{pmatrix}
&|&  \\
M_{II} &|& M_{IG} \\
 &|& \\
\hline
&|&  \\
M_{GI} &|& M_{GG} \\
 &|&  \\
\end{pmatrix}
\in
\qquad
\begin{matrix}
&|&  \\
\mathbb{R}^{(N-1)^2 \times (N-1)^2} &|& \mathbb{R}^{(N-1)^2 \times (N-1)} \\
 &|& \\
\hline
&|&  \\
\mathbb{R}^{(N-1) \times (N-1)^2} &|& \mathbb{R}^{(N-1) \times (N-1)} \\
 &|&  \\
\end{matrix}
\end{equation*}
where 
\begin{equation*}
{I}_{N-1} \in \mathbb{R}^{(N-1)\times(N-1)}
\end{equation*}
is the identity matrix,
\begin{equation*}
D=-\frac{1}{h^2}{I}_{N-1} \in \mathbb{R}^{(N-1)\times(N-1)}
\end{equation*}
and
\begin{equation*}
L = 
\frac{1}{h^2}\begin{pmatrix}
4 & -1 &  &  &    \\
  -1 & 4 & -1 &  &    \\
   & \ddots& \ddots & \ddots &      \\
   &  & -1&4 & -1    \\
   &  &  & -1 & 4  \\
\end{pmatrix} \in \mathbb{R}^{(N-1)\times(N-1)}.
\end{equation*}

The relaxation scheme can be written in matrix form as
\[
\begin{pmatrix}
\vec{u}_I\\
\vec{u}_G
\end{pmatrix}^{(k+1)}
=
\begin{pmatrix}
\vec{u}_I\\
\vec{u}_G
\end{pmatrix}^{(k)}
+
P^{-1}
\left( 
\begin{pmatrix}
\vec{f}_I\\
\vec{f}_G
\end{pmatrix}
-
\begin{pmatrix}
M_{II} & M_{IG} \\
M_{GI} & M_{GG}
\end{pmatrix}
\begin{pmatrix}
\vec{u}_I\\
\vec{u}_G
\end{pmatrix}^{(k)}
\right),
\]
with the preconditioner matrix 
\begin{equation}\label{eq:2Ddiri:P}
P =
\begin{pmatrix}
&|&  \\
L_{II} &|& O_{IG} \\
 &|& \\
\hline
&|&  \\
M_{GI} &|& \displaystyle \frac{1}{\Delta \tau} I_{N-1} \\
 &|&  \\
\end{pmatrix}
,
\end{equation}
where $L_{II} \in \mathbb{R}^{(N-1)^2 \times (N-1)^2}$ is the lower part (including the diagonal) of $M_{II}$ and $O_{IG} \in \mathbb{R}^{(N-1)^2 \times (N-1)}$ is the null matrix. 
To discriminate the contribution of the boundary relaxation from the internal equations, we analyze an alternative relaxation scheme obtained with
\begin{equation}\label{eq:2Ddiri:Pbc}
\widetilde{P} =
\begin{pmatrix}
&|&  \\
M_{II} &|& O_{IG} \\
 &|& \\
\hline
&|&  \\
M_{GI} &|& \displaystyle \frac{1}{\Delta \tau} I_{N-1} \\
 &|&  \\
\end{pmatrix}
.
\end{equation}
In other words, the internal variables are solved in block and
the ghost values $u^{(m+1)}_{N,j}$ are updated by the relaxation \eqref{relax2DV}.
This alternative relaxation scheme is of course impractical and used here only for theoretical studies.

We call {\tt REL-ALL} the TGCS associated with~\eqref{eq:2Ddiri:P} and {\tt REL-BC} the one associated with~\eqref{eq:2Ddiri:Pbc}.

\begin{figure}
\begin{minipage}{0.99\textwidth}
\begin{center}
(Dirichlet)\\
\includegraphics[width=0.79\textwidth]{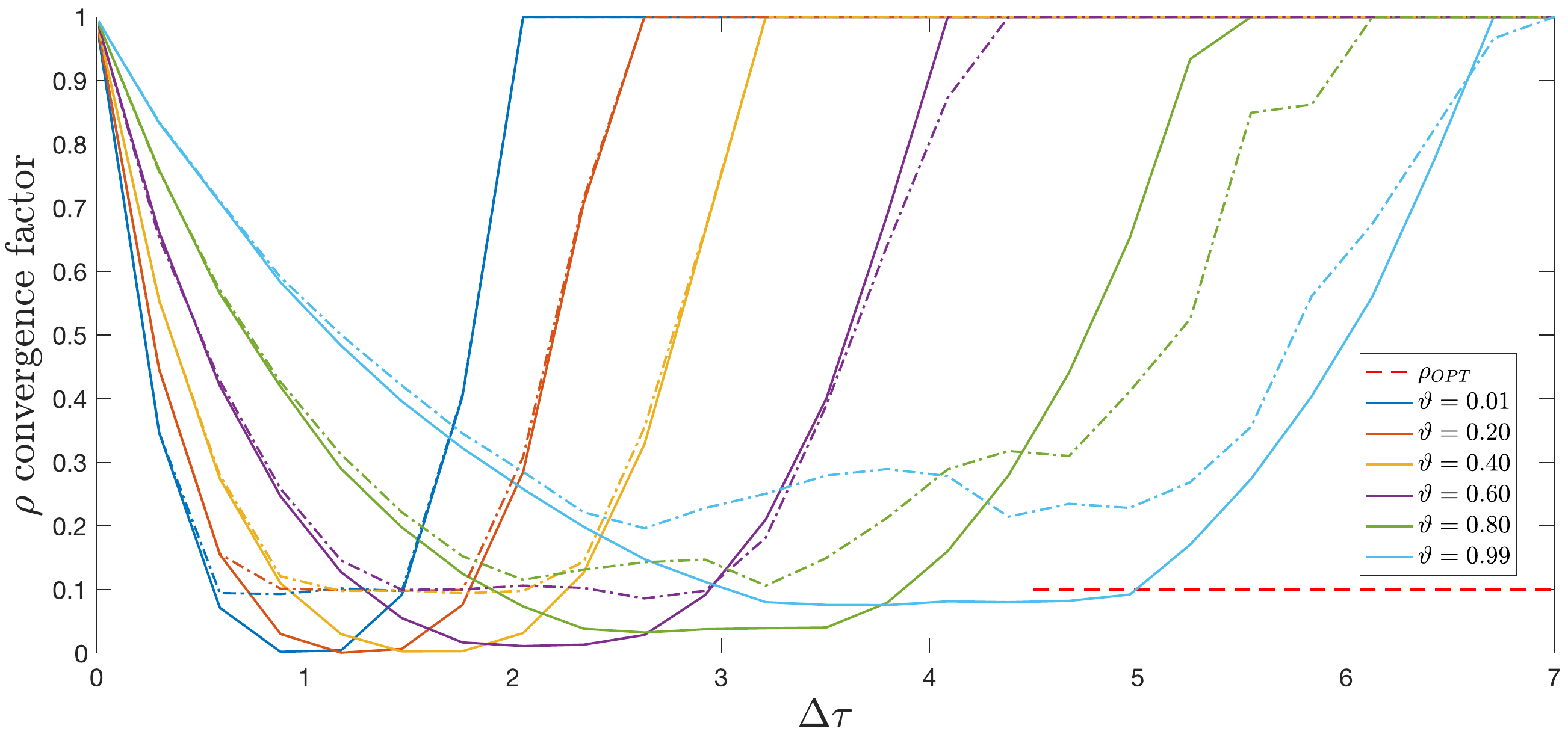}
\end{center}
\end{minipage}
\begin{minipage}{0.99\textwidth}
\begin{center}
(Neumann)\\
\includegraphics[width=0.79\textwidth]{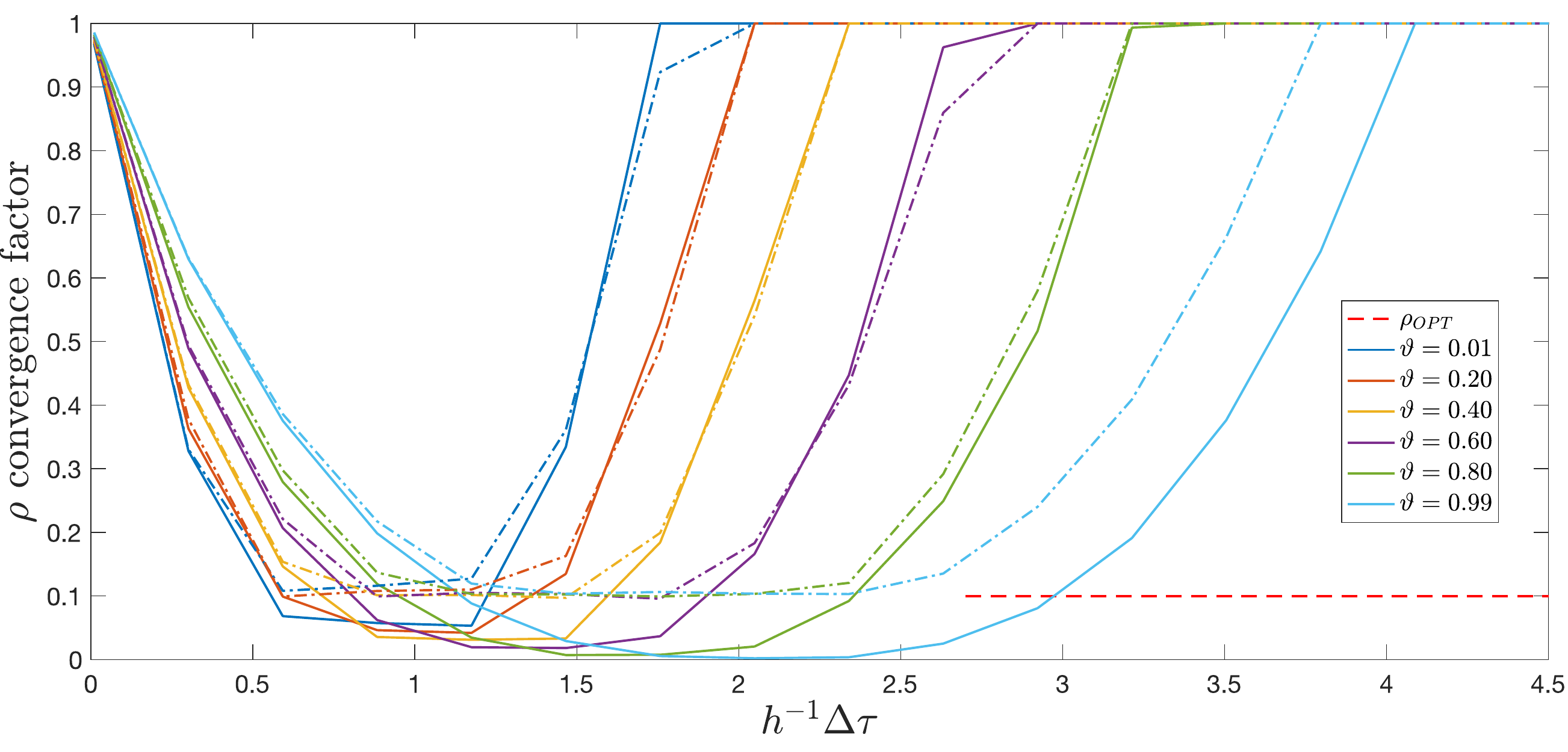}
\end{center}
\end{minipage}
\caption{TGCS convergence factor for 2D Dirichlet (top) and Neumann (bottom) boundary conditions in a rectangular domain (\S\ref{sect:rectangle}). The red dashed line is the optimal convergence factor without ghost nodes.}
\label{fig:2D_rhoD}
\end{figure}

In Fig.~\ref{fig:2D_rhoD} we compare the convergence factors obtained with {\tt REL-ALL} and {\tt REL-BC} for different values of $\vartheta \in [0,1]$.
We see that, for different values of $\vartheta$ (different colors), the convergence factor obtained with {\tt REL-BC} (solid line) performs similarly to {\tt REL-ALL} (dashed line) or better. If it does not converge, the convergence factor is cut artificially at $\rho=1$ for graphical purposes. Approximately, we observe that:
\begin{equation}\label{rhomax}
\rho_\text{REL-ALL} \approx \max \left\{ \rho_\text{REL-BC}, \rho_\text{OPT}  \right\},
\end{equation}
where $\rho_\text{OPT}$ is the convergence factor obtained for regular rectangular domains $\Omega = [-1,1]^2$ where the ghost point technique is not employed. For
$\nu_1=2$, $\nu_2=1$, GS-LEX as a smoother and FW restriction operator we have
\begin{equation}\label{eq:rhoopt2D}
\rho_\text{OPT} \approx 0.119,
\end{equation}
as predicted by the LFA (see \cite[Table 4.1]{Trottemberg:MG}).
Fig.~\ref{fig:2D_rhoD} suggests that for each value of $\vartheta$, there exists a range of possible values of $\Delta \tau$ such that 
$\rho_\text{REL-ALL} \approx \rho_\text{OPT}$. Therefore, it is not necessary to find $\Delta \tau$ that minimizes $\rho_\text{REL-BC}$, it is instead sufficient to identify a value of $\Delta \tau$ in that range. However, it is clear from the figure that $\Delta\tau\approx1.75$ would yield $\rho_\text{REL-ALL}\approx0.3$ for any $\vartheta$ (see Remark~\ref{remark:dtau} for an explanation of these values), but no value of $\Delta\tau$ can yield $\rho_\text{REL-ALL}\approx \rho_\text{OPT}$ for any $\vartheta$.

We want to find $\Delta \tau$ such that the smoothing property of the relaxation scheme is optimal.
To this end, we want to quantify how much the relaxation scheme reduces the high-oscillatory Fourier modes along the tangential direction of the boundary. Let
\begin{equation}\label{FourierModeDef}
\varphi(y) = e^{\iu \alpha y/h}
\end{equation}
be a Fourier mode along the tangential (vertical) direction, with $\alpha \in [0,\pi]$. Observe that we use $\iu$ to denote the imaginary unit $\iu^2=-1$, to be distinguished from $i$, that will be used later to denote some indices.
The Fourier mode $\varphi$ is a highly-oscillatory mode if $\alpha \in [\pi/2,\pi]$.

To perform the BLFA, we work on a semi-infinite domain:
\[
\Omega_\infty = \left\{ (x,y) \colon x \leq x_V\right\}
\]
and we analyze the smoothing of {\tt REL-BC} along the tangential direction.
As explained above, the relaxation {\tt REL-BC} is performed in two steps. The first step consists of solving the equations on the internal grid points in block when the ghost values are set to $\varphi$. Without loss of generality, we can work with the homogeneous problem, i.e.~$f_{i,j}=0$, which means solving the problem
\begin{equation}\label{eq:discpr}
\begin{split}
\frac{4 u_{ij} - (u_{i+1,j}+u_{i-1,j}+u_{i,j+1}+u_{i,j-1})}{h^2} = 0,
&
\text{ for } i = -\infty, \ldots, N-1, \quad j = -\infty, \ldots, +\infty,
\\
\lim_{i \rightarrow -\infty} u_{i,j} = 0,
\quad & 
u_{i,j} = \varphi(y_j) \text{ for } i=N.
\end{split}
\end{equation}
The exact solution of this difference equation is:
\begin{equation}\label{eq:uij}
u_{ij} = e^{p(\alpha) (x_i-1)/h}\varphi(y_j) 
\end{equation}
where 
\begin{equation}\label{alphaNalpha}
p(\alpha) = \cosh^{-1} (2 - \cos(\alpha)).
\end{equation}

The second step consists of relaxing the ghost values by \eqref{relax2DV} with $g_j = 0$:
\begin{equation}\label{eq:relD}
u_{N,j}^{(1)} = u_{N,j} - \Delta \tau \, \left( \frac{\vartheta (\vartheta-1)}{2} \, u_{N-2,j} + \vartheta (2-\vartheta) \, u_{N-1,j}+ \frac{(1-\vartheta)(2-\vartheta)}{2} u_{N,j} \right)
\end{equation}
\[
\Longrightarrow u_{N,j}^{(1)} = \varphi^{(1)}(y_j) = \varphi(y_j) - \Delta \tau \, \left( \frac{\vartheta (\vartheta-1)}{2} \, e^{-2 p(\alpha)}\varphi(y_j) + \vartheta (2-\vartheta) \, e^{-p(\alpha)}\varphi(y_j)+ \frac{(1-\vartheta)(2-\vartheta)}{2} \varphi(y_j) \right)
\]
\begin{equation}\label{LFA2Dampl}
\Longrightarrow \frac{\varphi^{(1)}(y_j) }{ \varphi(y_j)} = 1 - \Delta \tau \, \left( \frac{\vartheta (\vartheta-1)}{2} \, e^{-2 p(\alpha)} + \vartheta (2-\vartheta) \, e^{-p(\alpha)}+ \frac{(1-\vartheta)(2-\vartheta)}{2} \right).
\end{equation}
The right-hand side of \eqref{LFA2Dampl} is the amplification factor of the tangential Fourier mode $\varphi = e^{\iu \alpha y/h}$ for a given $\vartheta \in [0,1)$ and a given $\Delta \tau > 0$, that we denote as:
\begin{equation}\label{eq:G}
G(\alpha, \vartheta, \Delta \tau) =  1 - \Delta \tau \; G_0(\alpha,\vartheta)
\text{ with }
G_0(\alpha,\vartheta) = 
\sum_{r=0}^2 c_{-r} e^{-r \, p(\alpha)}
\end{equation}
where the coefficients $c_{-k}$ are given by \eqref{BC1D} and are
\begin{equation}\label{eq:coeff}
c_{-2}(\vartheta) = \frac{\vartheta (\vartheta-1)}{2}, \quad c_{-1}(\vartheta) = \vartheta (2-\vartheta), \quad c_{0}(\vartheta) = \frac{(1-\vartheta) (2-\vartheta)}{2}.
\end{equation}

Finding the optimal $\Delta \tau$ for a given $\vartheta$ requires to minimize the largest amplification factor for the high-oscillatory modes $\alpha \in [\pi/2,\pi]$ and thus
\begin{equation}\label{eq:def_dtau_opt}
\Delta \tau_\text{OPT} (\vartheta) = \arg \min_{\Delta \tau >0} \sup_{\pi/2 < \alpha < \pi} |G(\alpha, \vartheta, \Delta \tau)|.
\end{equation}
\begin{lemma}\label{lemma:fnostat}
Let $G_0(\alpha, \cdot) \colon (\pi/2,\pi ) \to \mathbb{R}$ be the continuous function of $\alpha$ defined in~\eqref{eq:G} with coefficients \eqref{eq:coeff}.
Then $G_0$ does not have stationary points in $(\pi/2,\pi)$. 
\end{lemma}
\noindent {\it Proof.} 
We have:
\[
\frac{\partial G_0}{\partial \alpha} = - e^{- p(\alpha)} \left( 2c_{-2} e^{-p(\alpha)} + c_{-1} \right) \cdot p'(\alpha). 
\]
Since $p'(\alpha) =\sin(\alpha) \left( (2-\cos(\alpha))^2-1 \right)^{-1/2} \neq 0$ for any $\alpha \in (0,\pi)$, then the only possible stationary points of $G_0(\alpha, \cdot)$ satisfy the equation
\begin{equation}\label{eq:alpha}
e^{-p(\alpha)} = \frac{c_{-1}}{-2c_{-2}} = 1 + \frac{1}{1-\vartheta}.
\end{equation}
Since $1 + (1-\vartheta)^{-1} > 1$ for $\vartheta\in[0,1)$ and $e^{-p(\alpha)}<1$ (because $p(\alpha)>0$ for every $\alpha$), then~\eqref{eq:alpha} is never satisfied, proving that $G_0(\alpha,\cdot)$ has no stationary points in $(\pi/2,\pi)$.
$\square$
%

\begin{lemma}\label{lemma:G0pos}
Let $G_0(\cdot, \vartheta) \colon [0,1) \to \mathbb{R}$ be the continuous function of $\vartheta$ defined in~\eqref{eq:G} with coefficients \eqref{eq:coeff}. Then $G_0(\cdot, \vartheta)>0$ for every $\vartheta \in [0,1)$.
\end{lemma}

\noindent {\it Proof.}
We have $G_0(\cdot, 0) = 1>0$ and $G_0(\cdot, 1) = e^{-p(\alpha)}>0$. To complete the proof, it is sufficient to show that $G_0(\cdot, \vartheta)$ has no stationary points in $(0,1)$.
We have
\[
\frac{\partial G_0}{\partial \vartheta} = \frac{1}{2} (e^{-p(\alpha)}-1) \left(  (e^{-p(\alpha)}-1) \vartheta - (e^{-p(\alpha)}-3) \right)
\]
Since $p(\alpha)>0$, we have
\[
\frac{\partial G_0}{\partial \vartheta} = 0 \Longleftrightarrow \vartheta = \frac{e^{-p(\alpha)}-3}{e^{-p(\alpha)}-1}.
\]
We observe that $(e^{-p(\alpha)}-3)(e^{-p(\alpha)}-1)^{-1} > 1$ (because $p(\alpha)>0$), proving that there are no stationary point of $G_0(\cdot, \vartheta)$ in $(0,1)$.
$\square$

\begin{lemma}\label{lemma:lines}
We consider linear functions $y_i(x) = 1-m_i x$ for $i=1, \ldots, l$, with $m_1 \geq m_2 \geq  \ldots \geq m_l \geq 0$. 
Then:
\[
\arg \min_{x>0} \max_{i} \left\{ |y_i(x)| \right\} = \frac{2}{m_1+m_l}.
\]
\end{lemma}

\begin{figure}
\begin{center}
\begin{tikzpicture}[scale=1, every node/.style={scale=1.0}]
\draw[thick,->] (0,-0.5) -- (0,6);
\draw[thick,->] (-0.5,0) -- (10,0);
\draw[-] (0,3) -- (1,0) -- (1+1*1.8,3*1.8);
\draw[-] (0,3) -- (2,0) -- (2+2*1.5,3*1.5);
\draw[-] (0,3) -- (6,0) -- (6+6*0.5,3*0.5);
\draw (6,2) node {$\ldots\; \ldots \; \ldots$};
\node[rotate=71.57, above] at (1+1*1.8,3*1.8) {$|y_1(x)| = |1-m_1 x| $};
\node[rotate=56.31, above] at (2+2*1.5,3*1.5) {$|y_2(x)| = |1-m_2 x| $};
\node[rotate=26.57, above] at (6+6*0.5,3*0.5) {$|y_l(x)| = |1-m_l x| $};
\draw[very thick,-] (0,3) -- (0.5714*3,2.1429) -- (0.5714*3 + 0.7*1,2.1429 + 0.7*3);
\node[rotate=0] at (1.1,3.5) {$\tilde{y}(x)$};
\draw[thick,dashed] (0.5714*3,0) -- (0.5714*3,2.1429);
\node[rotate=0] at (0.5714*3,-0.3) {$x^*$};
\draw[->] (1.0,3.2) -- (0.8,2.7);
\draw[->] (1.5,3.4) -- (1.95,3.1);
\node at (-0.3,6.0) {$y$};
\node at (10,-0.3) {$x$};
\end{tikzpicture}
\end{center}
\caption{Illustration of the minimax problem for linear functions (Lemma \ref{lemma:lines}).}
\label{fig:lemma_lines}
\end{figure}
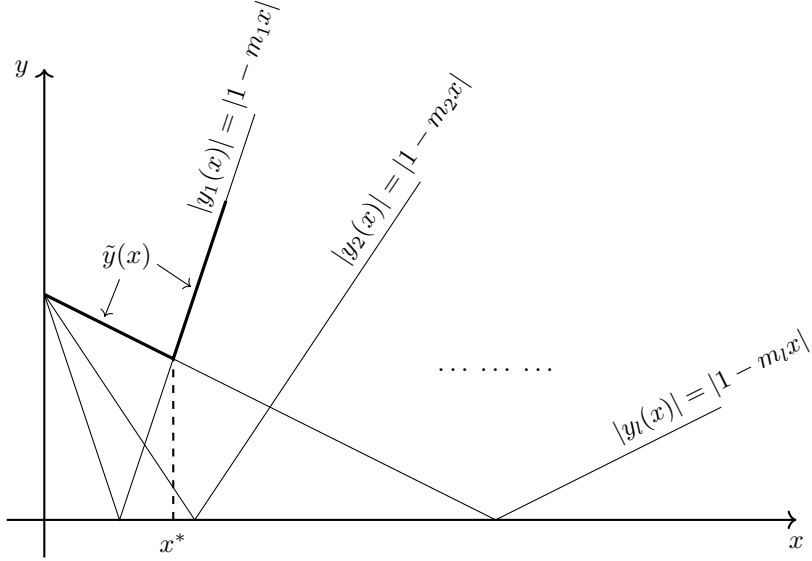

{\it Proof.}
Let $\tilde{y}(x) = \max_{i} \left\{ |y_i(x)| \right\}$.
It can be easily proved (see Fig.~\ref{fig:lemma_lines}) that
\[
\tilde{y}(x) = 
\left\{
\begin{matrix}
y_l(x) & \text{ if } x \leq x^* \\
-y_1(x) & \text{ if } x > x^* \\
\end{matrix}
\right. . 
\]
where $x^*$ is the intersection between the lines $y_1(x)$ and $-y_l(x)$, namely $x^* = 2 (m_1+m_l)^{-1}$.
Observe that $\tilde{y}(x)$ is decreasing for $x<x^*$ and increasing for $x>x^*$, proving that $\tilde{y}(x^*) \leq \tilde{y}(x)$ for every $x>0$.
$\square$

\begin{proposition}\label{prop:opt}
The optimal $\Delta \tau$ defined in~\eqref{eq:def_dtau_opt} is obtained by the following formula:
\begin{equation}\label{dtauoptLFA}
\Delta \tau_\text{OPT} (\vartheta) = \frac{2}{G_0(\pi/2, \vartheta)+G_0(\pi, \vartheta)}
\end{equation}
where $G_0$ is defined in~\eqref{eq:G} with coefficients \eqref{eq:coeff}.
\end{proposition}

{\it Proof.}
From Lemma \ref{lemma:fnostat} we have that $G(\alpha,\cdot, \cdot) = 1 - \Delta \tau G_0(\alpha, \cdot)$ has no stationary points for $\pi/2 < \alpha < \pi$. Then:
\[
\sup_{\pi/2 < \alpha < \pi} |G(\alpha, \vartheta, \Delta \tau)| = \max \left\{ |G(\pi/2, \vartheta, \Delta \tau)|, |G(\pi, \vartheta, \Delta \tau)| \right\}.
\]
Moreover, $G(\pi,\cdot, \Delta \tau)$ and $G(\pi/2,\cdot, \Delta \tau)$ are linear functions of $\Delta \tau$ and from Lemma \ref{lemma:G0pos} they have negative slopes with magnitude $G_0(\pi/2, \cdot)$ and $G_0(\pi, \cdot)$, respectively. Therefore, from Lemma \ref{lemma:lines}, we obtain \eqref{dtauoptLFA}.
$\square$

\begin{figure}
\begin{center}
\begin{tabular}{cc}
(Dirichlet) & (Neumann) \\
\includegraphics[width=0.45\textwidth]{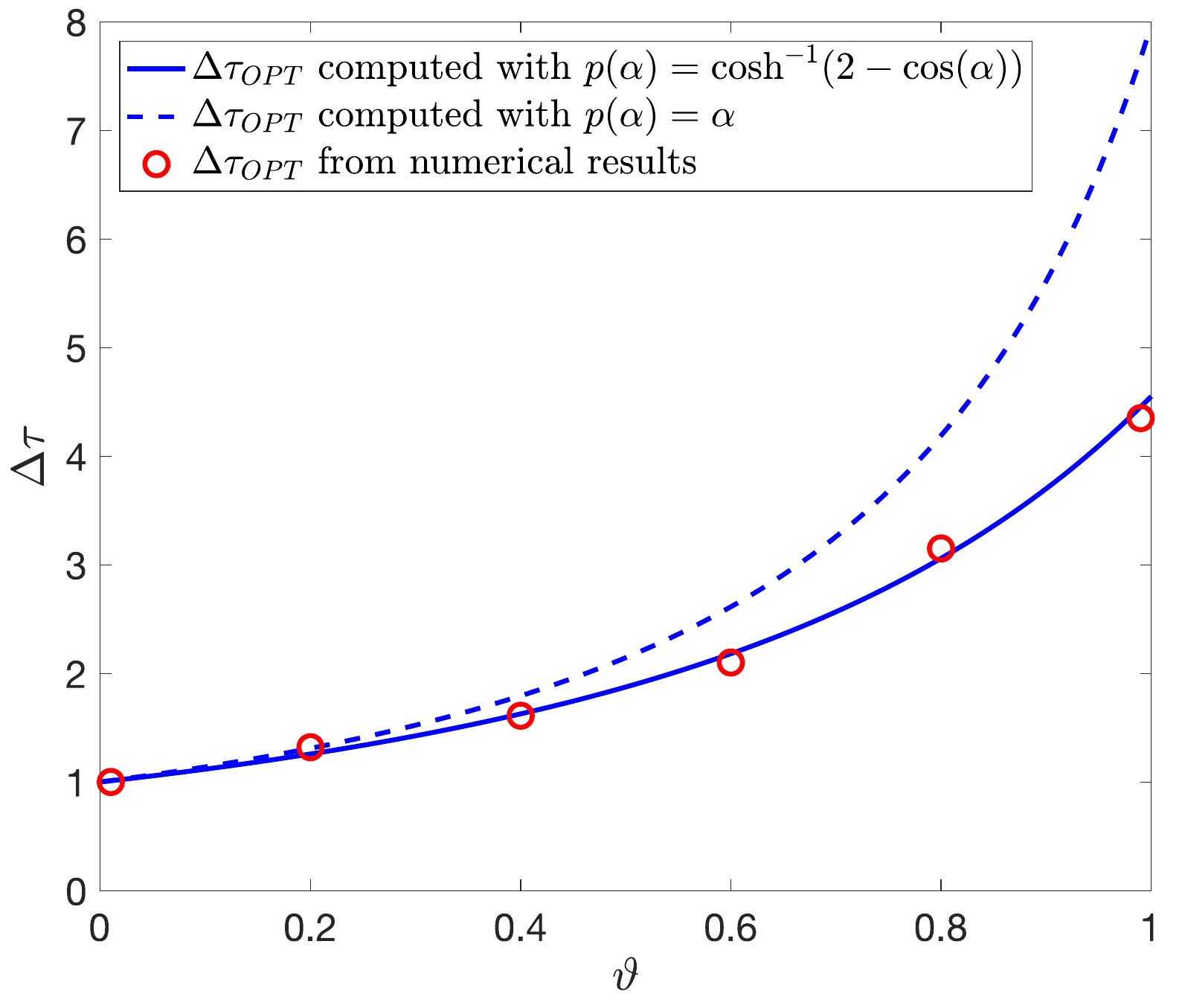}
&
\includegraphics[width=0.45\textwidth]{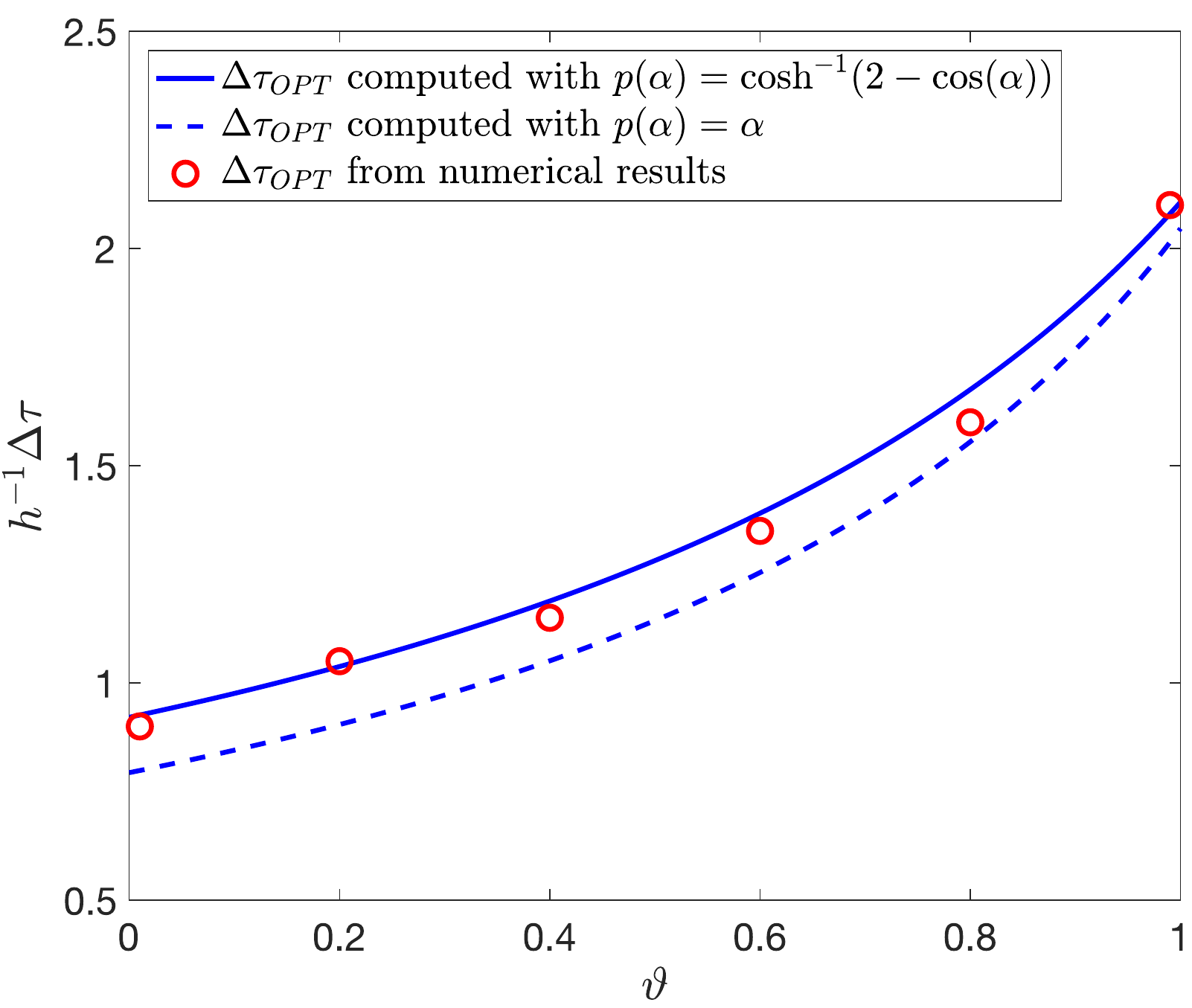}
\end{tabular}
\end{center}
\caption{Optimal relaxation parameter $\Delta \tau$ for the rectangular domain of \S\ref{sect:rectangle} from numerical experiments (red points, inferred from Fig.~\ref{fig:2D_rhoD}) and from the BLFA, namely Eqs.~\eqref{dtauoptLFA} and \eqref{dtauoptLFA_Neumann} for Dirichlet (left plot) and Neumann (right plot) boundary conditions, respectively. The dashed lines are the theoretical $\Delta \tau$ if the solutions to the continuous problems were used instead of the discretized ones (see Remark \ref{solnum2D}).
}
\label{fig:dtau_opt}
\end{figure}

In Fig.~\ref{fig:dtau_opt} we observe a rather satisfactory match between the optimal $\Delta \tau$ obtained from numerical results and from the BLFA \eqref{dtauoptLFA}.

\begin{remark}\label{solnum2D}
Instead of using the numerical solution \eqref{eq:uij} of \eqref{eq:discpr}, one can think to use the exact solution of the continuous problem
\[
- \Delta u = 0 \text{ in } \Omega_\infty, \quad \lim_{x \to \infty} u = 0, \quad u(x_V,y) =  \varphi(y),
\]
namely \eqref{eq:uij} with $p(\alpha)=\alpha$, that we denote by $u^\text{EXA}_{ij}$. We observe that this choice will lead to a completely wrong estimate of $\Delta \tau_\text{OPT}$. In fact, one can say that numerical and exact solution would differ by $u^\text{EXA}_{ij} - u_{ij} = \mathcal{O}(h^2)$, but this is not true for a fixed value of $\alpha$, since the boundary condition scale with $h$, namely \eqref{FourierModeDef} depends on $h$. This means that $u^\text{EXA}_{ij} - u_{ij} = \mathcal{O}(1)$ for $h \to 0$.
The statement $u^\text{EXA}_{ij} - u_{ij} = \mathcal{O}(h^2)$ would only be true if
the boundary condition did not scale with $h$, meaning that (from \eqref{FourierModeDef}) $\alpha/h = \text{const.}$, that is not the case of our interest for the BLFA. The discrepancy between the optimal $\Delta \tau_\text{OPT}$ obtained with \eqref{alphaNalpha} and the one obtained with $p(\alpha) = \alpha$ is observed numerically in Fig.~\ref{fig:dtau_opt} for the 2D case and Fig.~\ref{fig:dtau_opt3D} for the 3D case.
\end{remark}

\begin{figure}
\begin{center}
\begin{minipage}{0.3\textwidth}
\includegraphics[width=0.99\textwidth]{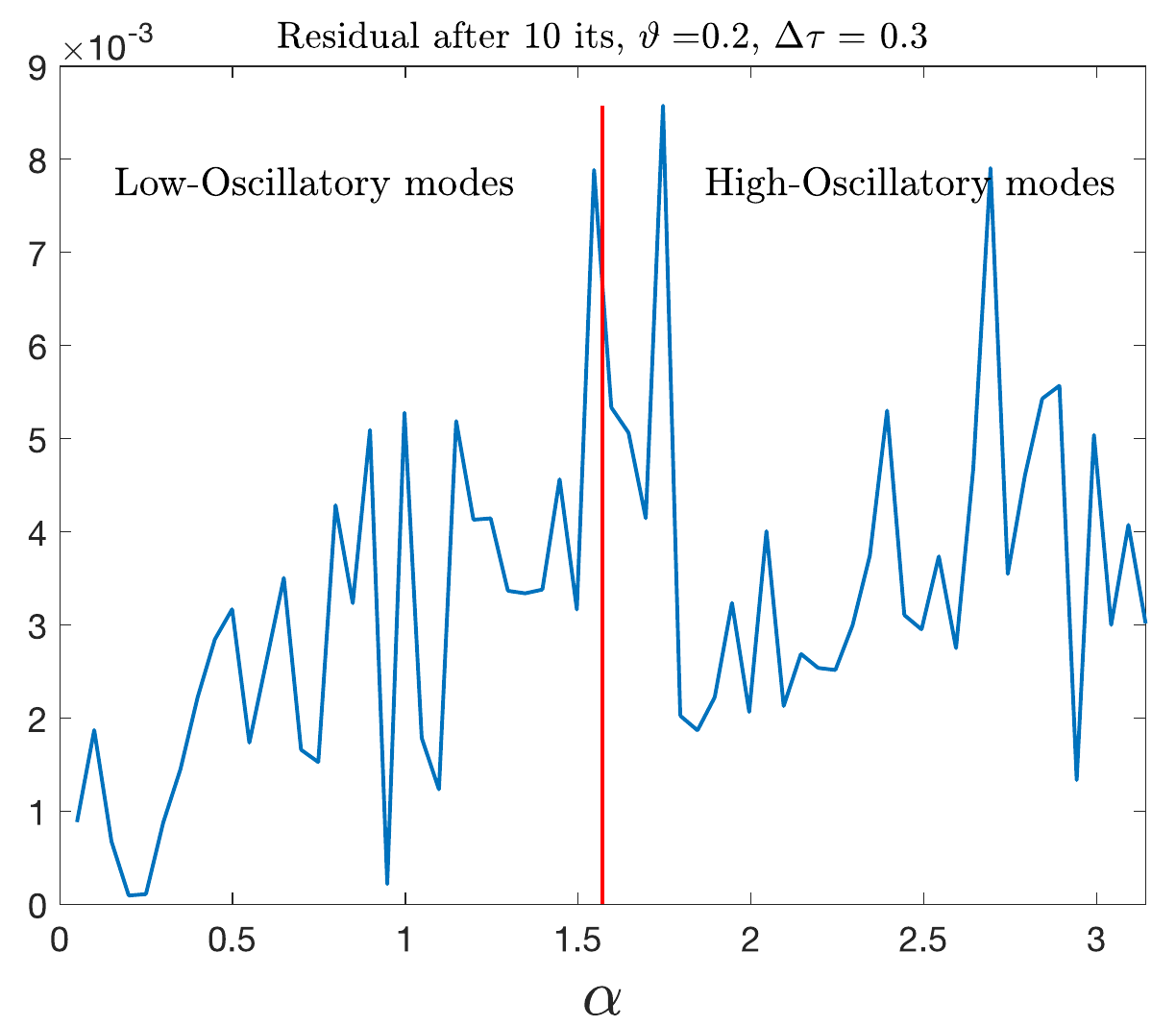}
\end{minipage}
\begin{minipage}{0.3\textwidth}
\includegraphics[width=0.99\textwidth]{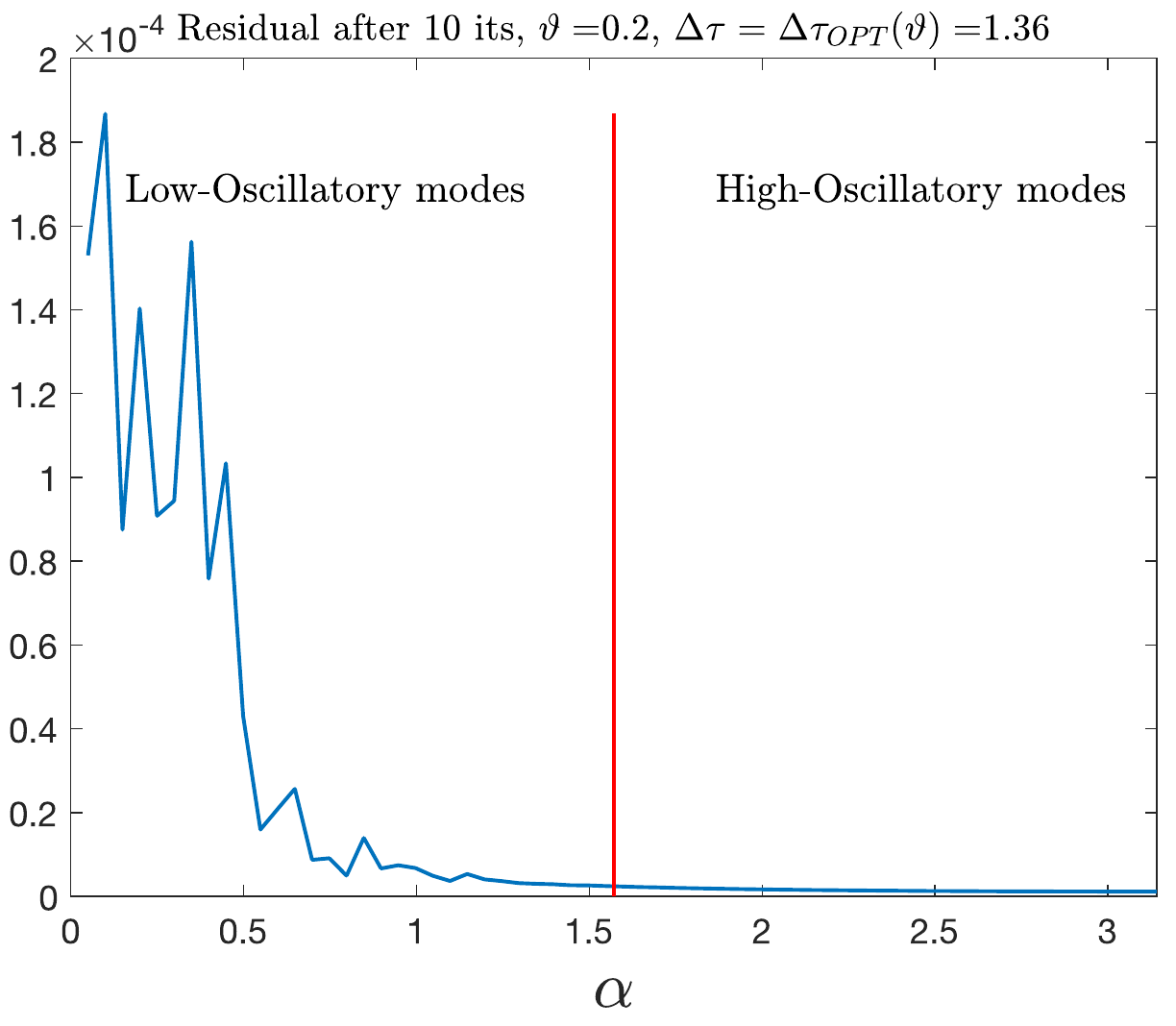}
\end{minipage}
\begin{minipage}{0.3\textwidth}
\includegraphics[width=0.99\textwidth]{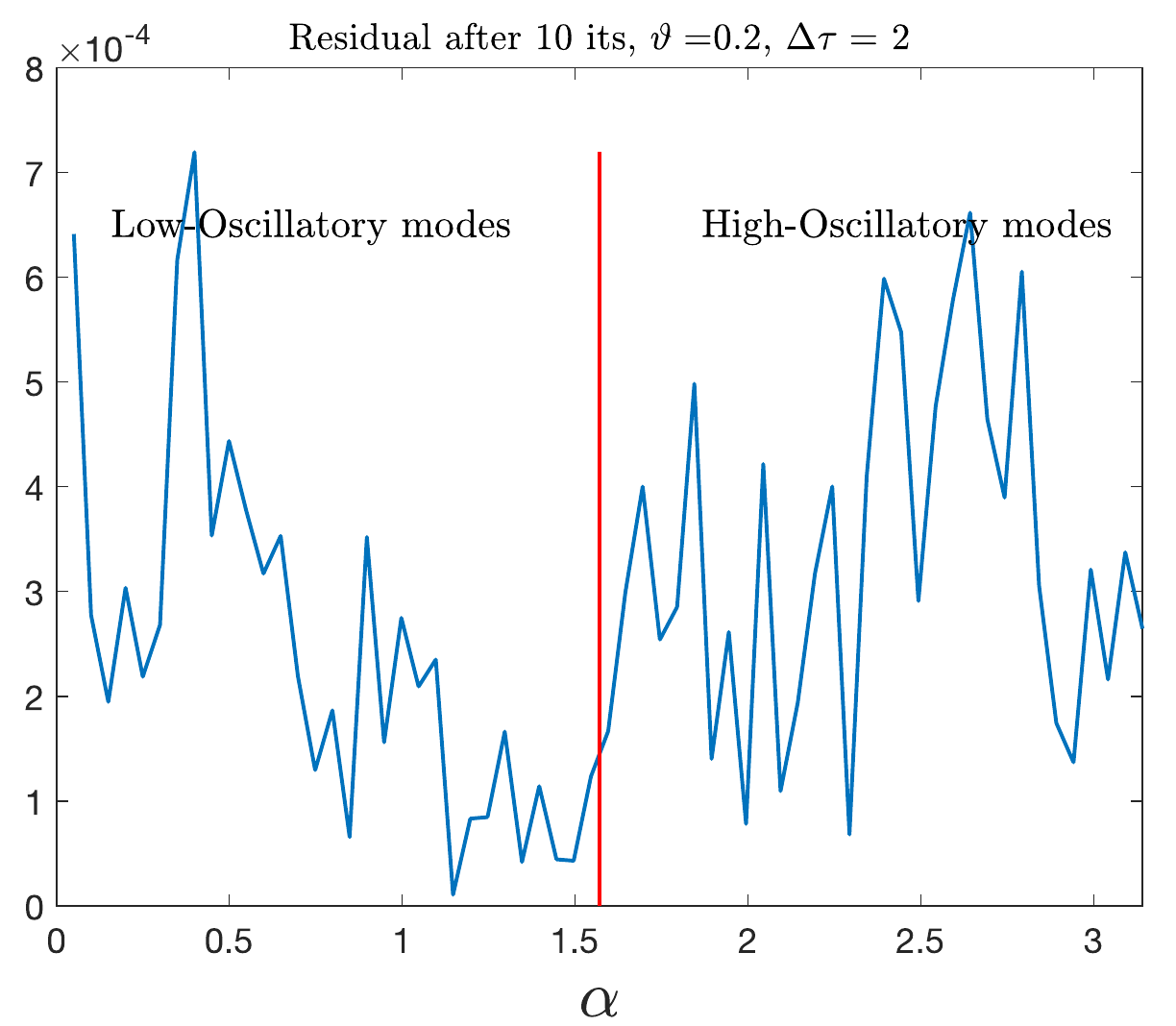}
\end{minipage}
\begin{minipage}{0.3\textwidth}
\includegraphics[width=0.99\textwidth]{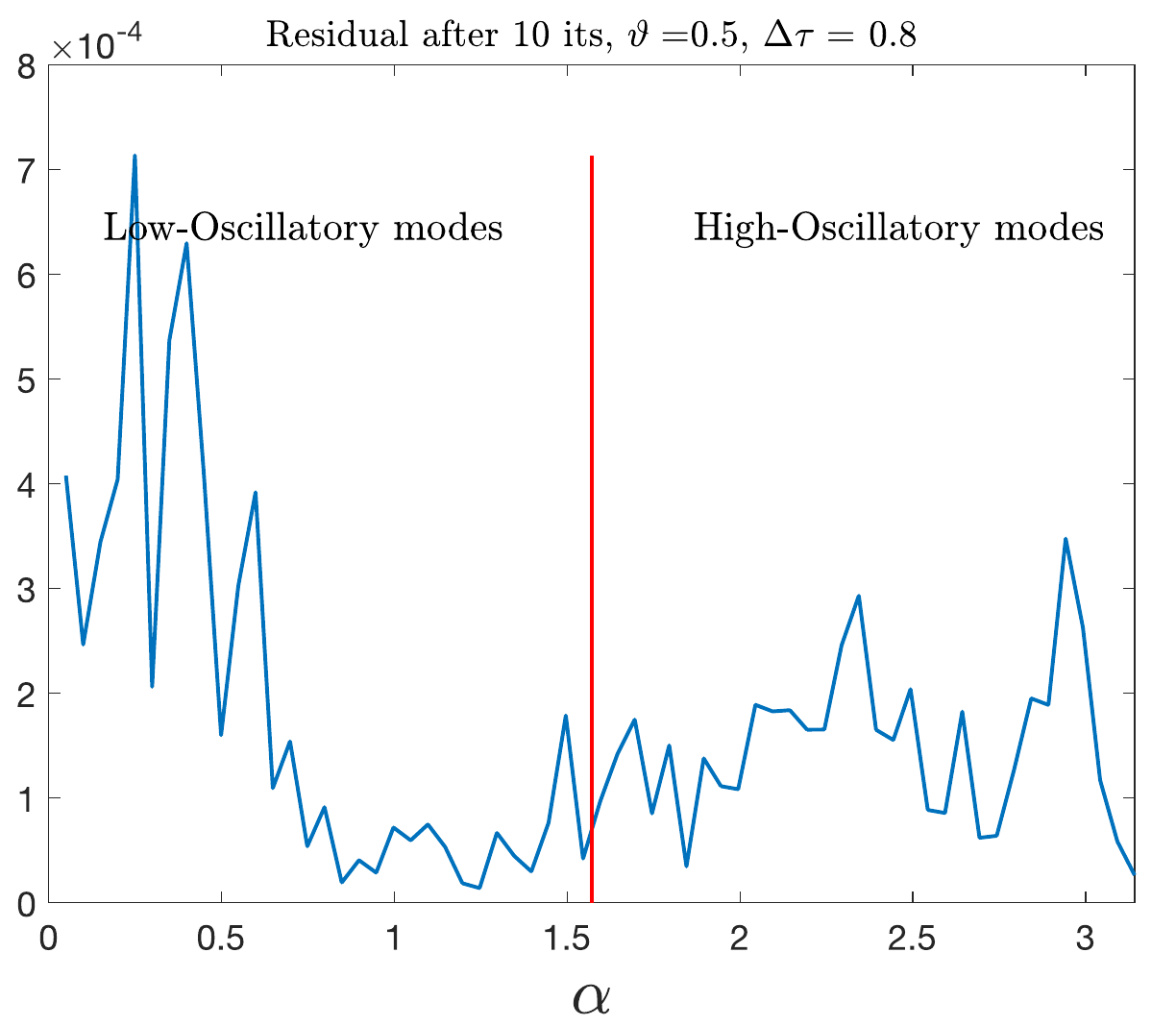}
\end{minipage}
\begin{minipage}{0.3\textwidth}
\includegraphics[width=0.99\textwidth]{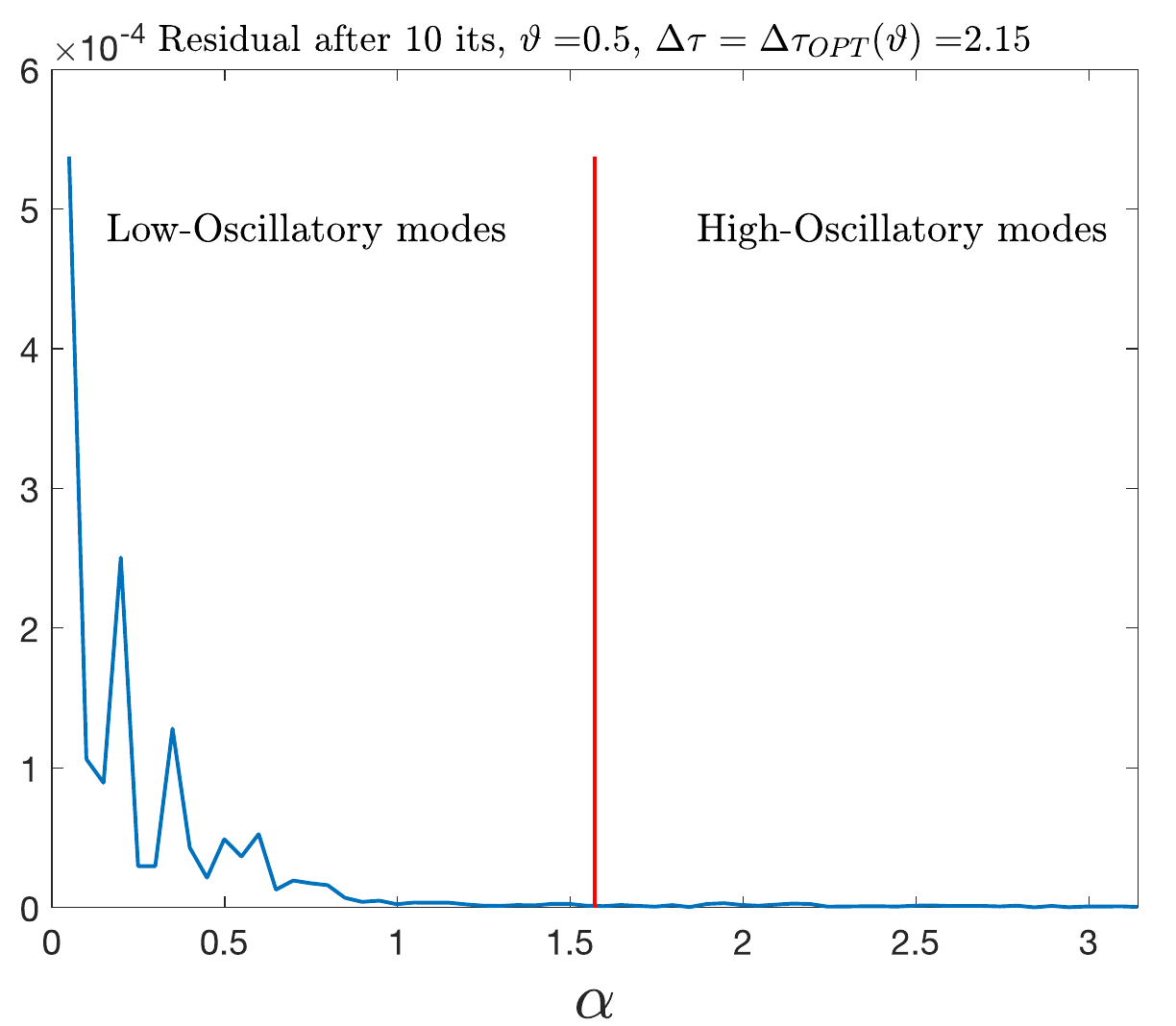}
\end{minipage}
\begin{minipage}{0.3\textwidth}
\includegraphics[width=0.99\textwidth]{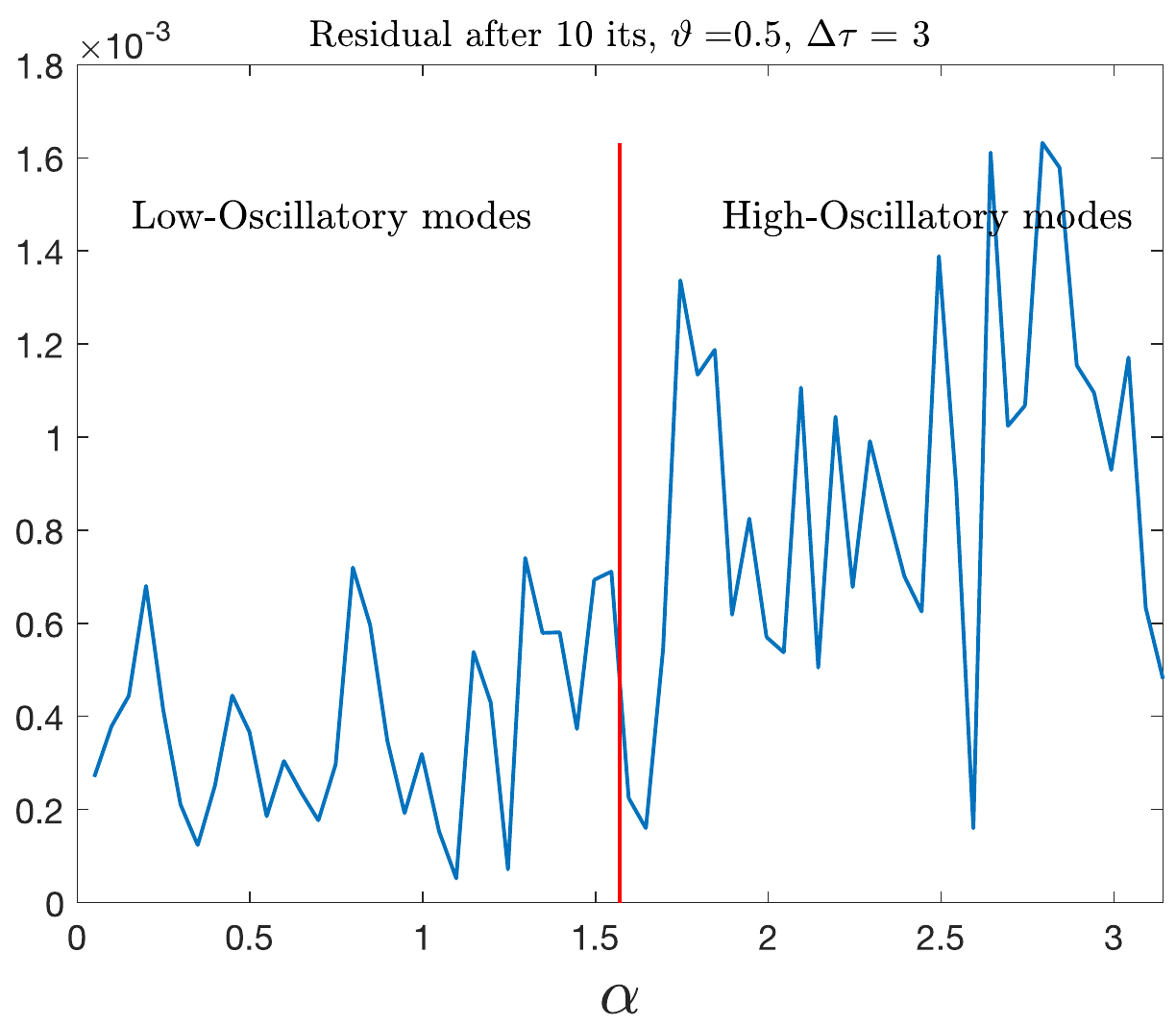}
\end{minipage}
\begin{minipage}{0.3\textwidth}
\includegraphics[width=0.99\textwidth]{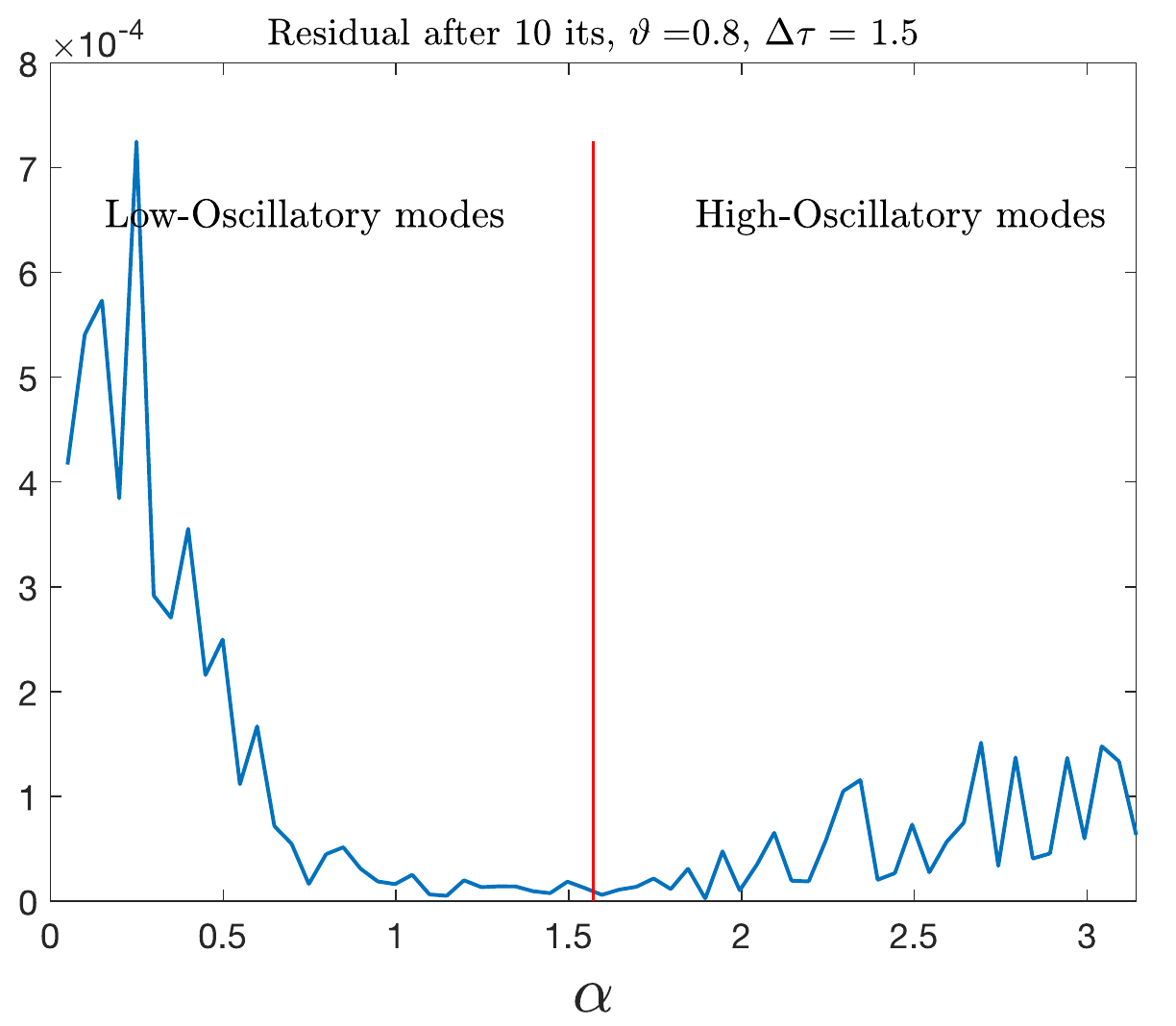}
\end{minipage}
\begin{minipage}{0.3\textwidth}
\includegraphics[width=0.99\textwidth]{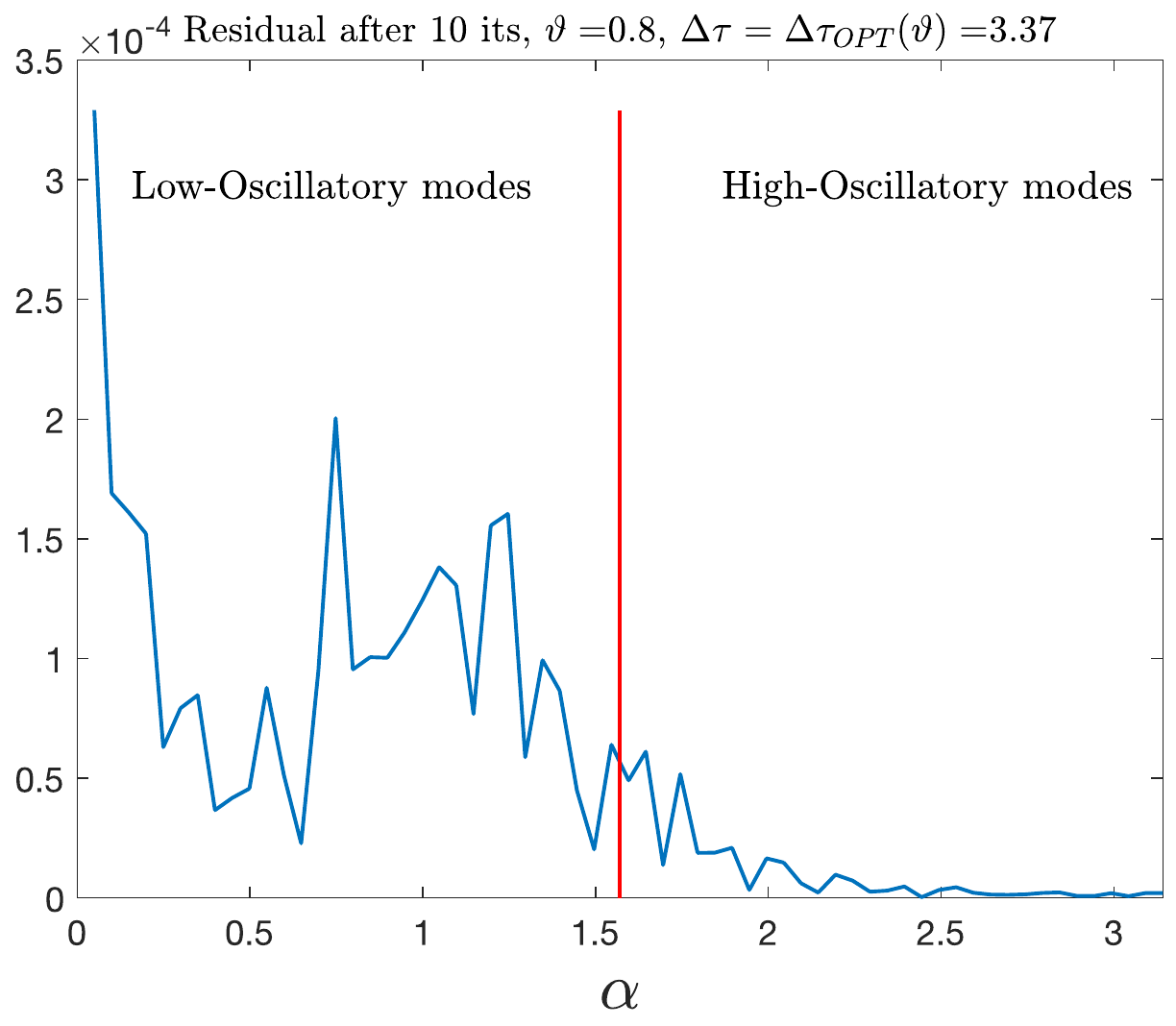}
\end{minipage}
\begin{minipage}{0.3\textwidth}
\includegraphics[width=0.99\textwidth]{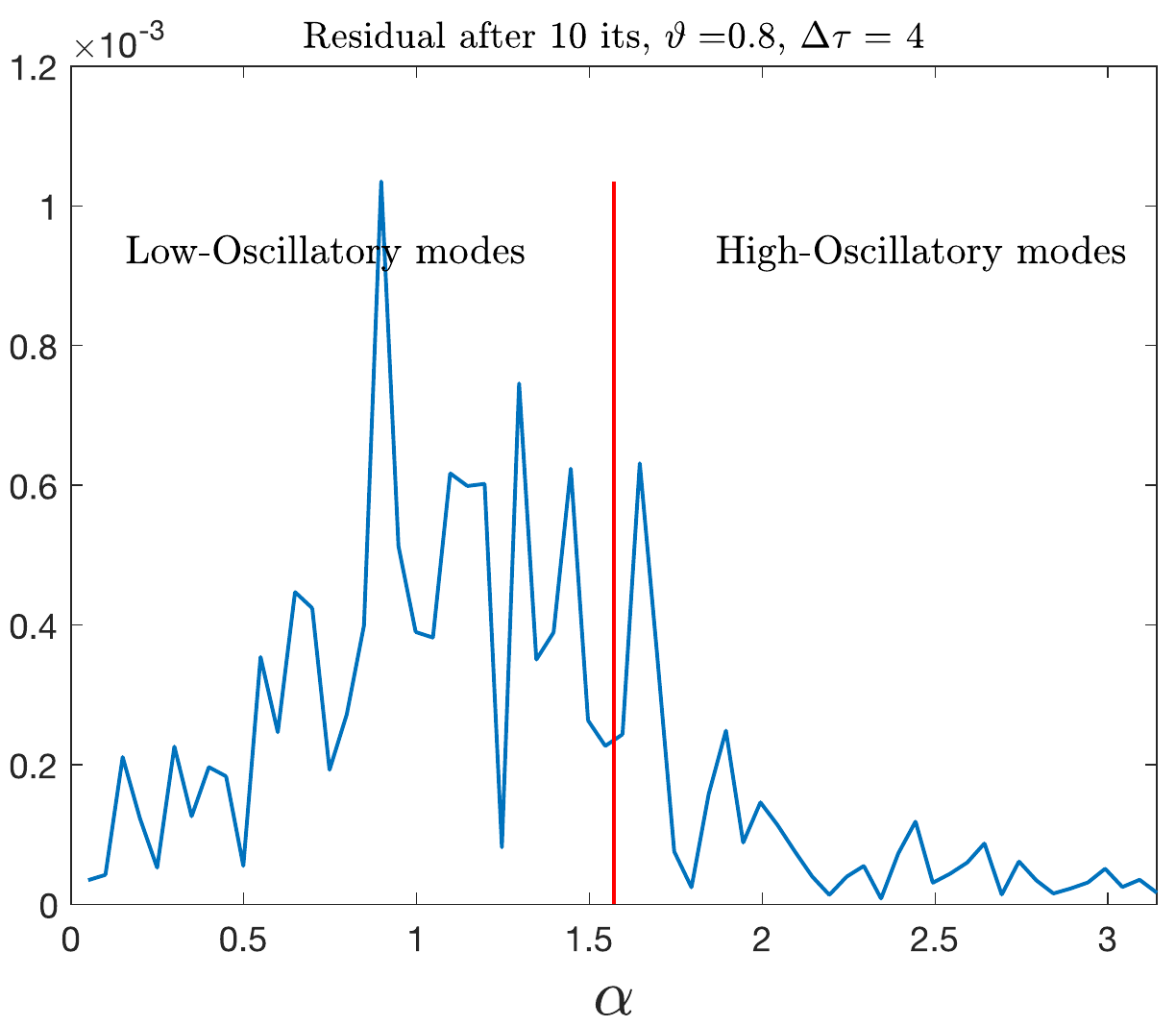}
\end{minipage}
\end{center}
\caption{
Frequency analysis of the residual along the ghost points. A vertical line at $\alpha = \pi/2 $ separates the low frequency components (left) and the high frequency components (right). 
Tests for $\vartheta = 0.2$ (first row), $\vartheta = 0.5$ (second row) and $\vartheta = 0.8$ (third row). 
Values of $\Delta \tau$: the optimal value (middle column), a lower value (left column) and a higher value (right column).
}
\label{fig:FFT_VLINE_DTAU}
\end{figure}

In order to show that high frequency components of the error are dumped much quicker when $\Delta \tau = \Delta \tau_\text{OPT}$, we compute the Discrete Fourier Transform of the residual along the ghost points (namely along the tangential direction of the boundary) and plot the components of the residual versus the frequency of the Fourier modes \eqref{FourierModeDef}. 
In detail, let $r_\text{BDY} \in \mathbb{R}^{N-1}$ be the residual along the ghost points defined as:
\[
(r_\text{BDY})_j = \vec{r}(x_N,y_j) \text{ for } j=1,\ldots, N-1. 
\]
Let $\vec{y} \in \mathbb{R}^{N-1}$ be the $y-$coordinate vector with components $y_j = \frac{j}{N}$, for $j=1,\ldots,N-1$.
The Fourier modes \eqref{FourierModeDef} are discretized as
\[
\varphi^{(\alpha_k)} (\vec{y}) = e^{\iu \alpha_k \vec{y}/h}, \; \text{ with } \alpha_k = -\pi+k\frac{2 \pi}{N-1} \in (-\pi,\pi), \text{ for } k=1,\ldots,N-1.
\]
The residual $r_B$ can be decomposed as a linear combination of Fourier modes 
\[
r_\text{BDY} = \sum_{k=1}^{N-1} R_k \, \varphi^{(\alpha_k)} (\vec{y})
\]
where the components $R_k \in \mathbb{C}$ are given by the Discrete Fourier Transform formula:
\[
R_k = \frac{1}{N-1} \sum_{j=1}^{N-1} (r_B)_j \, e^{-\iu \alpha_k y_j/h}, \text{ for } k=1,\ldots,N-1.
\]
Since $r_\text{BDY}$ is a real vector, then $\overline{R_{k_1}} = R_{k_2}$ if $\alpha_{k_1} = - \alpha_{k_2}$, where $\overline{R_{k_1}}$ denotes the complex conjugate of $R_{k_1}$.
Therefore, we can restrict our study to $\alpha_k \in (0,\pi)$. 
The Fourier modes $\varphi^{(\alpha_k)} (\vec{y})$ are divided in low-frequency modes ($\alpha_k \in (0,\pi/2))$ and high-frequency modes ($\alpha_k \in (\pi/2,\pi))$. 
For each value of $\alpha_k \in (0,\pi)$ we plot $2 \left| R_k \right|$ so that we include the contribution of the symmetric mode in $(-\pi,0)$.
Fig.~\ref{fig:FFT_VLINE_DTAU} shows the amplitudes $2 \left| R_k \right|$ versus $\alpha_k \in (0,\pi)$ where $r_\text{BDY}$ is the residual on the ghost points after 10 relaxations on the fine grid, starting from a random initial guess whose values are randomly uniformly distributed between $-1$ and $1$.
A vertical line at $\alpha = \pi/2 $ separates the low frequency components (left) from the high frequency components (right). We repeat the test for $\vartheta = 0.2$ (first row), $\vartheta = 0.5$ (second row) and $\vartheta = 0.8$ (third row). For each value of $\vartheta$ we choose three values of $\Delta \tau$: the optimal value (middle column), a lower value (left column) and a higher value (right column). We observe that the smoothing property is clearly satisfied only when $\Delta \tau = \Delta \tau_\text{OPT}$, namely the ratio between high- and low-frequency component amplitudes is smaller. The plots are obtained for $N=128$.

%

\subsubsection{Neumann boundary conditions: BLFA}\label{sect:2D_neuLFA}
We use the same approach as for Dirichlet boundary conditions. We replace the relaxation \eqref{eq:relD} with the Neumann boundary relaxation:
\[
u_{N,j}^{(1)} = u_{N,j} - \frac{\Delta \tau}{h} \, \left(  \left(\frac{1}{2}-\vartheta \right) \, u_{N-2,j} + 2\left( \vartheta-1 \right) \, u_{N-1,j}+ 
\left( \frac{3}{2}-\vartheta \right) u_{N,j} \right)
\]
\[
\Longrightarrow u_{N,j}^{(1)} = \varphi^{(1)}(y_j) = \varphi(y_j) - \frac{\Delta \tau}{h} \, \left(  \left(\frac{1}{2}-\vartheta \right) \, e^{-2 p(\alpha)} \varphi(y_j) + 2\left( \vartheta-1 \right) \, e^{-p(\alpha)} \varphi(y_j)+ \left( \frac{3}{2}-\vartheta \right) \varphi(y_j)  \right)
\]
\begin{equation}\label{LFA2DamplNeumann}
\Longrightarrow \frac{\varphi^{(1)}(y_j) }{ \varphi(y_j)} = 1 - \frac{\Delta \tau}{h} \, \left(  \left(\frac{1}{2}-\vartheta \right) \, e^{-2 p(\alpha)} + 2\left( \vartheta-1 \right) \, e^{-p(\alpha)} + \left( \frac{3}{2}-\vartheta \right) \right).
\end{equation}

The right-hand side of \eqref{LFA2DamplNeumann} is the amplification factor along the tangential direction of the Fourier mode $\varphi = e^{i \alpha y/h}$ for the Neumann boundary condition. Here
\begin{equation}\label{eq:GN}
G(\alpha, \vartheta, \Delta \tau) =  1 - \frac{\Delta \tau}{h} \; G_0(\alpha,\vartheta)
\text{ with }
G_0(\alpha,\vartheta) = 
\sum_{r=0}^2 c_{-r} e^{-r \, p(\alpha)}
\end{equation}
where
\begin{equation}\label{eq:coeffN}
c_{-2}(\vartheta) = \frac{1}{2}-\vartheta, \quad c_{-1}(\vartheta) = 2\left( \vartheta-1 \right), \quad c_{0}(\vartheta) = \frac{3}{2}-\vartheta.
\end{equation}
The definition of the optimal $\Delta \tau$ is the same as \eqref{eq:def_dtau_opt}.

\begin{lemma}\label{lemma:fnostatN}
Let $G_0(\alpha, \cdot) \colon (\pi/2,\pi ) \to \mathbb{R}$ be the continuous function of $\alpha$ defined in~\eqref{eq:GN} with coefficients \eqref{eq:coeffN}.
Then $G_0$ does not have stationary points in $(\pi/2,\pi)$. 
\end{lemma}

\noindent {\it Proof.} 
Similarly to Lemma \ref{lemma:fnostat}, we have:
\[
 \frac{\partial G_0}{\partial \alpha} = 0 \Longleftrightarrow 2c_{-2} e^{-p(\alpha)} + c_{-1} = 0. 
\]
If $\vartheta=1/2$, then $c_{-2} = 0$ and $c_{-1}=-1$, implying that $\displaystyle \frac{\partial G_0}{\partial \alpha} \neq 0$. If $\vartheta \neq 1/2$, then 
\[
 \frac{\partial G_0}{\partial \alpha} = 0 \Longleftrightarrow e^{-p(\alpha)} = \frac{c_{-1}}{-2c_{-2}} = 
 1 +  \displaystyle \frac{1}{1-2\vartheta}.
\]
We observe that $1 +  \displaystyle \frac{1}{1-2\vartheta} \notin (0,1)$ and $e^{-p(\alpha)} \in (0,1)$ (because $p(\alpha)>0$ for every $\alpha$), then  $G_0(\alpha,\cdot)$ has no stationary points in $(\pi/2,\pi)$.
$\square$

\begin{lemma}\label{lemma:G0posN}
Let $G_0(\cdot, \vartheta) \colon (0,1) \to \mathbb{R}$ be the continuous function of $\vartheta$ defined in~\eqref{eq:GN} with coefficients~\eqref{eq:coeffN}. Then $G_0(\cdot, \vartheta)>0$ for every $\vartheta \in [0,1)$.
\end{lemma}

\noindent {\it Proof.}
We have 
\[
G_0(\alpha,0) = 
\frac{1}{2} \, e^{-2 p(\alpha)} - 2 \, e^{-p(\alpha)} + \frac{3}{2} =
\frac{1}{2} \left( e^{-2 p(\alpha)} - 4 \, e^{-p(\alpha)} + 3 \right) =
\frac{1}{2} \left( e^{-p(\alpha)} - 1 \right) \left( e^{-p(\alpha)} - 3 \right) > 0
\]
and
\[
G_0(\alpha,1) = 
\frac{1}{2} \left(1-e^{-2 p(\alpha)} \right) > 0.
\]
To complete the proof, it is sufficient to show that $G_0(\cdot, \vartheta)$ does not have stationary points in $(0,1)$.
We have
\[
\frac{\partial G_0}{\partial \vartheta} = -e^{-2p(\alpha)} + 2 e^{-p(\alpha)} - 1 = - (e^{-p(\alpha)}-1)^2 < 0
\]
since $p(\alpha)>0$. Therefore, there are no stationary point of $G_0(\cdot, \vartheta)$ in $(0,1)$.
$\square$

\begin{proposition}\label{prop:optN}
The optimal $\Delta \tau$ for Neumann boundary conditions is obtained by the following formula:
\begin{equation}\label{dtauoptLFA_Neumann}
\frac{\Delta \tau_\text{OPT} (\vartheta)}{h} = \frac{2}{G_0(\pi/2, \vartheta)+G_0(\pi, \vartheta)}
\end{equation}
where $G_0$ is defined in~\eqref{eq:GN} with coefficients \eqref{eq:coeffN}.
\end{proposition}
The proof is analogous to Proposition \ref{prop:opt}.

\subsection{Non-rectangular domains}\label{sect:nonrec}
\begin{figure}
\begin{center}
\begin{minipage}{0.48\textwidth}
\includegraphics[width=0.99\textwidth]{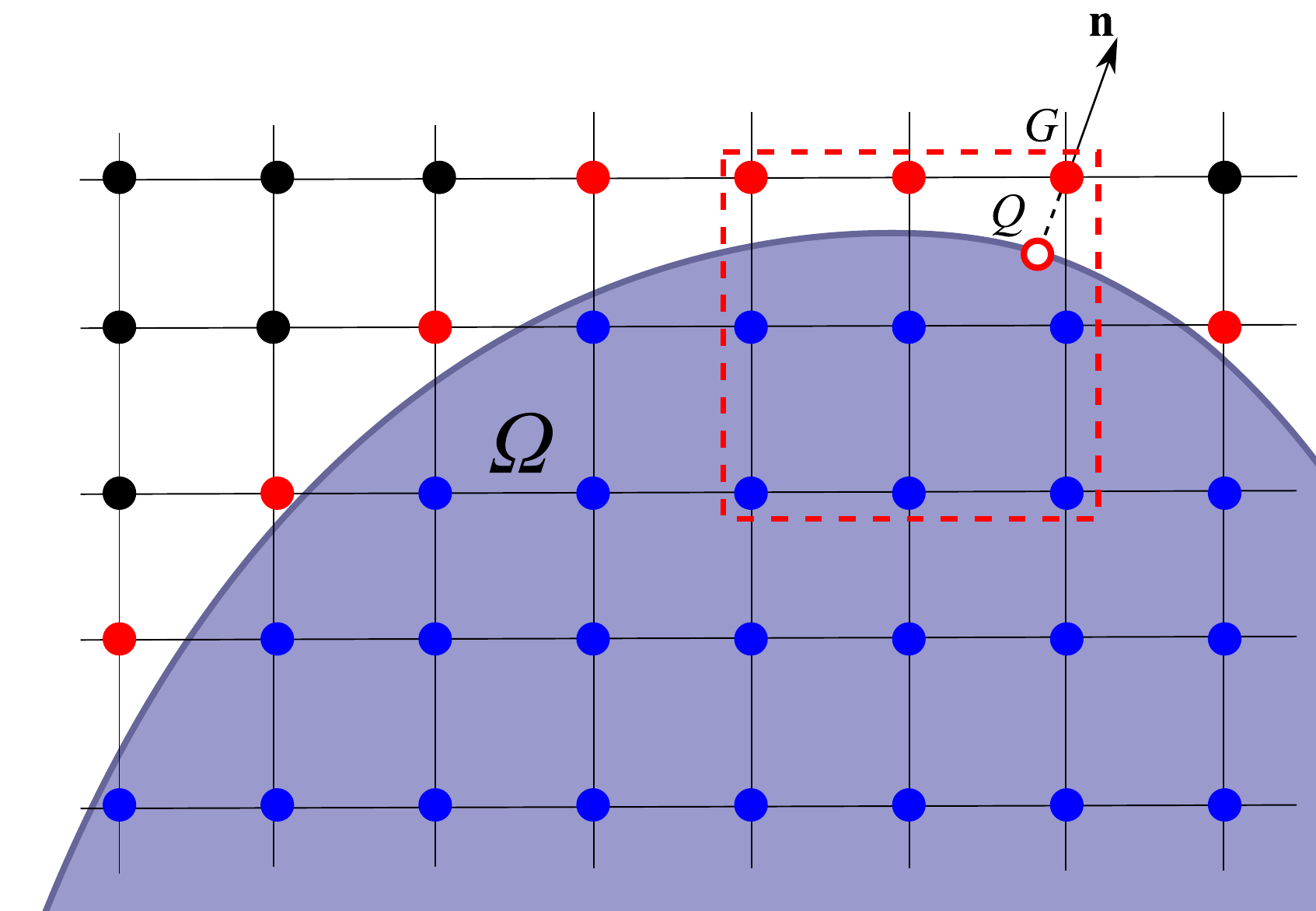}
\caption{
Partition of the computational grid points in \emph{internal points} (blue circles), \emph{ghost points} (red circles) and \emph{inactive points} (black circles). For each ghost point $G$ the orthogonal projection $Q$ is computed and the $3 \times 3$ stencil $ST-G$ (red dashed rectangle) is taken in upwind direction with respect to the outward normal $\vec{n}$.
}
\label{fig:Omega2D}
\end{minipage}
\hspace{0.1cm}
\begin{minipage}{0.48\textwidth}
\includegraphics[width=0.99\textwidth]{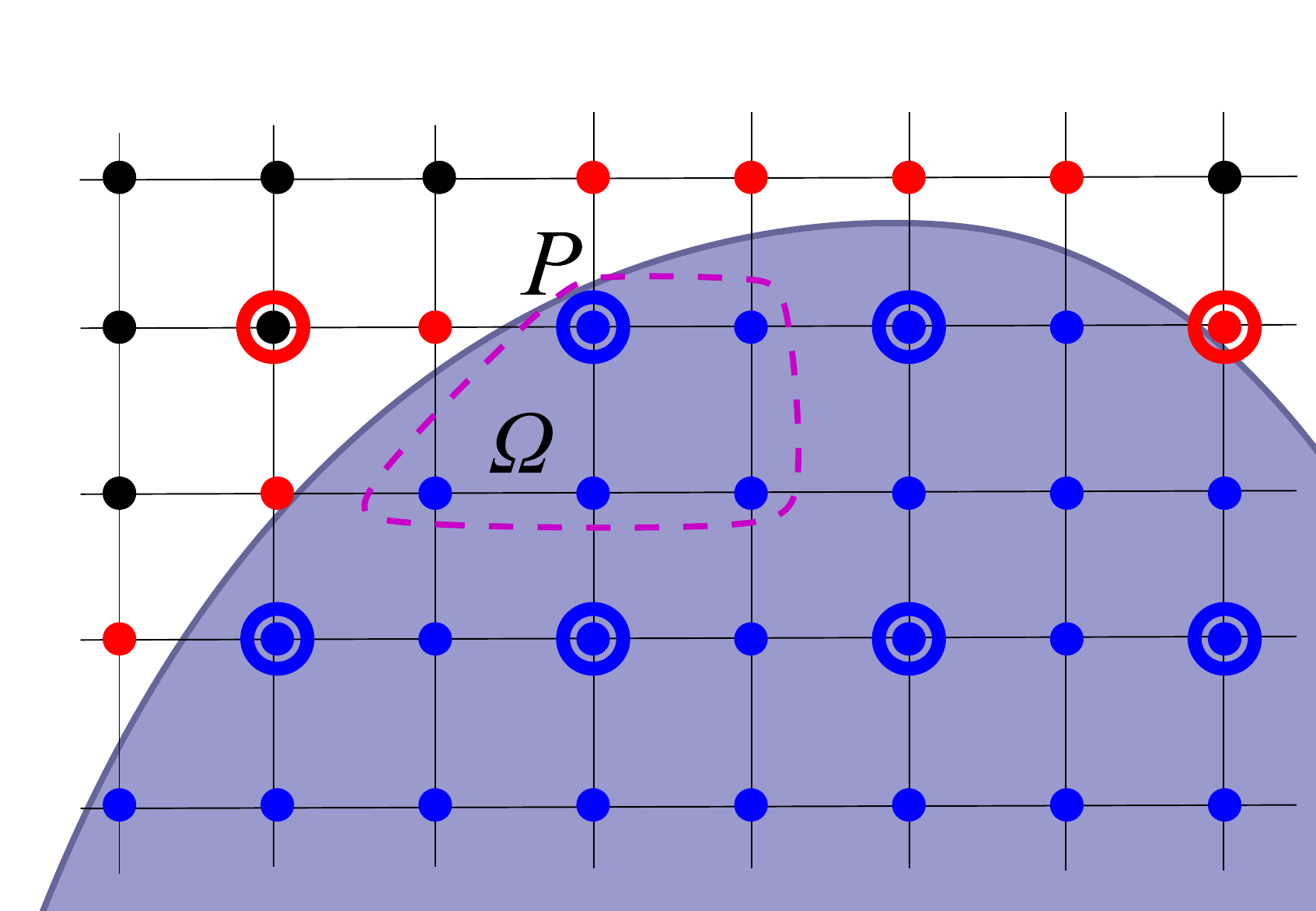}
\caption{
Restriction of the residual for a coarser internal point $P$ that is close to the boundary and then for which the FW nine-point stencil includes outside points. The large empty circles represent the grid points in the coarser grid (blue internal points, red ghost points and black inactive points). The restriction is performed as an average of the residual on the fine grid points enclosed by the purple dashed curve. Observe that some grid points can be inactive points in the fine grid and ghost points in the coarser grid. However, this behaviour does not affect the implementation of the overall BCMG method.
}
\label{fig:restriction2D}
\end{minipage}
\end{center}
\end{figure}

For general domains we adopt the ghost-point method and the level-set approach used 
in~\cite{CocoRusso:Elliptic}. In that paper, the authors adopted a multigrid approach with a fictitious-time dependent relaxation scheme using a constant $\Delta \tau$ along the boundary. Let us briefly recall the details.
As for the rectangular case, let $\mathcal{D} = [-1,1]^2$ be the computational domain covered with $(N+1)^2$ grid points $X_{ij} = (x_i,y_j)$, where $x_i = i\,h$ and $y_j=j\,h$ for $i,j=0,\ldots,N$ and $h=2/N$. Let $\mathcal{D}_h = \left\{ X_{ij} \text{ for } i,j=0,\ldots,N  \right\}$ be a Cartesian grid covering $\mathcal{D}$.
Let the domain $\Omega \subseteq \mathcal{D}$ be implicitly defined by a level-set function $\phi(x,y)$ as
\[
\Omega = \left\{ (x,y) \in \mathcal{D} \text{ such that } \phi(x,y)<0 \right\}
\]
and let $\Gamma = \partial \Omega$.
We solve the elliptic problem \eqref{mainProblem2D},
where $B$ represents a boundary operator, such as $Bu=u$ for Dirichlet boundary conditions and $Bu=\nabla u \cdot \vec{n}$ for Neumann boundary conditions, where $\vec{n}$ is the outward unit normal to $\Gamma$.
$\mathcal{D}_h$ is partitioned in three subsets: $\mathcal{D}_h = \Omega_h \cup \Gamma_h \cup I_h$, defined as follows. $\Omega_h$ is the set of \emph{internal grid points}, defined as $\Omega_h = \mathcal{D}_h \cap \Omega$; $\Gamma_h$ is the set of \emph{ghost points}, defined as:
\[
\Gamma_h = \left\{ (x_i,y_j) \in \mathcal{D}_h \text{ such that } \mathcal{N}_{ij} \cap \Omega_h \neq \emptyset   \right\},
\]
where $\mathcal{N}_{ij} = \left\{ (x_{i\pm1},y_{j}), (x_{i},y_{j\pm1}) \right\}$ is the set of four 
neighbour grid points of $(x_i,y_j)$; $I_h$ is the set of inactive points, defined as $I_h = \mathcal{D}_h \backslash (\Omega_h \cup \Gamma_h)$.
We order the internal and ghost points with, for example, a lexicographic order (a different order can also be chosen, as the method presented in this paper does not rely on a specific ordering).
Then $\Omega_h = \left\{ X_1, X_2, \ldots, X_{N_I} \right\}$ and
$\Gamma_h = \left\{ X_{N_I+1}, X_{N_I+2}, \ldots, X_{N_I+N_G} \right\}$, where
$N_I$ and $N_G$ are the total number of internal and ghost points, respectively.
The discrete solution of $\eqref{mainProblem2D}$ is represented by a column vector
$\vec{u}_h = \left[ U_1, U_2, \ldots, U_{N_I+N_G}\right]^T$.
For simplicity, we refer to the discrete value of the solution in a grid point as either $U_k$ or $u_{i,j}$, depending on the context.
For each internal grid point $X_k \in \Omega_h$ we discretize $-\Delta u = f$ using standard five-point stencil central differencing:
\begin{equation}\label{eq:Laph}
\frac{4u_{i,j}-u_{i-1,j}-u_{i+1,j}-u_{i,j-1}-u_{i,j+1}}{h^2}.
\end{equation}
For each ghost point $G=X_k \in \Gamma_h$ we discretize the boundary condition $Bu = g_\Gamma$ as follows (details can be found in~\cite{CocoRusso:Elliptic}). Firstly, we compute
the orthogonal projection $Q$ of $G$ onto the boundary $\Gamma$ (see Fig.~\ref{fig:Omega2D}) 
by solving $\phi(G - \eta \, h \, \vec{n}) = 0$ for $\eta \in [0,1)$, where $\vec{n}$ is the outward unit normal vector that can be approximated by a finite-difference discretization of
$\vec{n} = \nabla \phi / |\nabla \phi|$ centered at the ghost point.
Secondly, we identify the $3 \times 3$ stencil $ST_G$ that encloses $Q$ and has $G$ in one of its four 
corners. In detail, $ST_G$ is taken in an upwind fashion with respect to the normal direction (see Fig.~\ref{fig:Omega2D}).
Such a large stencil ensures second order accuracy for the solution $u$ and for the reconstruction of the gradient $\nabla u$ for both Dirichlet and Neumann boundary conditions. If one is interested only in the accuracy of the solution and not of the gradient, a $2 \times 2$ stencil will suffice to achieve second order accuracy for Dirichlet boundary conditions (see~\cite{CocoRusso:Elliptic}).

Finally, the discretized boundary condition is  
\begin{equation}\label{eq:BCh}
\left. B \tilde{u}_h \right|_Q = g_\Gamma(Q)
\end{equation}
where $\tilde{u}_h$ is the biquadratic interpolant on the stencil $ST_G$.
Observe that, provided that the grid is sufficiently fine, $ST_G$ does not contain any inactive point.

Equations \eqref{eq:BCh} cannot be explicitly solved for the ghost values, since the stencil $ST_G$ can contain  ghost points other than $G$ and then these equations would be fully coupled each other. Equations \eqref{eq:Laph} and \eqref{eq:BCh} for internal and ghost values are then collected into a sparse linear system, that we denote as $A_h \vec{u}_h = \vec{f}_h$, with $M_h \in \mathbb{R}^{N_I+N_G,N_I+N_G}$, $\vec{f}_h \in \mathbb{R}^{N_I+N_G}$.
This linear system is solved by a multigrid approach. 
The main structure of the TGCS is the same as in the beginning of \S\ref{MG1D}. 
Smoother and transfer operators are defined as follows.

\subsubsection{Relaxation scheme.}
As in the rectangular case, the smoother consists of a Gauss-Seidel scheme on internal points and of a relaxation of the boundary conditions on each ghost point $G$:
\begin{equation}\label{eq:rel}
u^{(m+1)}_G = u^{(m)}_G + \Delta \tau (f_{i,j} - B \tilde{u}^{(m)}).
\end{equation}
In~\cite{CocoRusso:Elliptic}, the value of $\Delta \tau$ was chosen to satisfy a stability condition and ensure convergence of the iterative scheme to the numerical solution. 
However, choosing a fixed $\Delta \tau$ in 2D and 3D (e.g.~$\Delta \tau = 1$ and $\Delta \tau = h$ for Dirichlet and Neumann boundary conditions, respectively) for all ghost iterations, although it would lead to a convergent scheme, it would degrade the performance of the multigrid strategy since boundary effects will be quickly propagated into the interior of the domain. This effect is mainly caused by the fact that adjacent ghost values are updated at different speeds depending on their distance from the boundary, generating highly-oscillatory modes that are not easily dumped from the multigrid strategy. On the other hand, for rectangular domains ghost points are located at the same distance from the boundary and a degradation effect is also observed if $\Delta \tau$ is not properly chosen, as observed in \S\ref{sect:rectangle}.
A possible solution to overcome the issue is presented in~\cite{CocoRusso:Elliptic} by over-relaxing the ghost value equations, namely adding some extra relaxation sweeps on ghost points per each interior one. The computational effort introduced by over relaxations is negligible in uniform grids, but it might become more evident in adaptive grids where the majority of grid points could be  concentrated close to the boundary. 

In this paper, we adopt a heuristic approach based on the BLFA performed in \S\ref{sect:rectangle} for rectangular domains.
In detail, we perform the ghost iterations \eqref{eq:rel} by choosing (see~\eqref{dtauoptLFA}):
\begin{equation}\label{DtauOptTheta}
\Delta \tau = \Delta \tau_\text{OPT} (\tilde{\vartheta}) \text{ with }
\tilde{\vartheta} = \frac{\left\| G - Q \right\|_2 }{ h }
\end{equation}
where $\left\| \cdot \right\|_2$ is the Euclidean norm (distance between $G$ and $Q$) and $\Delta \tau_\text{OPT}(\tilde{\vartheta})$ is defined in \eqref{dtauoptLFA} and \eqref{dtauoptLFA_Neumann} for Dirichlet and Neumann boundary conditions, respectively.

\begin{figure}
\begin{center}
\begin{minipage}{0.69\textwidth}
\includegraphics[width=0.99\textwidth]{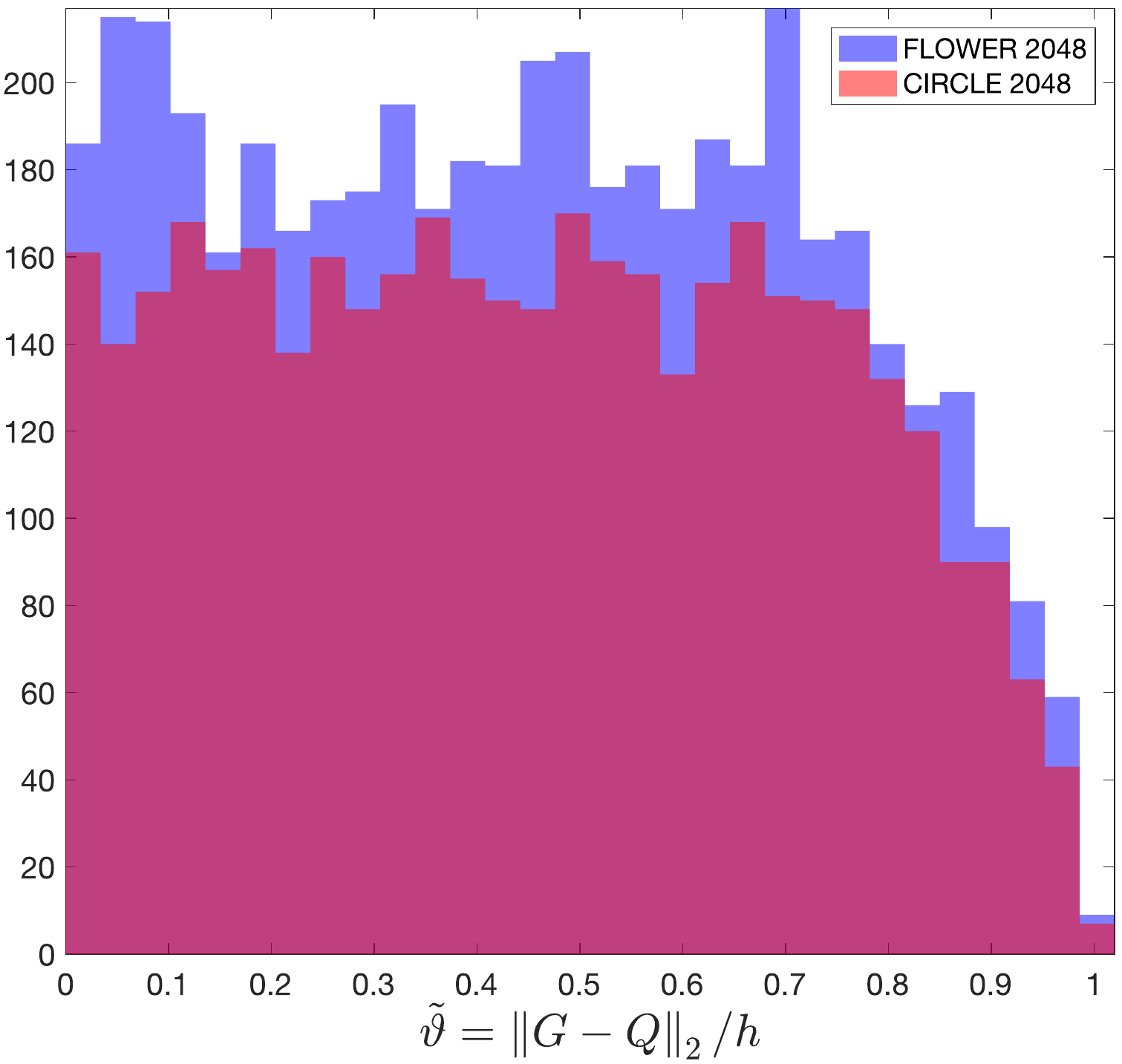}
\end{minipage}
\end{center}
\caption{
Distribution of $\tilde{\vartheta}$ between 0 and 1 for the circular and flower-shaped domains ($N=2048$).
}
\label{fig:hist_thetaxy}
\end{figure}

\begin{remark}\label{thetaxy}
For a curved domain, the values of $\tilde{\vartheta}$ are fully distributed between 0 and 1, highlighting the benefit of using a variable $\Delta \tau$ that depends on $\tilde{\vartheta}$ (see Fig.~\ref{fig:hist_thetaxy}).
We point out that the case of $\tilde{\vartheta}$ close to 1 is observed only for very specific configurations. Therefore, the distribution is almost uniform for $\tilde{\vartheta}$ far from 1 and smoothly decreases when $\tilde{\vartheta}$ approaches 1.
\end{remark}

\begin{remark}\label{exptheta2D}
We observe that \eqref{eq:G} and \eqref{eq:GN} can be thought as quadratic interpolations of $\exp(p(\alpha) x)$ and $p(\alpha) \exp(p(\alpha) x)$, respectively, on the boundary $x=x_V$. This suggests that we can use
$G_0(\alpha, \vartheta) = \exp(-p(\alpha) \vartheta)$ and 
$G_0(\alpha, \vartheta) = p(\alpha) \exp(-p(\alpha) \vartheta)$
for Dirichlet and Neumann boundary conditions instead of \eqref{eq:G} and \eqref{eq:GN}. 
In fact, \eqref{eq:G} and \eqref{eq:GN} where computed from the exact discretizations on the rectangular domain and there is no specific benefit to use them in a non-rectangular domain, where the discretization is more articulated. The approximations proposed in this remark are therefore adopted for non-rectangular domains and we observed numerically that they actually perform slightly better than \eqref{eq:G} and \eqref{eq:GN} (not shown).
\end{remark}

\subsubsection{Transfer operators.}
The residual of internal iterations is transferred to the coarser grid by the FW operator \eqref{eq:coarseFW} if all nine points of the FW stencil are all internal points, otherwise it is performed as an average of the residual only on the internal points of the nine-point FW stencil (see Fig.~\ref{fig:restriction2D}).
The rationale is that the residual of internal and ghost iterations scale with a different power of $h$ and then one should not mix them up, otherwise efficiency and sometimes even convergence cannot be guaranteed.

The residual of the ghost iterations is firstly extrapolated along the normal direction to the inactive points and then transferred to the coarser grid in the same way as the internal residual, namely by performing the FW if the nine points are all ghost or inactive points or as an average of external points only if any of the nine points is an internal one. Extrapolation to inactive points can be achieved, for example, by performing few iterations to solve
$
\frac{\partial r}{\partial \sigma} + \frac{\partial r}{\partial \vec{n}}=0
$
on each inactive point (or at least in a narrow band close to the boundary), discretized by Euler explicit in $\sigma$ and central difference for the normal derivative.

Interpolation from the coarse to the fine grid is performed by the standard bilinear interpolation operator on the whole computational domain without distinguishing between internal and ghost points (after extrapolating the ghost values along the normal direction to the inactive points). The rational is that the error to be interpolated is the solution of the coarse grid residual equation and then it is expected to scale with the same order for internal and ghost points.

\begin{remark}\label{remark:dtau}
If we want to use a constant $\Delta \tau$ for all ghost points and aim at achieving, if possible, a convergence factor of $\rho_\text{TARGET}$, 
the best that we can infer from the numerical results of the rectangular case is checking whether there exists a value of $\Delta \tau$ such that $\rho_\text{REL-BC}(\vartheta, \Delta \tau) \leq \rho_\text{TARGET}$ for all $\vartheta \in [0,1)$, where $\rho_\text{REL-BC}(\vartheta, \Delta \tau)$
is the convergence factor of $\texttt{REL-BC}$ in the rectangular domain of \S\ref{sect:rectangle}.
From Fig.~\ref{fig:2D_rhoD}, we can see that the minimum $\rho_\text{TARGET}$ for which we have a value of $\Delta \tau$ that satisfies the condition above is $\rho_\text{TARGET} \approx 0.3$, achieved for $\Delta \tau \approx 1.75$ for Dirichlet boundary conditions and 
$\rho_\text{TARGET} \approx 0.15$, achieved for $h^{-1} \Delta \tau \approx 1$ for Neumann boundary conditions.
This is confirmed by the numerical results presented in this paper on non-rectangular domains, for example Figs.~\ref{fig:2D_RHO_CIRCLE_FLOWER}, \ref{fig:2D_LINE30_N} and \ref{fig:2D_CIRCLE_FLOWER_N}. In particular, the use of the optimal $\Delta \tau$ proposed in the \textsc{BCMG} method is more important for Dirichlet than for Neumann boundary conditions.
\end{remark}

\subsubsection{Numerical results in 2D: Dirichlet boundary conditions}\label{test:nonrec_D}
\begin{figure}
\begin{center}
\begin{minipage}{0.49\textwidth}
\includegraphics[width=0.99\textwidth]{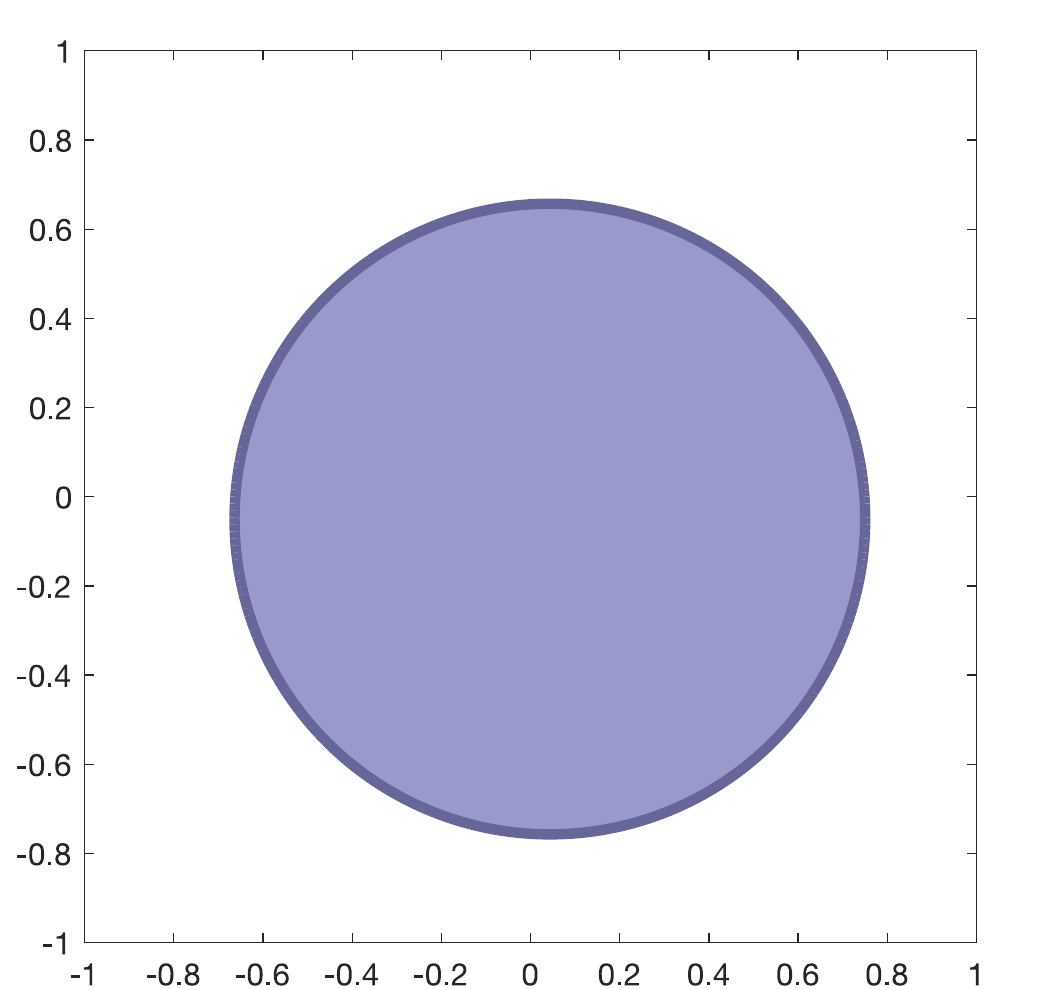}
\end{minipage}
\begin{minipage}{0.49\textwidth}
\includegraphics[width=0.99\textwidth]{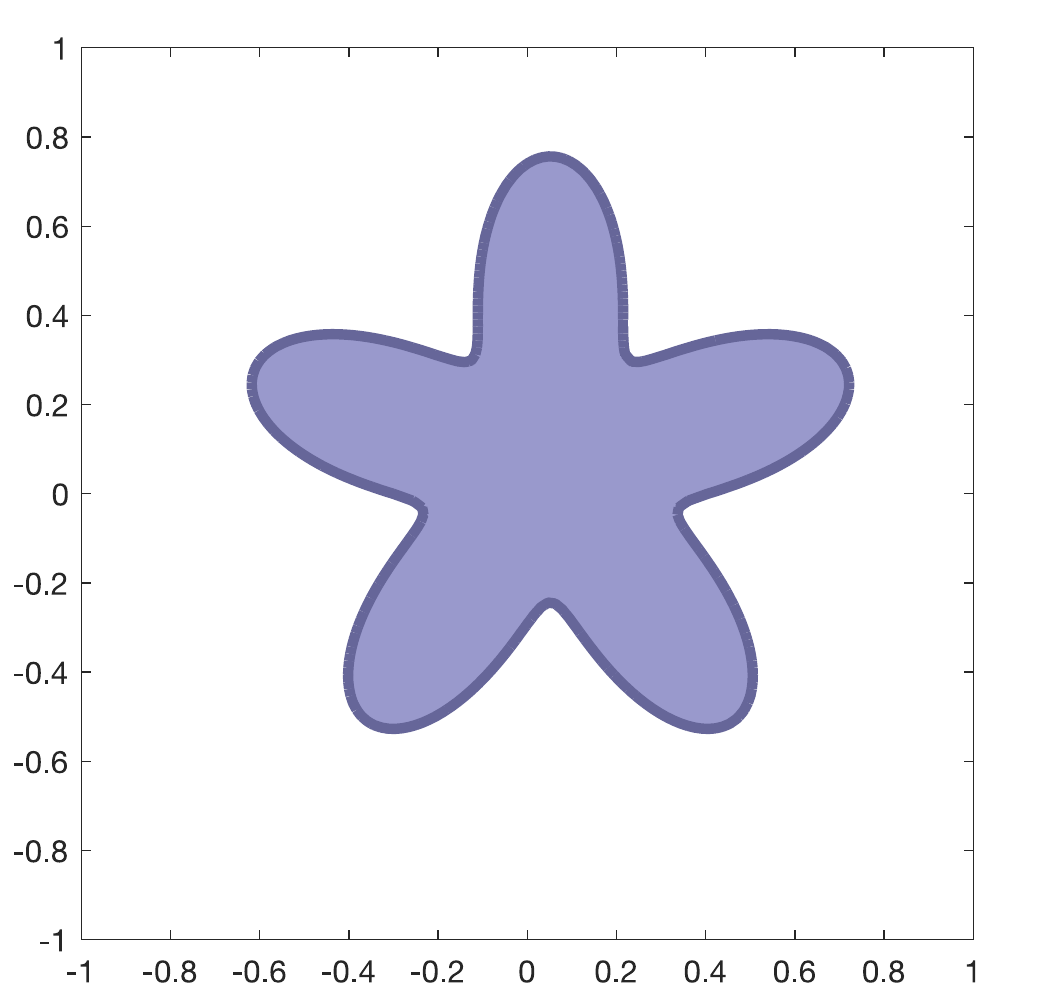}
\end{minipage}
\end{center}
\caption{
Domains used in the two-dimensional tests.
}
\label{fig:domains}
\end{figure}

We perform two tests. The first employs a circular domain with level-set function (see left panel of Fig.~\ref{fig:domains})
\begin{equation}\label{dom2D:circle}
 \text{\sc{(Circle)}} \qquad \phi(x,y) = (x-x_0)^2 + (y-y_0)^2 - R^2 \qquad
\text{ with } \; x_0=\sqrt{3}/40, \; y_0=-\sqrt{4}/40, \; R=\sqrt{2}/2
.
\end{equation}
The second one is a flower-shaped domain defined by
\begin{equation}\label{dom2D:flower}
\begin{aligned}
     \text{\sc{(Flower)}} \qquad  \phi(x,y) &= r-r_1-\frac{r_2}{r^5} \left( (y-y_0)^5+5(x-x_0)^4 (y-y_0)-10(x-x_0)^2 (y-y_0)^3 \right), \\
    & \text{ with } r = \sqrt{x^2+y^2}, \quad r_1=0.5, \quad r_2=0.2,
 \end{aligned}
\end{equation}
that is the level-set representation of a boundary whose parametric representation is (see~\cite{CocoRusso:Elliptic,Gibou:Ghost})
\begin{equation*}
\left\{
\begin{aligned}
x(\theta) &= r(\gamma) \cos(\gamma) + x_0 \\
y(\theta) &= r(\gamma) \sin(\gamma) + y_0
\end{aligned}
\right.
\end{equation*}
with $\gamma \in [0,2\pi]$ and $r(\gamma) = r_1 + r_2 \sin(5 \gamma)$.

\begin{figure}
\begin{center}
\begin{minipage}{0.49\textwidth}
\includegraphics[width=0.99\textwidth]{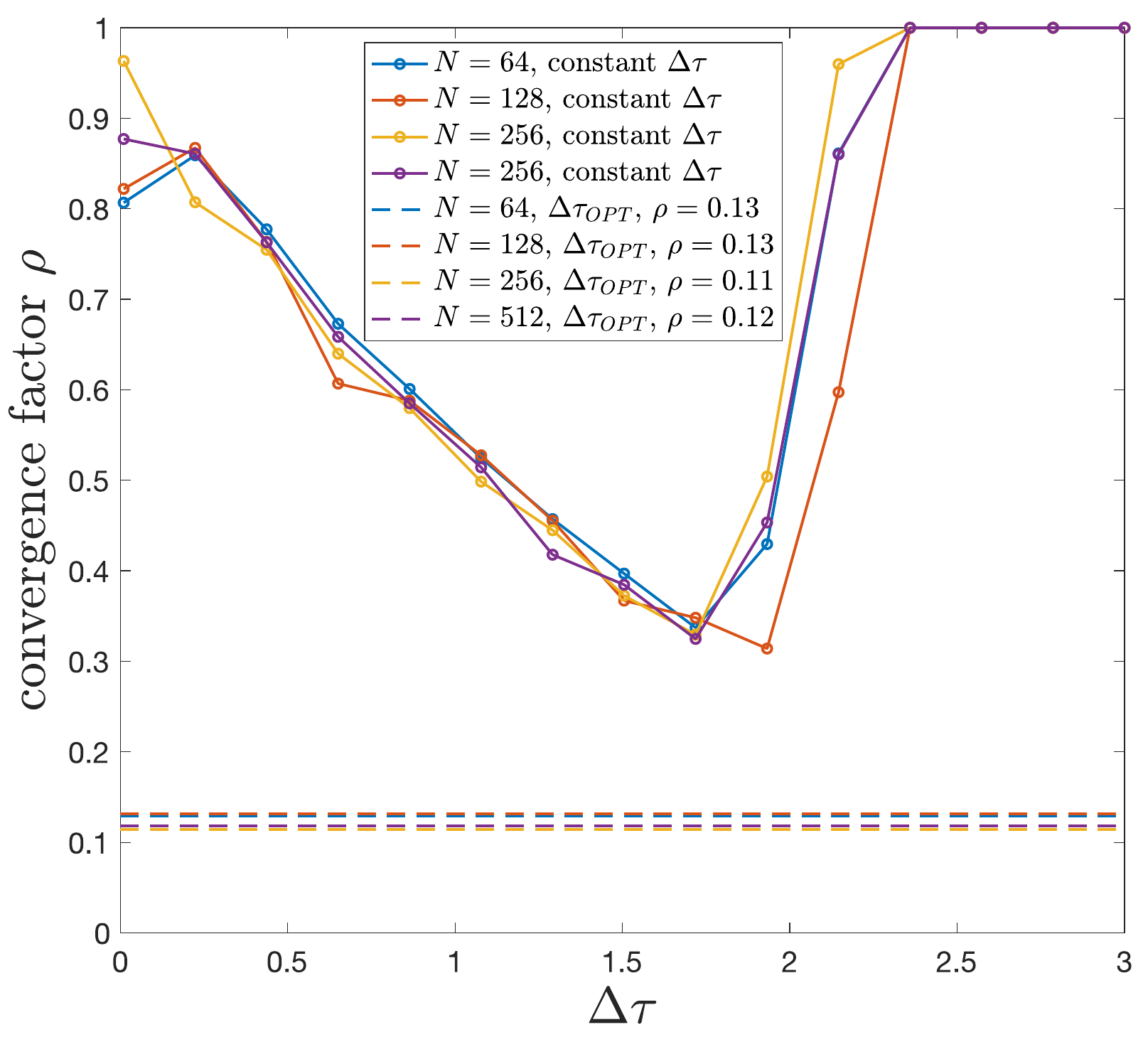}
\end{minipage}
\begin{minipage}{0.49\textwidth}
\includegraphics[width=0.99\textwidth]{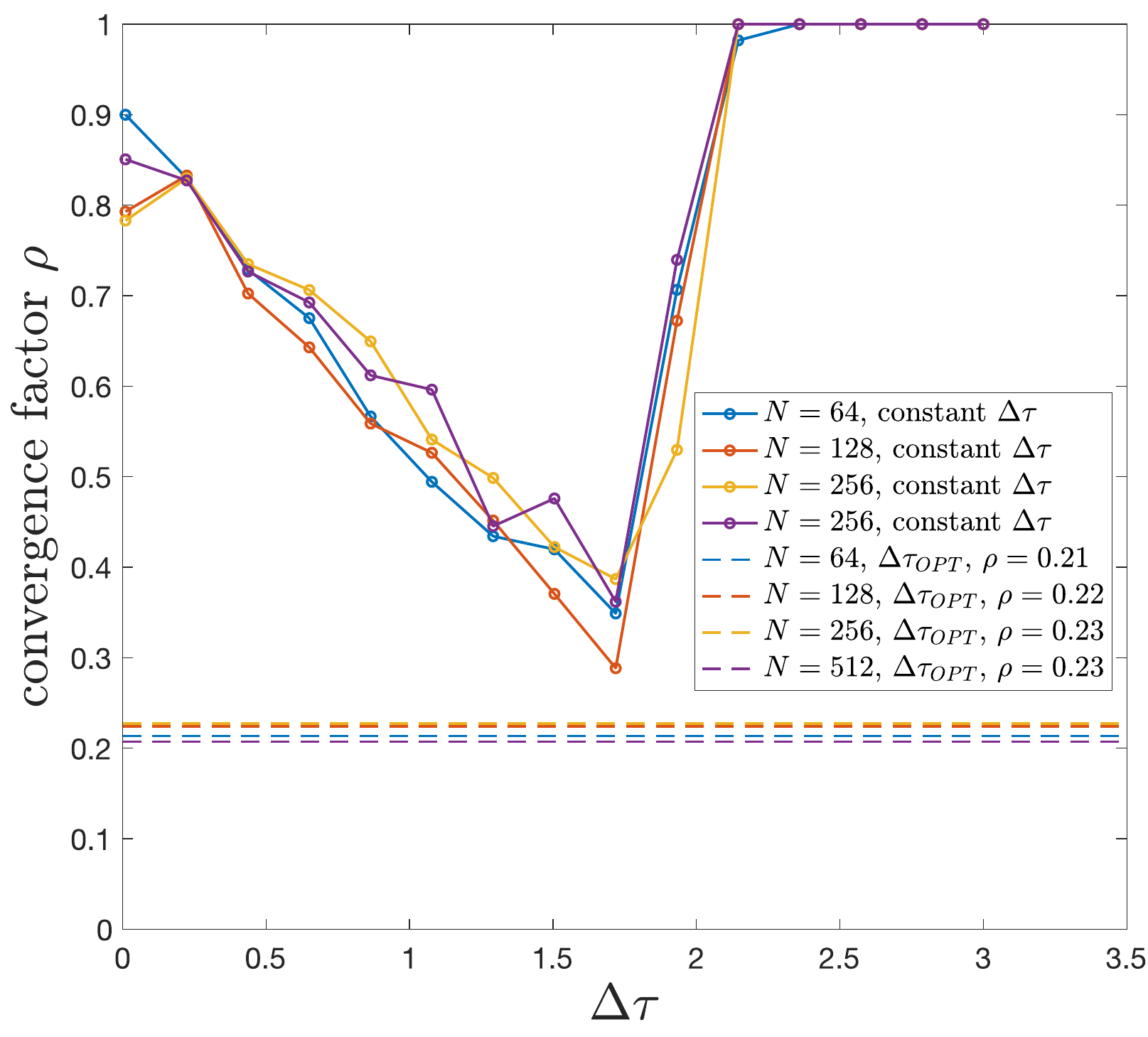}
\end{minipage}
\end{center}
\caption{
Two-dimensional problems with Dirichlet boundary conditions (\S\ref{test:nonrec_D}).
Left: circular domain (\textsc{Circle}). Right: flower-shaped domain (\textsc{Flower}).
Convergence factors of the TGCS (with setup from \S\ref{sect:setup}) for different values of $N$. Solid lines: tests with constant $\Delta \tau$. Dashed lines: with the optimal $\Delta \tau (\tilde{\vartheta})$ defined in Eq.~(\ref{DtauOptTheta}).
}
\label{fig:2D_RHO_CIRCLE_FLOWER}
\end{figure}

In the left panel of Fig.~\ref{fig:2D_RHO_CIRCLE_FLOWER} we plot the convergence factors of the TGCS (with set-up from \S\ref{sect:setup}) for different values of $N$. The solid line represents the test when $\Delta \tau$ is constant for all ghost points, while the dashed line is the test with the optimal $\Delta \tau (\tilde{\vartheta})$ defined in~\eqref{DtauOptTheta}. If we choose a constant $\Delta \tau$ for all ghost points, the best convergence factor that we can achieve is $\rho \approx 0.3$ obtained with $\Delta \tau \approx 1.75$ (see Remark~\ref{remark:dtau} for an explanation of these values). On the other hand, the implementation of a variable $\Delta \tau (\tilde{\vartheta})$ defined in~\eqref{DtauOptTheta} considerably improves the efficiency. The plot in fact suggests that the convergence factor approaches (when $N$ goes to infinity) the optimal one predicted by the LFA for the standard case $\Omega = [-1,1]^2$, which is $\rho_\text{OPT} \approx 0.119$ (when $\nu_1=2$, $\nu_2=1$ and FW restriction operator is adopted).
Similar results are observe for the flower-shaped domain (right panel of Fig.~\ref{fig:2D_RHO_CIRCLE_FLOWER}), although the convergence factors $\rho$ for $\Delta \tau_\text{OPT}$ are larger than the ones observed in a circular domain.
However, we believe that $\rho$ approaches the optimal value $\rho_\text{OPT} \approx 0.119$ when $N$ goes to infinity.
The discrepancy between the two tests can be explained by the fact that the boundary curvatures are not the same. In other words, the convergence factors between that two tests can be compared only when the ratio between boundary curvature and spatial step is similar. Also, it might be that the heuristic approach is a first approximation that neglects some additional features that are not accounted for in the present BCMG method, such as the relation with the local curvature of the boundary. A comprehensive theoretical analysis of the non-rectangular case is much more cumbersome and out of scope here.

In summary, we can suggest that the convergence factors of the TGCS is not much degraded by the ghost-point technique, provided that $N$ is sufficiently large (or, in other words, the resolution is greater near the boundary as ensured, for example, in the AMR framework).

\subsubsection{Comparison with the AGMG method}
\begin{figure}
\begin{center}
\begin{minipage}{0.49\textwidth}
\includegraphics[width=0.99\textwidth]{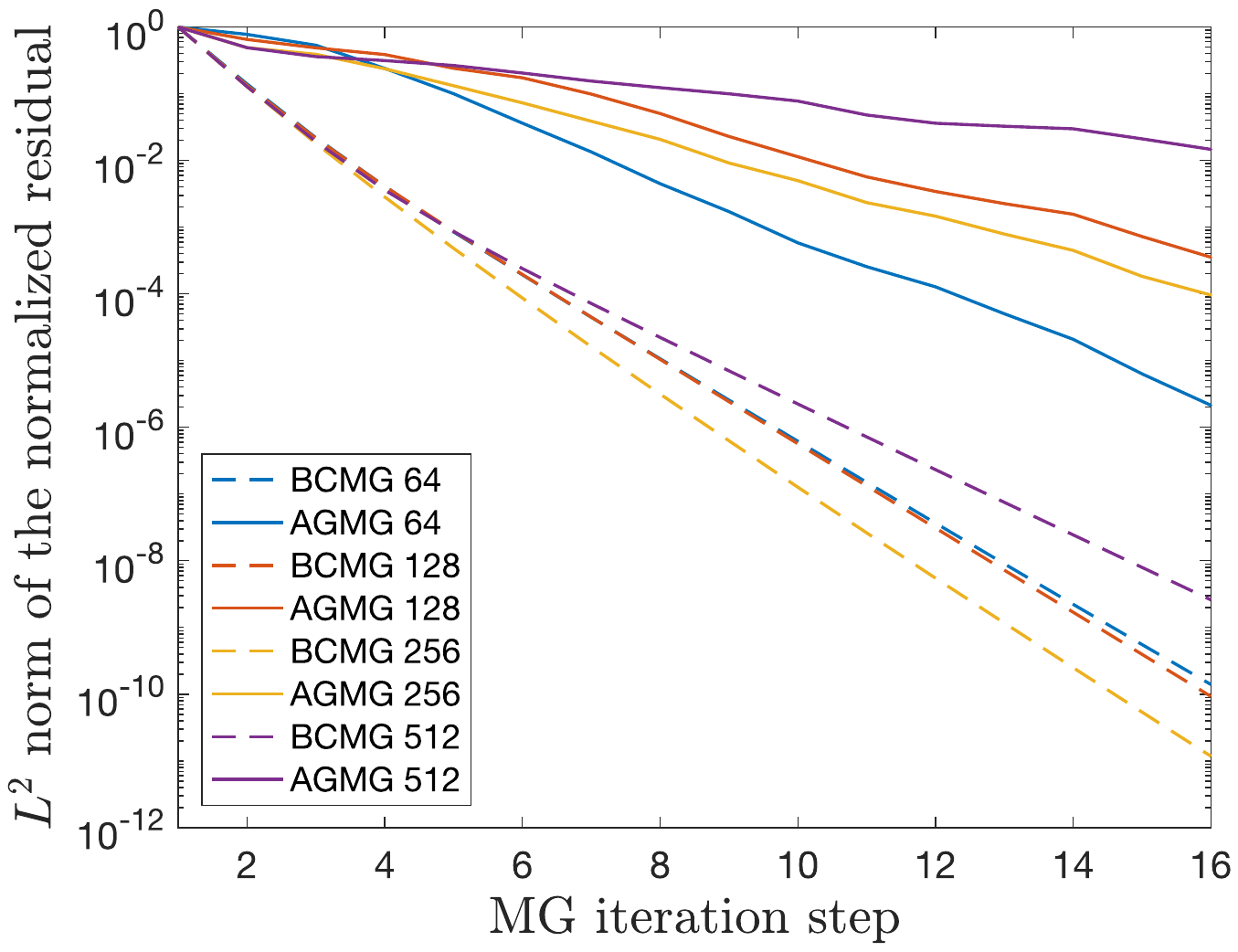}
\end{minipage}
\begin{minipage}{0.49\textwidth}
\includegraphics[width=0.99\textwidth]{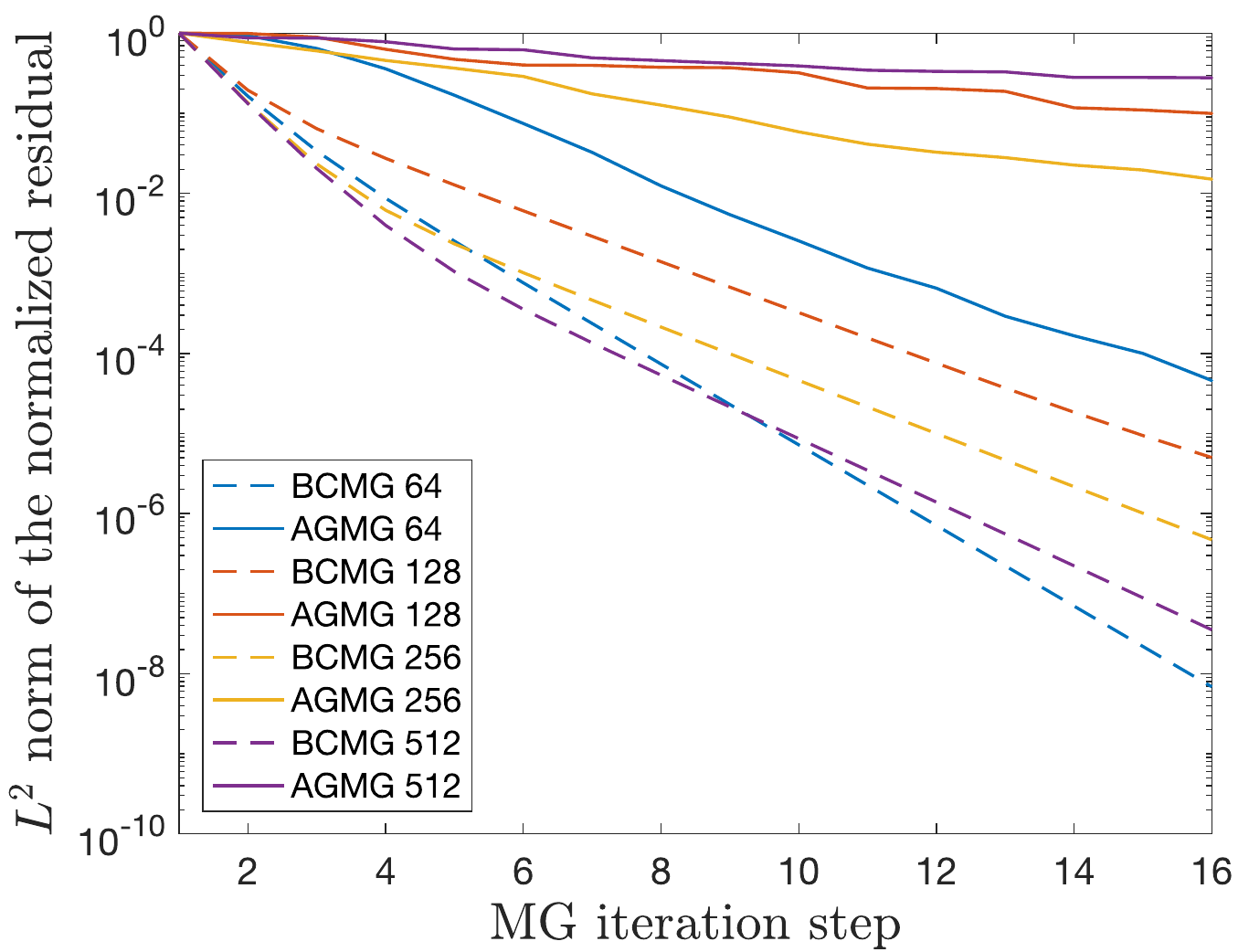}
\end{minipage}
\end{center}
\caption{
Comparison between BCMG and AGMG on the circular (left) and flower-shaped (right) domains of Fig.~\ref{fig:domains} with Dirichlet boundary conditions. One pre- and post- smoothing iterations are performed for both methods. The relative residual (normalized with respect to the first iteration) is plotted against each MG iteration for different sizes of the problem: $N = 64,128,256,512$.
}
\label{fig:comparisonAGMG}
\end{figure}

As a benchmark, we compare the performance of our BCMG method with that of an advanced algebraic multigrid method, namely, AGMG, cf.~\cite{notay-sw,notay1,notay2,notay3}.
AGMG is an aggregation-based algebraic multigrid
method that  performs the so-called K-cycle, i.e., two
Krylov accelerated iterations at each intermediate level, and one forward and
one backward Gauss-Seidel sweep for pre- and post-smoothing, respectively. 
AGMG is a powerful software designed to work as a black-box algorithm to solve efficiently sparse linear systems arising from discretizations of elliptic-like differential equations. Great flexibility is allowed, since it does not rely on the type of mesh (structured/unstructured, fitted/unfitted, etc.) or the particular discretization of the boundary conditions.
    This flexibility can come to a cost of under-performing against solvers that are designed for a specific problem. In Fig.~\ref{fig:comparisonAGMG} we compare AGMG with BCMG, showing that BCMG performs better than AGMG both on the circular and flower-shaped domains of Fig.~\ref{fig:domains}.
    The asymptotic slope of the lines is an approximation of the asymptotic convergence factor, that is observed numerically to be around $0.2$ for BCMG and $0.5$ for AGMG for the circular domain, while it is around $0.3$ for BCMG and $0.8$ for AGMG for the flower-shaped domain. Moreover, the initial slope of BCMG is even steeper, suggesting that the initial convergence factor is dictated by the internal relaxations, while the asymptotic convergence factor is affected by small boundary effects (this behaviour was actually observed numerically in all tests performed in this paper). On the other hand, AGMG starts slowly, with a shallower slope compared to later iteration steps. All these aspects suggest that BCMG will fall below a desired tolerance on the residual much earlier than AGMG, as shown in Fig.~\ref{fig:comparisonAGMG}.

\subsubsection{Numerical results in 2D: Neumann boundary conditions}\label{test:nonrec_N}
If Neumann boundary conditions are prescribed on the whole boundary, then the problem does not admit a unique solution. In particular, it can admit either infinite solutions (up to an additive constant) or no solutions (when the boundary condition and the source are not compatible). This property is also reflected at the discrete level, where the matrix of coefficients is singular.
The standard technique to approach these problems is to deal with an augmented linear system, with an additional unknown (that represents a measure of the compatibility condition misfit) and an additional equation (usually the average of the solution must vanish to guarantee uniqueness), leading to a non-singular matrix of coefficients. These aspects, together with the standard multigrid techniques to face these scenarios, are detailed, for example, in~\cite{Trottemberg:MG}. In this section we want to focus on the optimality of the relaxation parameters for Neumann boundary conditions, without generating unnecessary complications from the singularity of the problem that can be treated by standard techniques as anticipated above. To this purpose, we circumvent the singularity aspect in two ways: (i) by prescribing mixed boundary conditions, namely both Dirichlet and Neumann boundary conditions are applied on two different parts of the boundary, respectively (\textsc{Test1N}); (ii) by adding a reaction term to the Poisson equation, namely the differential equation of \eqref{mainProblem2D} becomes: $u - \Delta u = f$ (\textsc{Test2N}).

\begin{figure}
\begin{center}
\begin{minipage}{0.49\textwidth}
\includegraphics[width=0.99\textwidth]{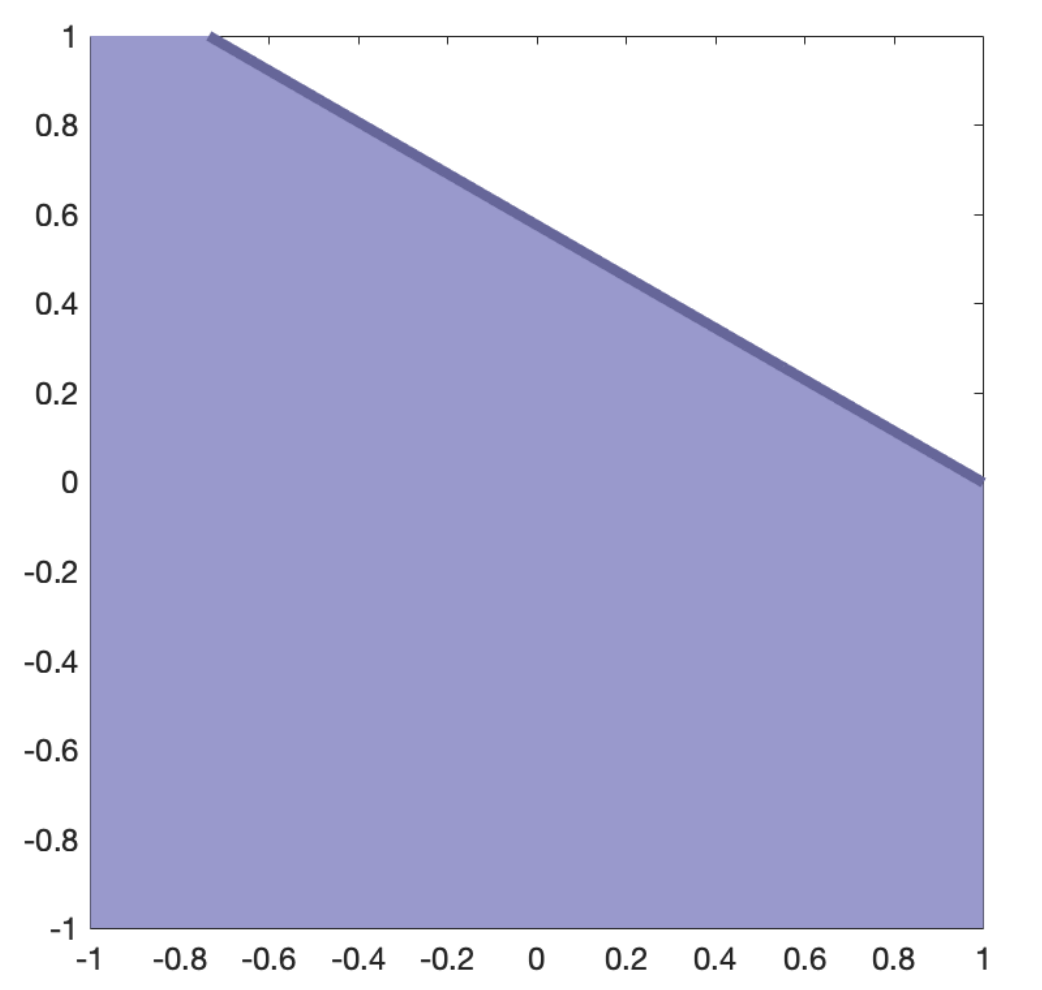}
\end{minipage}
\begin{minipage}{0.49\textwidth}
\includegraphics[width=0.99\textwidth]{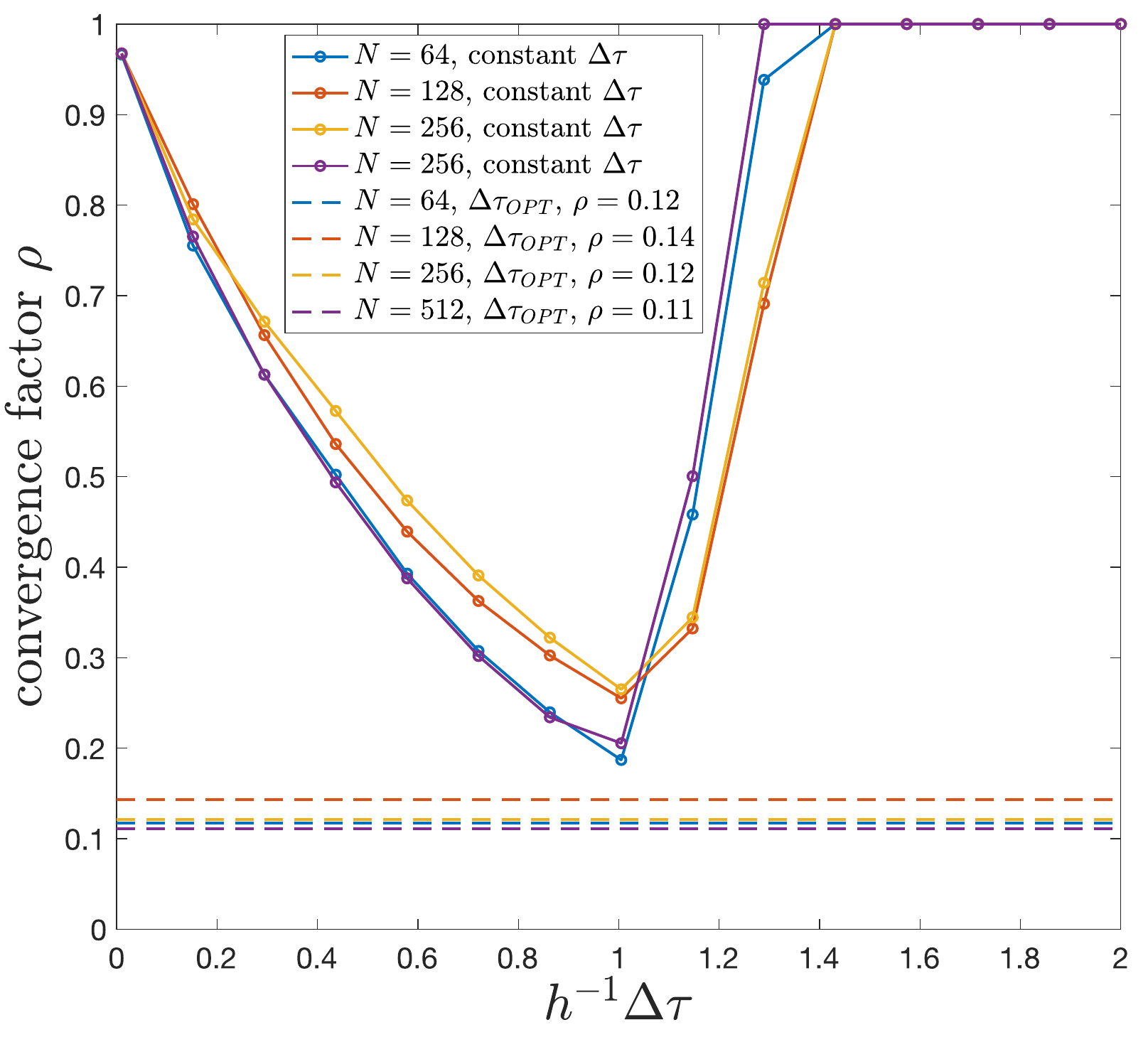}
\end{minipage}
\end{center}
\caption{
{\scshape Test1N}. 
Left: domain (Neumann boundary conditions are applied on the thick oblique line). 
Right: asymptotic TGCS convergence factors.
}
\label{fig:2D_LINE30_N}
\end{figure}

\textsc{Test1N} can be achieved in the problem \eqref{mainProblem2D} considering a domain for which $\Gamma_\mathcal{B} \neq \emptyset$ and on which Dirichlet boundary conditions are prescribed, while we apply Neumann boundary conditions on $\Gamma$.
In particular, we choose a level-set $\phi$ that represents a straight line passing through the point of coordinates $(1-7h/10,0)$ and with a negative slope of $-1/\sqrt{3}$:
\[
\phi(x,y) = \frac{1}{2} \left( x+\sqrt{3}y -1+\frac{7}{10} h \right).
\]
With this choice, the distance between ghost points and the boundary is not constant; in fact it is highly varying along the tangential direction, ensuring that the main difficulties posed by non-rectangular domains are present in this numerical test. 
The domain is represented in Fig.~\ref{fig:2D_LINE30_N}, left panel, while in the right panel we plot the convergence factors of the TGCS (with set-up from \S\ref{sect:setup}) for different values of $N$. The solid line represents the test when $\Delta \tau$ is constant for all ghost points, while the dashed line is with the optimal $\Delta \tau (\tilde{\vartheta})$ defined in Eq.~\eqref{DtauOptTheta}. With a constant $\Delta \tau$ for all ghost points, the best convergence factor that we can achieve is $\rho \approx 0.25$ obtained with $h^{-1} \, \Delta \tau \approx 1$ (see Remark~\ref{remark:dtau} for an explanation of these values). Efficiency is improved with a variable $\Delta \tau (\tilde{\vartheta})$ as defined in \eqref{DtauOptTheta}, with a convergence factor approaching the optimal one \eqref{eq:rhoopt2D} predicted by the LFA for the standard case $\Omega = [-1,1]^2$, that is $\rho_\text{OPT} \approx 0.119$.

\begin{figure}
\begin{center}
\begin{minipage}{0.49\textwidth}
\includegraphics[width=0.99\textwidth]{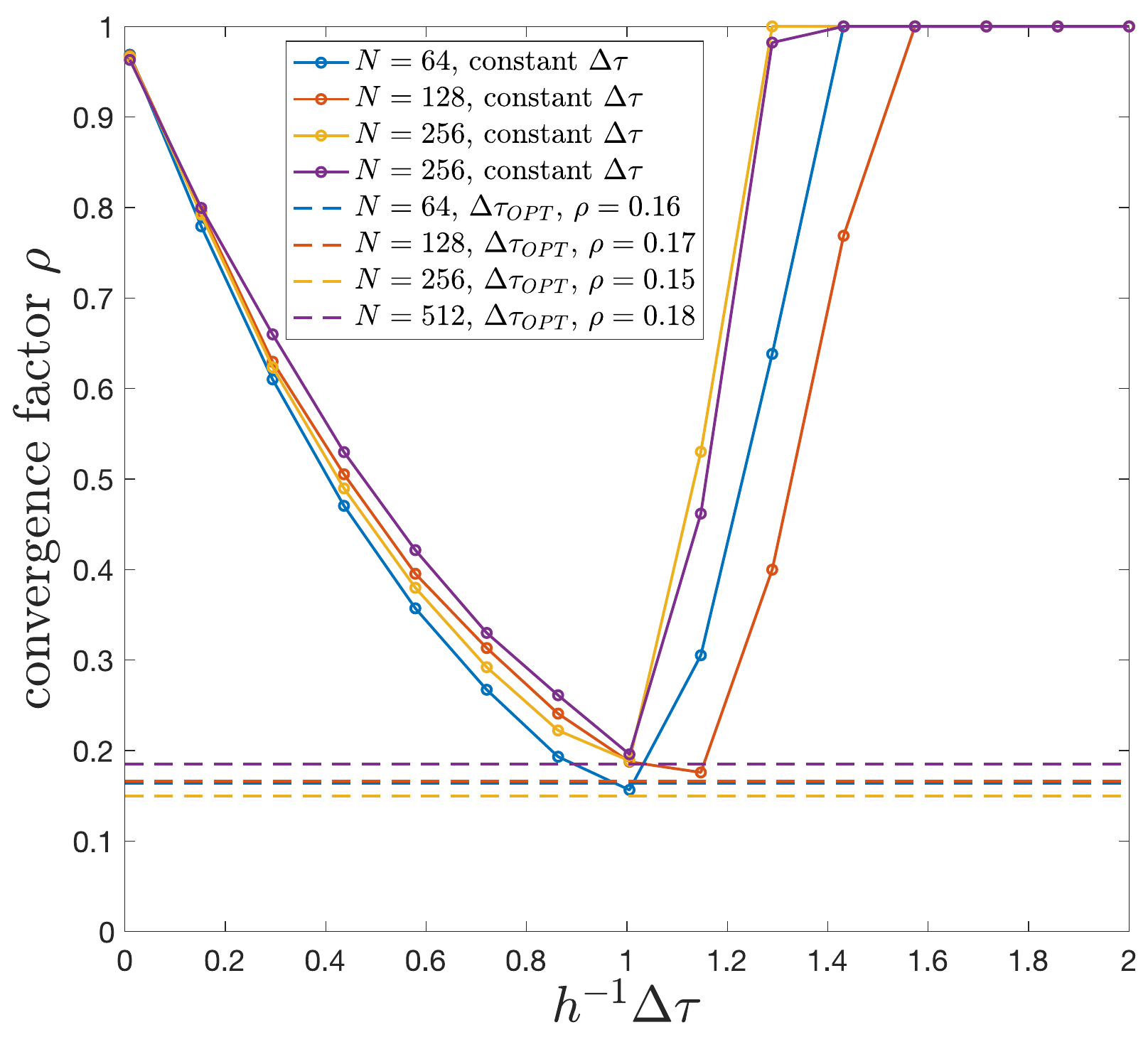}
\end{minipage}
\begin{minipage}{0.49\textwidth}
\includegraphics[width=0.99\textwidth]{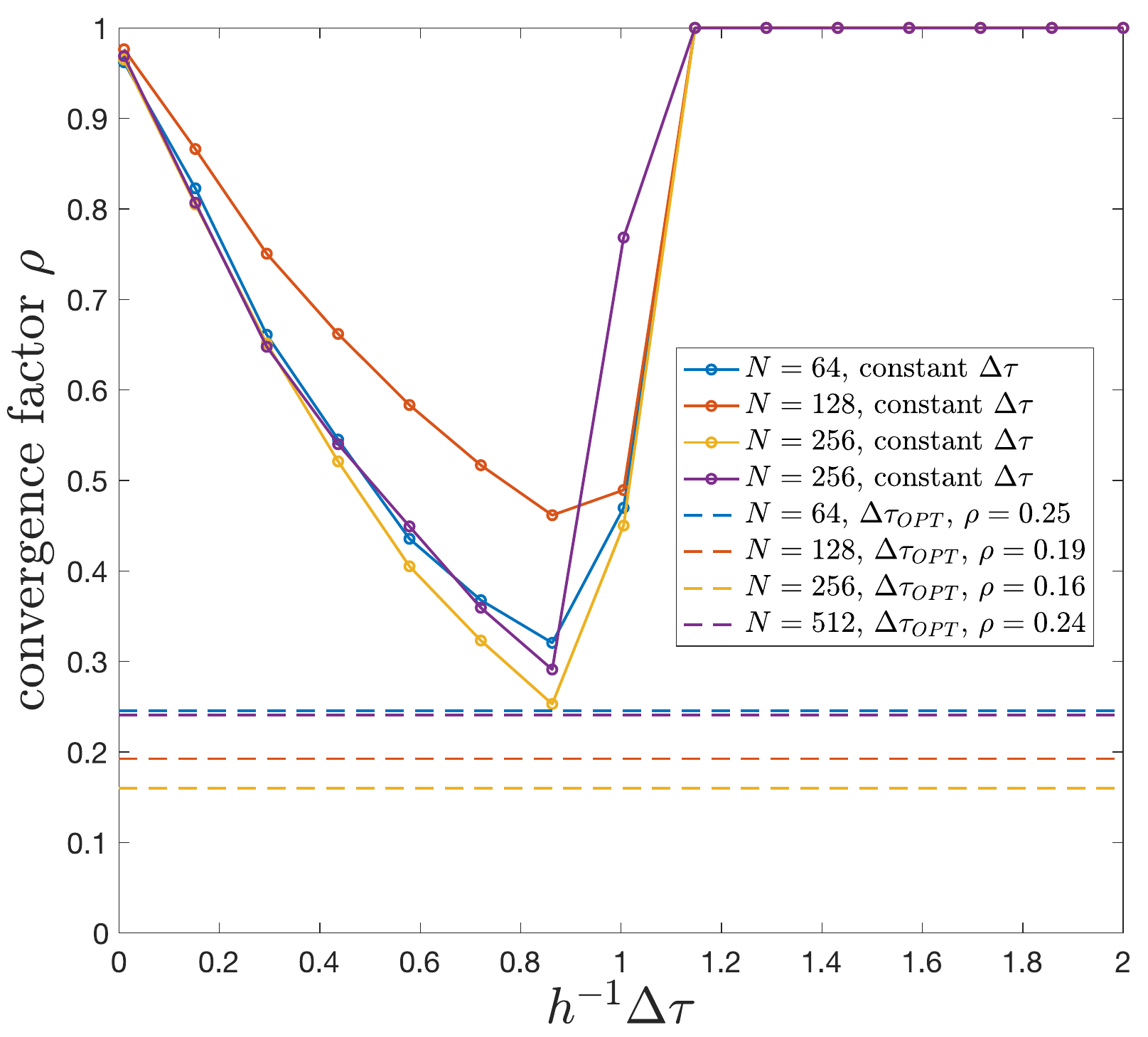}
\end{minipage}
\end{center}
\caption{
{\scshape Test2N}. 
Asymptotic TGCS convergence factors for Neumann buondary conditions on a reaction-diffusion equation. Left: circular domain. Right: flower shaped domain.
}
\label{fig:2D_CIRCLE_FLOWER_N}
\end{figure}

In \textsc{Test2N} we consider a reaction-diffusion equation with homogeneous Neumann boundary conditions on the boundary of the circular domain \eqref{dom2D:circle} and the flower-shaped domain \eqref{dom2D:flower}. Results are presented in Fig.~\ref{fig:2D_CIRCLE_FLOWER_N}, with a similar behaviour of \textsc{Test1N} discussed above, although with slightly larger asymptotic convergence factors, especially for the flower-shaped domain. As for the Dirichlet case, this is due to the high curvature of the boundary and the convergence factor would improve for finer grids.

\section{Three-Dimensional case}\label{sect:3D}
As in the previous section, we use a uniform Cartesian grid with $h=\Delta x = \Delta y = \Delta z$. The computational domain is $\mathcal{D} = [-1,1]^3$, covered by a grid with $(N+1)^3$ grid points $X_{ijk} = (x_i,y_j,z_k)$, where $x_i = -1+i\,h$, $y_j=-1+j\,h$ and $z_k=-1+k\,h$ for $i,j,k=0,\ldots,N$ and $h=2/N$.
The elliptic problem is formally the same as~\eqref{mainProblem2D}.

As in 2D, we aim at comparing the convergence factors with the optimal one $\rho_\text{OPT}$ that is obtained for regular rectangular domains $\Omega = [-1,1]^3$ where the ghost point technique is not necessary. For $\nu_1=2$, $\nu_2=1$, GS-LEX as a smoother and FW restriction operator we have obtained numerically that
\begin{equation}\label{eq:rhoopt3D}
\rho_\text{OPT} \approx 0.155.
\end{equation}

We mimic the approach adopted for the analysis in two space dimension, but in 3D the boundary is a surface and then the residual of the boundary conditions is a discrete approximation of a two-dimensional function instead of a one-dimensional function as in the 2D case. This means that the BLFA will be based on a two-dimensional analysis.


\subsection{Rectangular domains}\label{sect:rectangle3D}
The boundary to be treated with the ghost point technique, in this case is a plane $x=x_V$. In particular we consider the rectangular domain
\[
\Omega = \left\{ (x,y,z) \colon 0 \leq x \leq x_V, 0\leq y,z \leq 1 \right\},
\]
where $1-h<x_V \leq 1$.
We call $\vartheta = (1-x_V)/h$. Then $0 < \vartheta \leq 1$.
The multigrid ingredients, such as relaxation scheme, restriction and interpolation operators, can be easily extended from the 2D treatment described in \S\ref{sect:rectangle}. Full details are provided in~\cite{CocoRusso:Elliptic}.

\subsubsection{Dirichlet boundary conditions: BLFA}\label{sect:3D_dirLFA}
\begin{figure}
\begin{center}
\begin{tikzpicture}[scale=1.5, every node/.style={scale=1.2}]
\filldraw[fill=gray!50]
(3.14*0.5,0) -- (3.14,0) -- (3.14,3.14) -- (0,3.14) -- (0,3.14*0.5) -- (3.14*0.5,3.14*0.5) -- cycle;
\filldraw[fill=gray!20]
(3.14*0.5,0) -- (3.14*0.5,3.14*0.5) -- (0,3.14*0.5) -- (0,0) -- cycle;
        \draw[thick,->] (0,0) -- (4.3,0);
        \draw[thick,->] (0,0) -- (0,4.3);
        \draw[very thick,-] (3.14*0.5,0) -- (3.14,0);
        \draw[very thick,-] (3.14,0) -- (3.14,3.14);
        \draw[very thick,-] (0,3.14) -- (3.14,3.14);
        \draw[very thick,-] (0,3.14*0.5) -- (0,3.14);        
        \draw[very thick,-] (0,3.14*0.5) -- (3.14*0.5,3.14*0.5);                
        \draw[very thick,-] (3.14*0.5,0) -- (3.14*0.5,3.14*0.5);                        
        \foreach \Point/\PointLabel in {(3.14*0.5,0)/\frac{\pi}{2}, (3.14,0)/\pi}
        \draw[fill=black] \Point circle (0.1) node[below=0.1cm] {$\displaystyle \PointLabel$};
        \foreach \Point/\PointLabel in {(0,3.14*0.5)/\frac{\pi}{2}, (0,3.14)/\pi}
        \draw[fill=black] \Point circle (0.1) node[left=0.1cm] {$\displaystyle \PointLabel$};         
        \foreach \Point in {(3.14*0.5,3.14*0.5), (3.14,3.14)}
        \draw[fill=black] \Point circle (0.1);       
        \foreach \Point/\PointLabel in {(3.14*0.5,0)/A, (3.14,0)/D, (3.14,3.14)/C, (3.14*0.5,3.14*0.5)/B, (0,3.14)/E, (0,3.14*0.5)/F}        
	\draw \Point node[above right] {$\PointLabel$};
	\draw (4.5,0) node[right] {$\alpha_1$};
	\draw (0,4.5) node[above] {$\alpha_2$};	
	\draw (3.14*0.25,3.14*0.25) node {$\text{LOW}$};		
	\draw (3.14*0.75,3.14*0.75) node {$\text{HIGH}$};			
    \end{tikzpicture}
\end{center}
\caption{Frequency domain splitting for bi-variate BLFA.}
\label{fig:2Dalpha}
\end{figure}
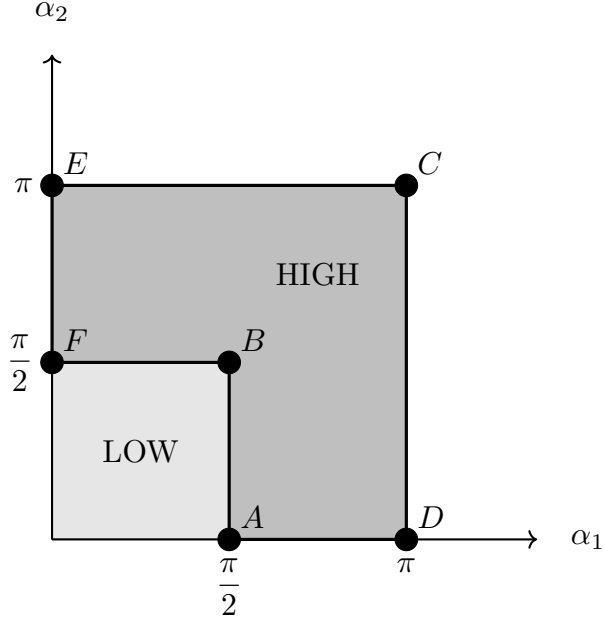

We start analysing the Dirichlet boundary condition. 
Following the same approach as in 2D, see \S\ref{2D_dirLFA}, we consider a semi-infinite domain
\[
\Omega_\infty = \left\{ (x,y,z) \colon x \leq x_V\right\}
\]
and the 2D Fourier modes
\[
\varphi(y,z) = e^{\iu \alpha_1 y/h} e^{\iu \alpha_2 z/h},
\]
with $\alpha_1, \alpha_2 \in [0,\pi]$.
The Fourier mode $\varphi$ is a highly-oscillatory mode if $\max \left\{ \alpha_1, \alpha_2 \right\} \in [\pi/2,\pi]$. We split the frequency space $[0,\pi]^2$, see Fig.~\ref{fig:2Dalpha}, into 
\[
\text{LOW} = [0,\pi/2]^2, \qquad \text{HIGH} = [0,\pi]^2 \backslash \text{LOW}.
\]

As in the 2D case (see \S\ref{2D_dirLFA}), the relaxation {\tt REL-BC} is performed in two steps. The first step consists of solving the internal grid points in block with ghost values $\varphi$:
\[
\frac{4 u_{i,j,k} - (u_{i+1,j,k}+u_{i-1,j,k}+u_{i,j+1,k}+u_{i,j-1,k})+
u_{i,j,k+1}+u_{i,j,k-1}}{h^2} = 0,
\quad 
\text{ for } i \leq N-1, \quad j,k\in\mathbb{Z},
\]
\[
\lim_{i \rightarrow -\infty} u_{i,j,k} = 0,
\qquad
u_{i,j,k} = \varphi(y_j,z_k) \text{ for } i=N.
\]
The exact solution of this difference equation is
\begin{equation}\label{eq:exactnum3D}
u_{i,j,k} = e^{p(\alpha_1,\alpha_2) (x_i-1)/h}\varphi(y_j,z_k) =e^{p(\alpha_1,\alpha_2) (x_i-1)/h} e^{\iu \alpha_1 y_j/h} e^{\iu \alpha_2 z_k/h}
,
\end{equation}
where 
\begin{equation}\label{alphaNalpha3D}
p(\alpha_1,\alpha_2) = \cosh^{-1} (3 - \cos(\alpha_1) - \cos(\alpha_2)).
\end{equation}

The second step consists of relaxing the ghost values by the 3D version of \eqref{relax2DV} with $g_{jk} = 0$:
\[
u_{N,j,k}^{(1)} = u_{Njk} - \Delta \tau \, \left( \frac{\vartheta (\vartheta-1)}{2} \, u_{N-2,j,k} + \vartheta (2-\vartheta) \, u_{N-1,j,k}+ \frac{(1-\vartheta)(2-\vartheta)}{2} u_{N,j,k} \right)
\]
\[
\Longrightarrow
\varphi^{(1)}(y_j,z_k) = \varphi(y_j,z_k) - \Delta \tau \, \left( \frac{\vartheta (\vartheta-1)}{2} \, e^{-2 p(\alpha_1,\alpha_2)}\varphi(y_j,z_k) + \vartheta (2-\vartheta) \, e^{-p(\alpha_1,\alpha_2)}\varphi(y_j,z_k)+ \frac{(1-\vartheta)(2-\vartheta)}{2} \varphi(y_j,z_k) \right)
\]
\begin{equation}\label{LFA3Dampl}
\Longrightarrow \frac{\varphi^{(1)}(y_j,z_k) }{ \varphi(y_j,z_k)} = 1 - \Delta \tau \, \left( \frac{\vartheta (\vartheta-1)}{2} \, e^{-2 p(\alpha_1,\alpha_2)} + \vartheta (2-\vartheta) \, e^{-p(\alpha_1,\alpha_2)}+ \frac{(1-\vartheta)(2-\vartheta)}{2} \right).
\end{equation}

The right-hand side of \eqref{LFA3Dampl} is the amplification factor of the Fourier mode $\varphi = e^{i \alpha y/h}$ for a given $\vartheta \in [0,1)$ and a given $\Delta \tau > 0$, that we denote as:
\begin{equation}\label{eq:G3D}
G(\alpha_1, \alpha_2, \vartheta, \Delta \tau) =  1 - \Delta \tau \; G_0(\alpha_1, \alpha_2, \vartheta)
\text{ with }
G_0(\alpha_1, \alpha_2, \vartheta) = 
\sum_{r=0}^2 c_{-r} e^{-r \, p(\alpha_1,\alpha_2)}
\end{equation}
with coefficients \eqref{eq:coeff}.

Finding the optimal $\Delta \tau$ means minimizing the largest amplification factor for the high-frequency modes:
\begin{equation}\label{eq:def_dtau_opt3D}
\Delta \tau_\text{OPT} (\vartheta)  = \arg \min_{\Delta \tau >0} \sup_{ (\alpha_1,\alpha_2) \in \text{HIGH} } |G(\alpha_1, \alpha_2, \vartheta, \Delta \tau)|.
\end{equation}

\begin{lemma}\label{lemma:fnostat3D}
Let $G_0(\alpha_1, \alpha_2, \cdot) \colon (\pi/2,\pi )^2 \to \mathbb{R}$ be the continuous function defined in~\eqref{eq:G3D} with coefficients \eqref{eq:coeff}.
Then $G_0$ has no stationary points in \text{HIGH}. 
\end{lemma}
\noindent {\it Proof.}

Let $r=1$ or $r=2$. We have:
\[
\frac{\partial G_0}{\partial \alpha_r} = - e^{- p(\alpha_1,\alpha_2)} \left( 2c_{-2} e^{-p(\alpha_1,\alpha_2)} + c_{-1} \right)
\frac{\partial p }{\partial \alpha_r}
\]
Since $\partial p / \partial \alpha_r =\sin(\alpha_r) \left( (3-\cos(\alpha_1)-\cos(\alpha_2))^2-1 \right)^{-1/2} \neq 0$ for any $(\alpha_1,\alpha_2) \in (0,\pi)^2$, then the only possible stationary points of $G_0(\alpha_1,\alpha_2, \cdot)$ satisfy the equation
\begin{equation}\label{eq:alpha3D}
e^{-p(\alpha_1,\alpha_2)} = \frac{c_{-1}}{-2c_{-2}} = 1 + \frac{1}{1-\vartheta}.
\end{equation}
Since $1 + (1-\vartheta)^{-1} > 1$ (because $0 \leq \vartheta <1$) and $e^{-p(\alpha_1,\alpha_2)}<1$ (because $p(\alpha_1,\alpha_2)>0$ for every $(\alpha_1,\alpha_2)$), then~\eqref{eq:alpha3D} is never satisfied, proving that $G_0(\alpha_1,\alpha_2,\cdot)$ has no stationary points in \text{HIGH}.
$\square$

\begin{lemma}\label{lemma:G0pos3D}
Let $G_0(\cdot, \cdot, \vartheta) \colon (0,1) \to \mathbb{R}$ be the continuous function of $\vartheta$ defined in~\eqref{eq:G3D} with coefficients \eqref{eq:coeff}. Then $G_0(\cdot, \cdot, \vartheta)>0$ for every $\vartheta \in [0,1)$.
\end{lemma}

\noindent {\it Proof.}
We have $G_0(\cdot, \cdot, 0) = 1>0$ and $G_0(\cdot, \cdot, 1) = e^{-p(\alpha_1,\alpha_2)}>0$. To complete the proof, it is sufficient to show that $G_0(\cdot, \cdot, \vartheta)$ has no stationary points in $(0,1)$.
We have
\[
\frac{\partial G_0}{\partial \vartheta} = \frac{1}{2} (e^{-p(\alpha_1,\alpha_2)}-1) \left(  (e^{-p(\alpha_1,\alpha_2)}-1) \vartheta - (e^{-p(\alpha_1,\alpha_2)}-3) \right)
\]
Since $p(\alpha_1,\alpha_2)>0$, we have that
\[
\frac{\partial G_0}{\partial \vartheta} = 0 \Longleftrightarrow \vartheta = \frac{e^{-p(\alpha_1,\alpha_2)}-3}{e^{-p(\alpha_1,\alpha_2)}-1}.
\]
We observe that $(e^{-p(\alpha_1,\alpha_2)}-3)(e^{-p(\alpha_1,\alpha_2)}-1)^{-1} > 1$ (because $p(\alpha_1,\alpha_2)>0$), proving that there are no stationary point of $G_0(\cdot, \cdot, \vartheta)$ in $(0,1)$.
$\square$

\begin{proposition}\label{prop:opt3D}
The optimal $\Delta \tau$ defined in~\eqref{eq:def_dtau_opt3D} is given by
\begin{equation}\label{dtauoptLFA3D}
\Delta \tau_\text{OPT} (\vartheta) = \frac{2}{G_0^\text{min}(\vartheta)+G_0^\text{max}(\vartheta)}
\end{equation}
with
\begin{equation}\label{eq:G0minmax}
G_0^\text{min/max}(\vartheta) = \min / \max \left\{ |G_0(\pi/2,0,\vartheta)|, |G_0(\pi/2,\pi/2,\vartheta)|, |G_0(\pi,\pi,\vartheta)| \right\}
\end{equation}
and $G_0$ defined in~\eqref{eq:G3D} with coefficients \eqref{eq:coeff}.
\end{proposition}

{\it Proof.}
From Lemma \ref{lemma:fnostat3D} we have that $G(\alpha_1, \alpha_2,\cdot, \cdot) = 1 - \Delta \tau \, G_0(\alpha_1, \alpha_2, \cdot)$ has no stationary points for $(\alpha_1,\alpha_2) \in \text{HIGH}$. Then
\begin{equation}\label{eq:3DsupG}
\sup_{ (\alpha_1,\alpha_2) \in \text{HIGH} } |G(\alpha_1, \alpha_2, \vartheta, \Delta \tau)| = 
\max_{V \in \left\{ A,B,C,D,E,F \right\}} \left\{ |G(V, \vartheta, \Delta \tau)|  \right\},
\end{equation}
where $A=(\pi/2, 0)$, $B=(\pi/2, \pi/2)$, $C=(\pi, \pi)$, $D=(\pi, 0)$, $E=(0,\pi)$, $F=(0,\pi/2)$ (see Fig.~\ref{fig:2Dalpha}). Observe that $G(\alpha_1, \alpha_2 \cdot, \cdot)$ is symmetric about the line $\alpha_1=\alpha_2$, so points $E$ and $F$ are reduntants in \eqref{eq:3DsupG} (since $G(E, \cdot, \cdot)=G(B, \cdot, \cdot)$ and $G(F, \cdot, \cdot)=G(A, \cdot, \cdot)$). Also, from Eq.~\eqref{alphaNalpha3D} we have that
$\alpha_N(B) = \alpha_N(D)= \cosh^{-1}(3)$,
meaning that $G(B, \cdot, \cdot)=G(D, \cdot, \cdot)$. 
Therefore, Eq.~\eqref{eq:3DsupG} becomes
\begin{equation}\label{eq:3DsupG_ABC}
\sup_{ (\alpha_1,\alpha_2) \in \text{HIGH} } |G(\alpha_1, \alpha_2, \vartheta, \Delta \tau)|
=
\max_{V \in \left\{ A,B,C \right\}} \left\{ |G(V, \vartheta, \Delta \tau)|  \right\}.
\end{equation}

Moreover, for every $V=A,B,C$ we have that $G(V,\cdot, \Delta \tau)$ is a linear function of $\Delta \tau$ and from Lemma \ref{lemma:G0pos3D} it has negative slope with magnitude $G_0(V, \cdot)$. Therefore, from Lemma \ref{lemma:lines} we obtain \eqref{dtauoptLFA}.
$\square$

\begin{figure}
\begin{center}
\begin{tabular}{cc}
(Dirichlet) & (Neumann) \\
\includegraphics[width=0.45\textwidth]{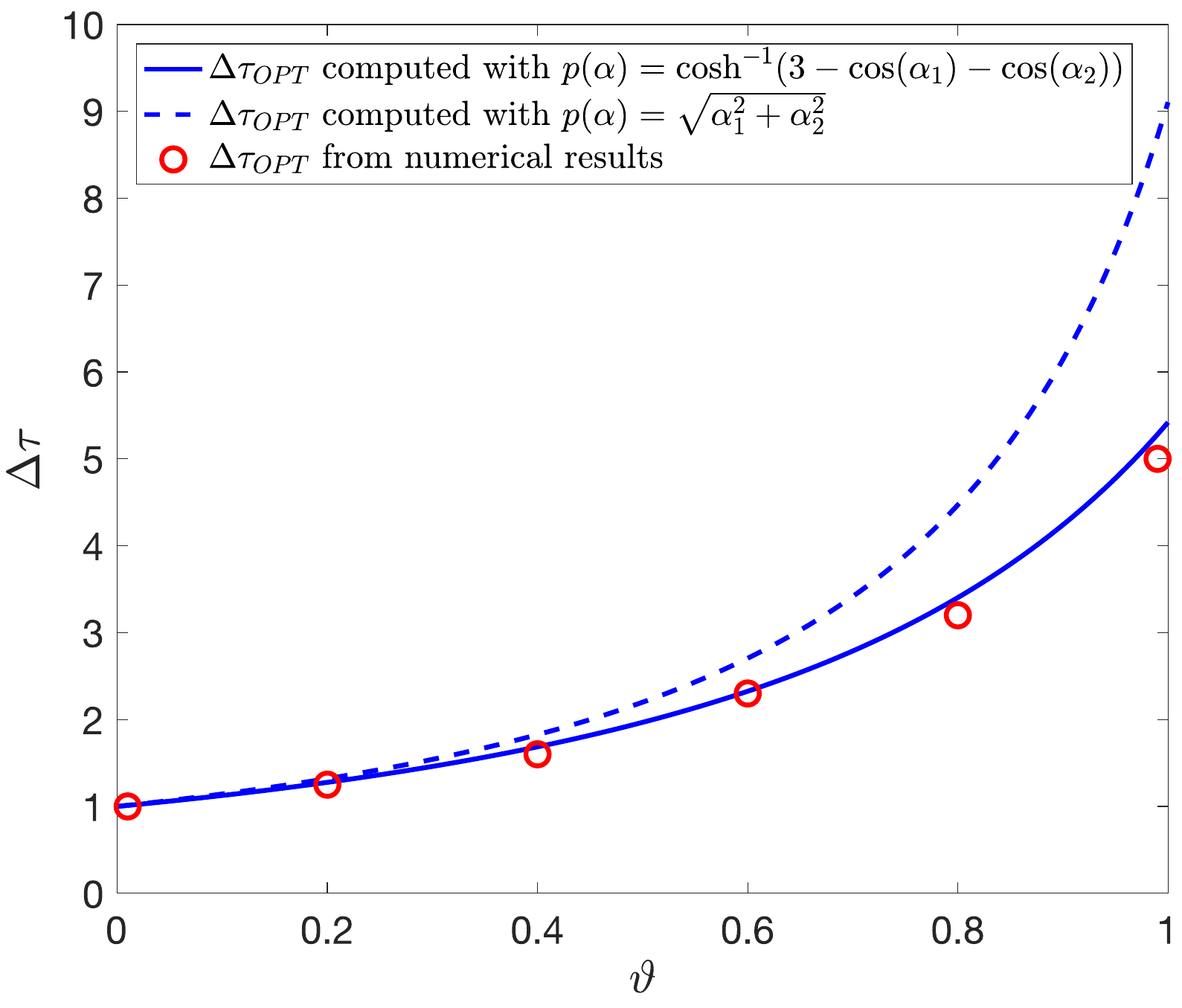}
&
\includegraphics[width=0.45\textwidth]{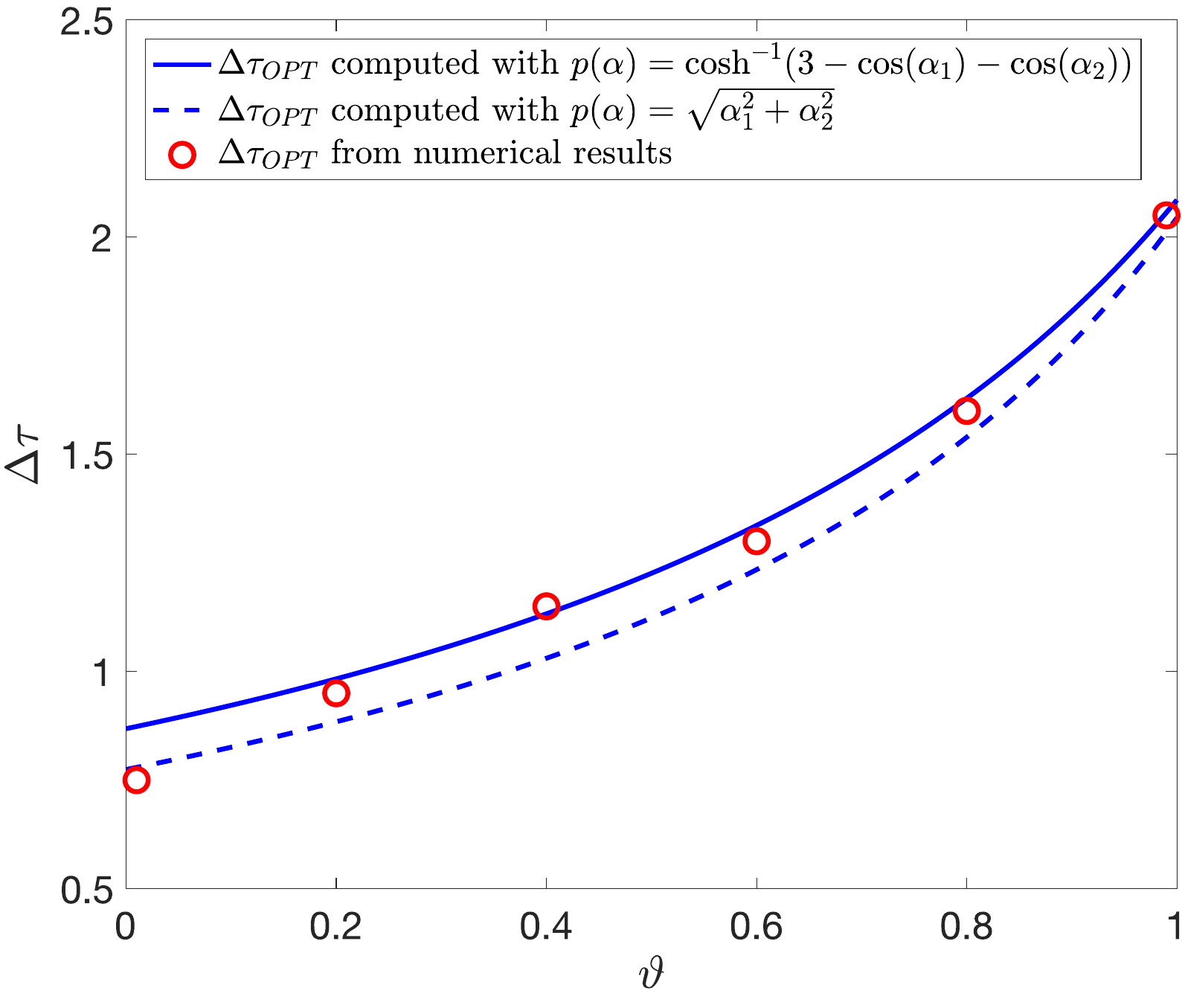}
\end{tabular}
\end{center}
\caption{Optimal relaxation parameter $\Delta \tau$ for the rectangular domain of \S\ref{sect:rectangle3D} from numerical experiments (red points, inferred from Fig.~\ref{fig:3D_rho}) and from the BLFA, namely Eqs.~\eqref{dtauoptLFA3D} and \eqref{dtauoptLFA3D_Neumann} for Dirichlet (left plot) and Neumann (right plot) boundary conditions, respectively. The dashed lines are the theoretical $\Delta \tau$ if the solutions to the continuous problems were used instead of the discretized ones (see Remark \ref{solnum3D}).
}
\label{fig:dtau_opt3D}
\end{figure}

\begin{remark}\label{solnum3D}
A similar argument of Remark \ref{solnum2D} applies to the 3D case, where the exact solution is obtained with $p(\alpha_1,\alpha_2) = \alpha_1 + \alpha_2$. The discrepancy between the optimal $\Delta \tau_\text{OPT}$ is observed in Fig.~\ref{fig:dtau_opt3D}.
\end{remark}

\subsubsection{Neumann boundary conditions: BLFA}\label{3D_neuLFA}
The BLFA for the Neumann boundary condition is analogous to \S\ref{sect:3D_dirLFA} for the 3D aspect of the analysis and to \S\ref{sect:2D_neuLFA} for the treatment of the boundary condition. Therefore, using the same approaches the following proposition can be easily proved for the Neumann case.

\begin{proposition}\label{prop:opt3DN}
The optimal $\Delta \tau$ for Neumann boundary conditions in 3D is obtained by the following formula:
\begin{equation}\label{dtauoptLFA3D_Neumann}
\frac{\Delta \tau_\text{OPT} (\vartheta)}{h} = \frac{2}{G_0^\text{min}(\vartheta)+G_0^\text{max}(\vartheta)}
\end{equation}
with $G_0^\text{min}(\vartheta)$ and $G_0^\text{max}(\vartheta)$ defined as in \eqref{eq:G0minmax}
and $G_0$ defined in~\eqref{eq:G3D} with coefficients \eqref{eq:coeffN}.
\end{proposition}

\begin{figure}
\begin{minipage}{0.99\textwidth}
\begin{center}
(Dirichlet)\\
\includegraphics[width=0.79\textwidth]{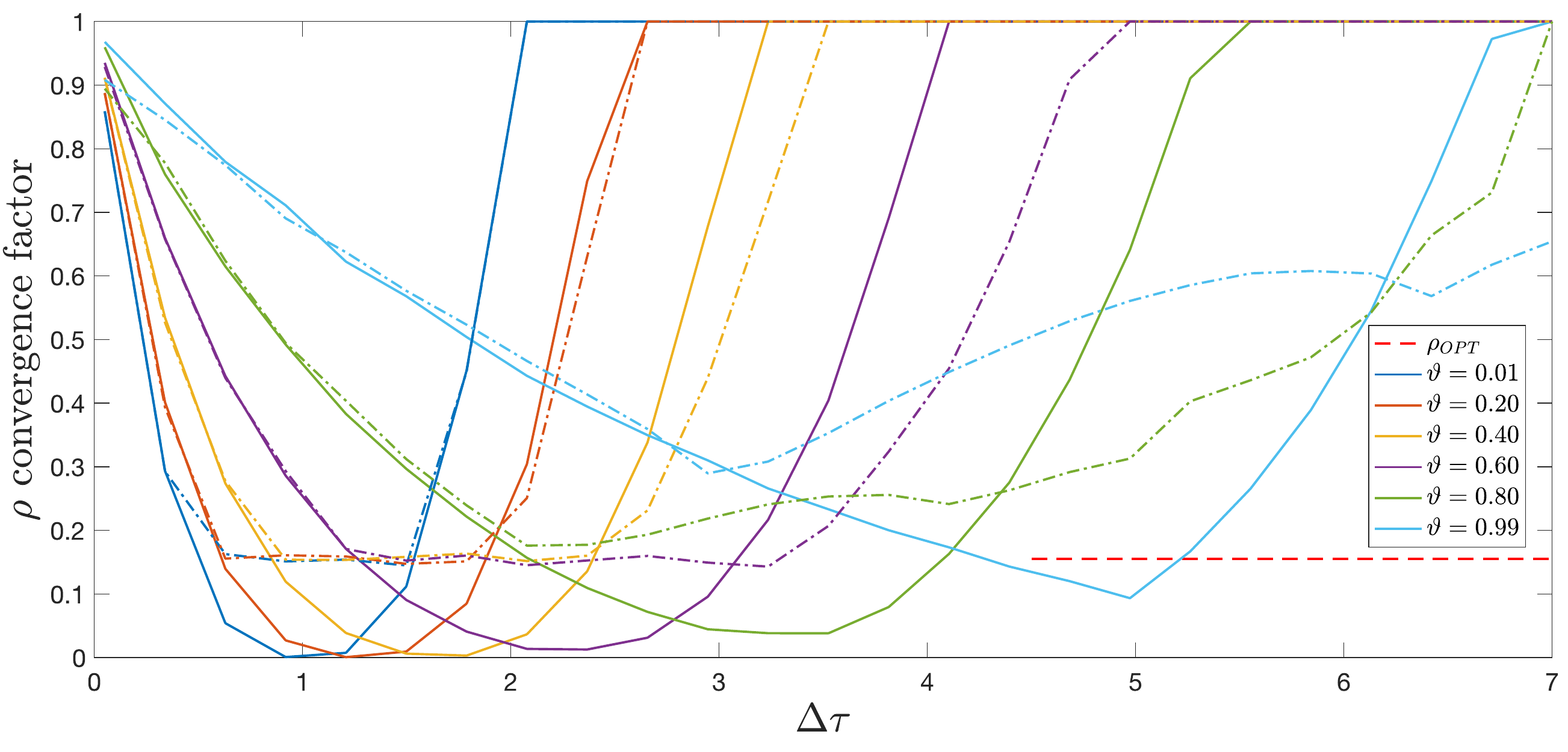}
\end{center}
\end{minipage}
\begin{minipage}{0.99\textwidth}
\begin{center}
(Neumann)\\
\includegraphics[width=0.79\textwidth]{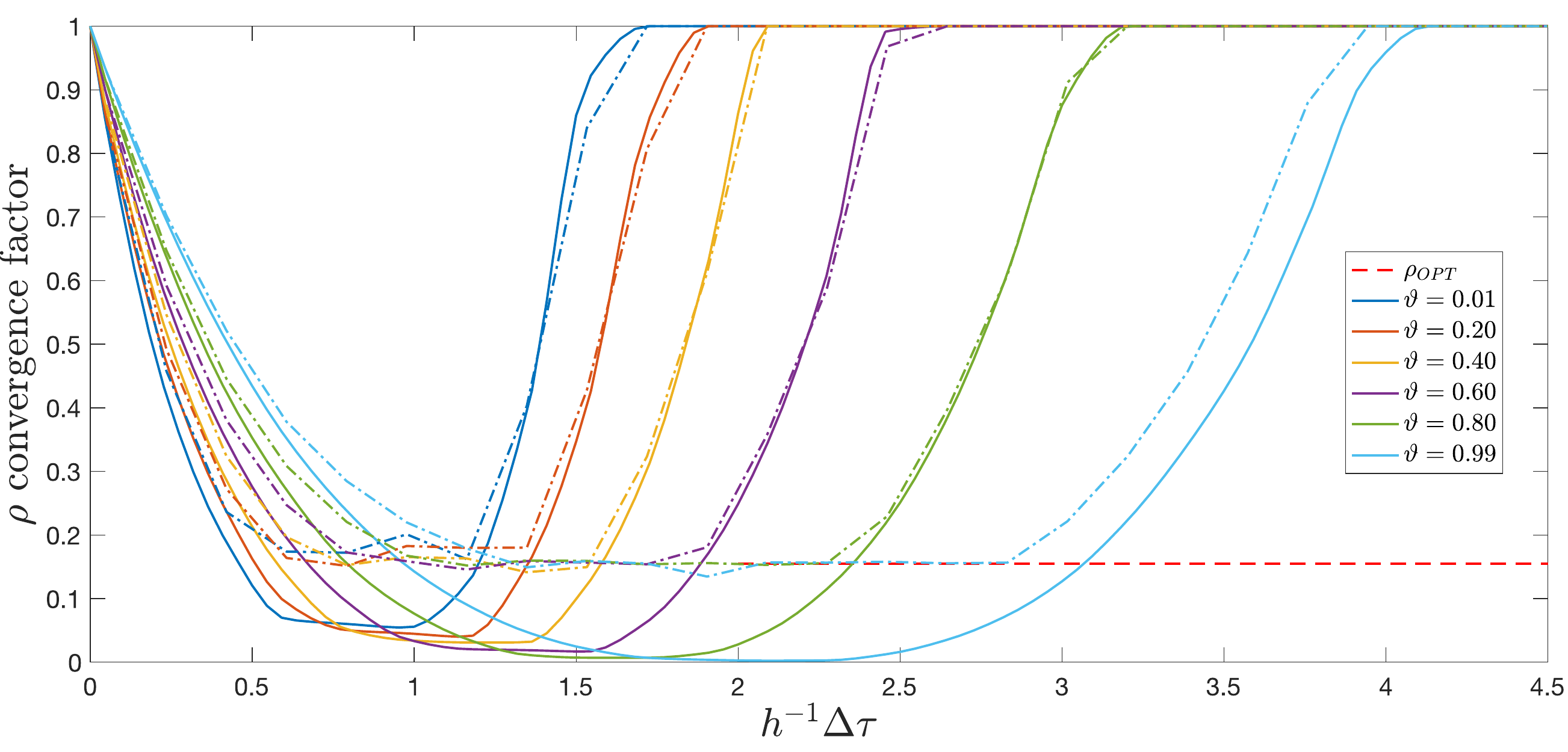}
\end{center}
\end{minipage}
\caption{TGCS convergence factor for 3D Dirichlet (top) and Neumann (bottom) boundary conditions in a rectangular domain (\S\ref{sect:rectangle3D}). The red dashed line is the optimal convergence factor without ghost nodes.}
\label{fig:3D_rho}
\end{figure}

\begin{remark}
In the general $d$-dimensional space, the BLFA is performed on a $(d-1)$-dimensional hyperplane. The exact numerical solution \eqref{eq:exactnum3D} is generalized as
\[
u_{i_1,\ldots,i_d} = e^{p(\alpha_1,\ldots,\alpha_{d-1}) (x_{1,i_1}-1)/h}\varphi(x_{2,i_2},\ldots,x_{d,i_d}) = e^{p(\alpha_1,\ldots,\alpha_{d-1}) (x_{1,i_1}-1)/h} e^{\iu \alpha_1 x_{2,i_2}/h} \cdot \ldots \cdot e^{\iu \alpha_{d-1} x_{d,i_d}/h}
\]
with
\[
p(\alpha_1,\ldots,\alpha_{d-1}) = \cosh^{-1} \left( d-\sum_{i=1}^{d-1} \cos(\alpha_i) \right).
\]
The optimal $\Delta \tau$ is defined as in \eqref{prop:opt3D} and \eqref{prop:opt3DN} for Dirichlet and Neumann boundary conditions, respectively, with
\begin{equation}\label{eq:G0minmax_dD}
G_0^\text{min/max}(\vartheta) = \min / \max \left\{ |G_0(A_1,\vartheta)|,
|G_0(A_2,\vartheta)|,
\ldots,
|G_0(A_{d},\vartheta)|
\right\},
\end{equation}
where $A_1,\ldots,A_d$ are $(d-1)$-dimensional points:
\[
A_1 = (\pi/2,0,\ldots,0), \quad
A_2 = (\pi/2,\pi/2, 0, \ldots, 0), \quad
\ldots \quad,
A_{d-1} = (\pi/2,\pi/2, \ldots, \pi/2), \quad
A_{d} = (\pi,\pi, \ldots, \pi).
\]
\end{remark}

\subsection{Non-rectangular domains}
The approach is easily extended from the 2D treatment of general domains of \S\ref{sect:nonrec}. Relaxation of boundary conditions is performed by \eqref{eq:rel} with \eqref{DtauOptTheta}. The optimal relaxation parameter is determined by \eqref{dtauoptLFA3D} and \eqref{dtauoptLFA3D_Neumann} for Dirichlet and Neumann boundary conditions, respectively.
Also, we use the same approach suggested in Remark \ref{exptheta2D} for the 2D case.

\subsubsection{Numerical results in 3D: Dirichlet boundary conditions}\label{test:nonrec3D_D}
In this section we test the multigrid efficiency on an ellipsoidal domain with Dirichlet boundary conditions. The level-set function is
\[
\phi(x,y,z) = \left( \frac{x-\sqrt{2}/20}{0.686} \right)^2+\left( \frac{y-\sqrt{3}/40}{0.386} \right)^2+\left( \frac{z+\sqrt{2}/50}{0.586} \right)^2-1.
\]

In the left panel of Fig.~\ref{fig:3D_RHO_ELLIPSOID_DandN} we plot the convergence factors of the TGCS (with set-up from \S\ref{sect:setup}) for different values of $N$. The solid line represents the test when $\Delta \tau$ is constant for all ghost points, while the dashed line is the test with the optimal $\Delta \tau (\tilde{\vartheta})$ defined in~\eqref{DtauOptTheta} with \eqref{dtauoptLFA3D}. If we choose a constant $\Delta \tau$ for all ghost points, the best convergence factor that we can achieve is $\rho \approx 0.5$ obtained with $\Delta \tau \approx 1.75$ (a possible explanation would be based on the same argument that was discussed in Remark~\ref{remark:dtau} for the 2D case). A variable $\Delta \tau (\tilde{\vartheta})$ defined in \eqref{DtauOptTheta} with \eqref{dtauoptLFA3D} improves the efficiency, approaching the optimal convergence factor for the standard case $\Omega = [-1,1]^3$, that is
$\rho_\text{OPT} \approx 0.155$ (see \eqref{eq:rhoopt3D}).

\subsubsection{Numerical results in 3D: Neumann boundary conditions}\label{test:nonrec3D_N}
\begin{figure}
\begin{center}
\begin{minipage}{0.49\textwidth}
\includegraphics[width=0.99\textwidth]{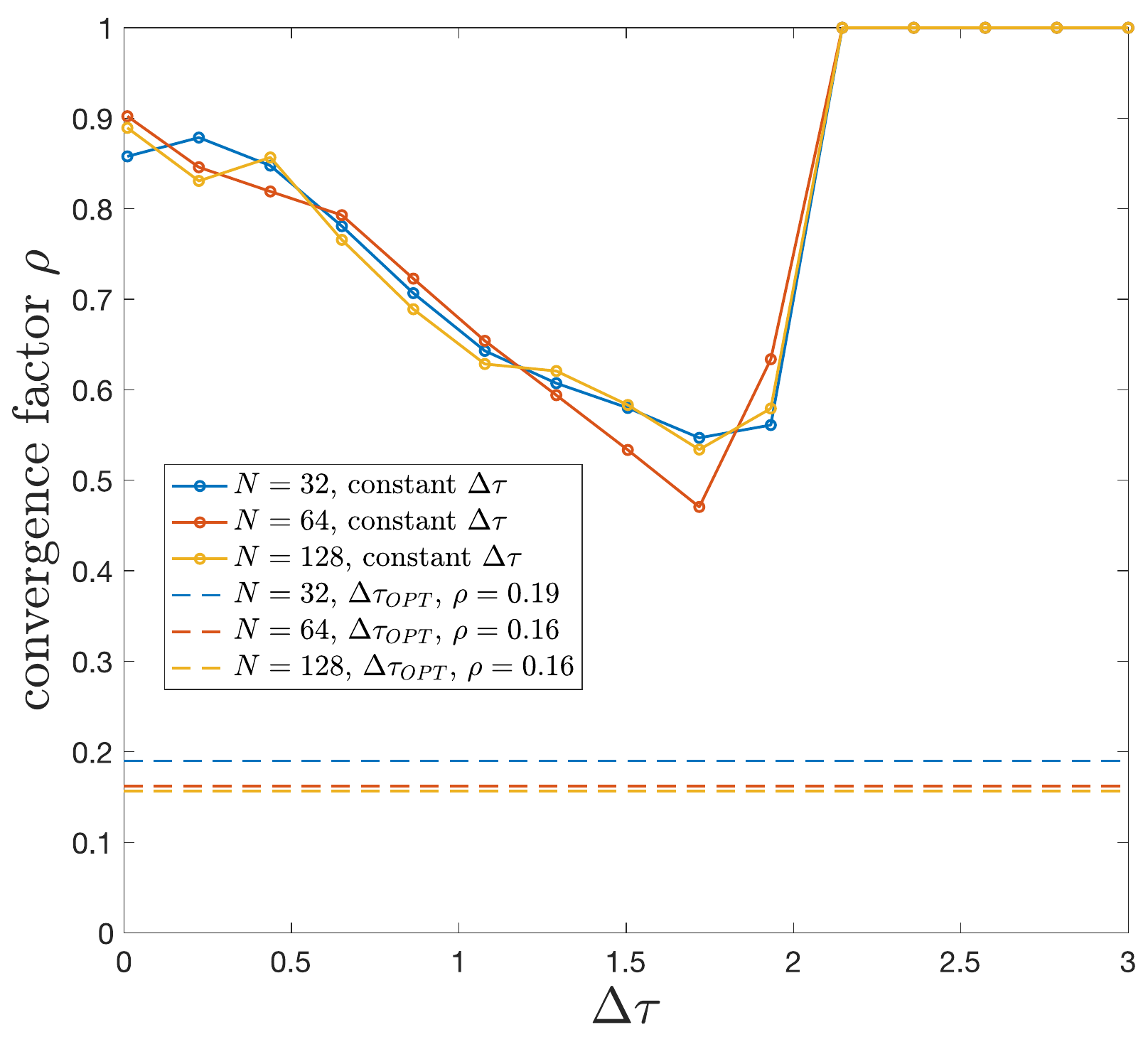}
\end{minipage}
\begin{minipage}{0.49\textwidth}
\includegraphics[width=0.99\textwidth]{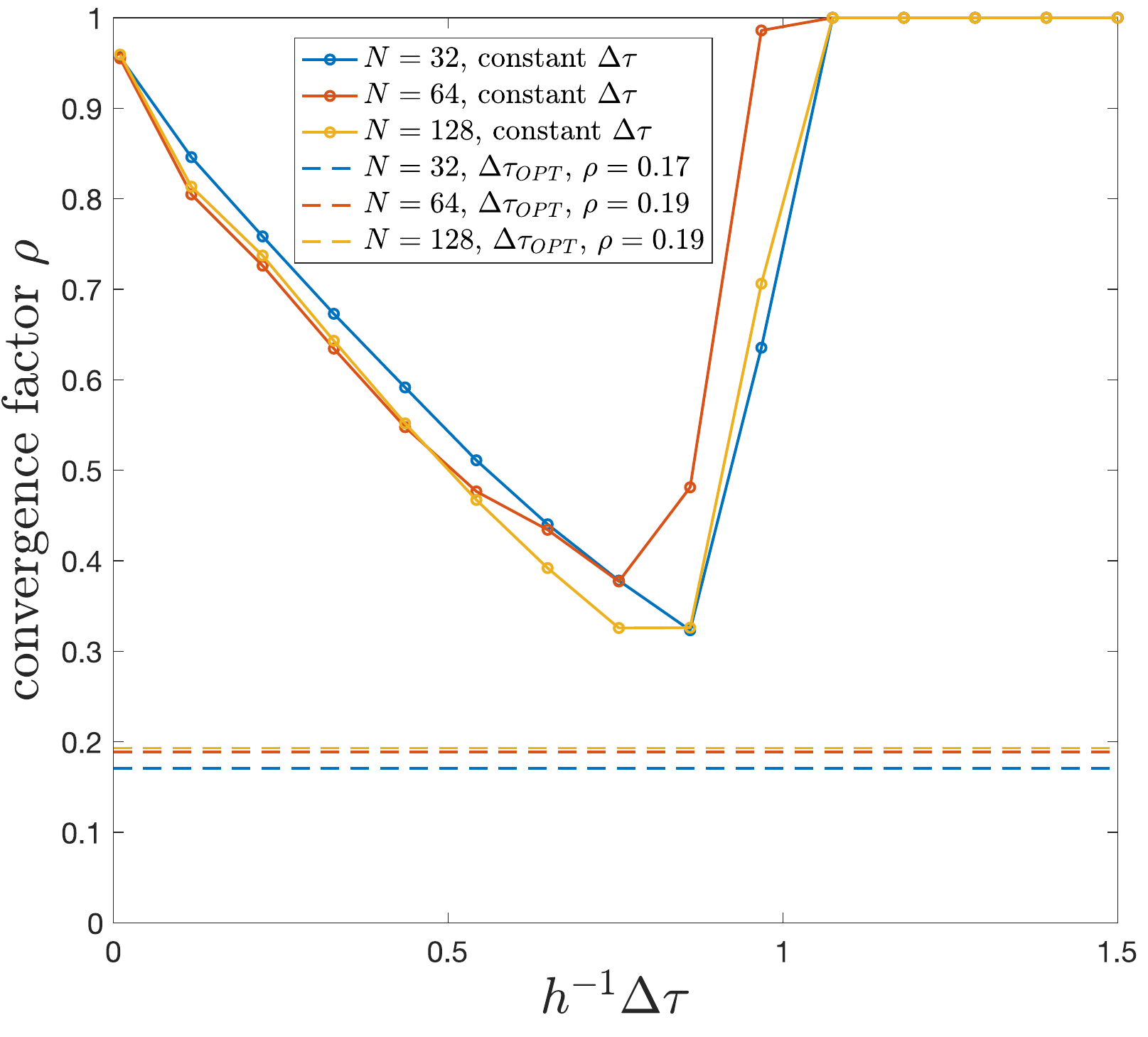}
\end{minipage}
\end{center}
\caption{
Asympototic convergence factors of the TGCS for the three-dimensional problem in an ellipsoidal domain. Left: Dirichlet boundary conditions (\S\ref{test:nonrec3D_D}). Right: Neumann boundary conditions on a problem with the reaction term (\S\ref{test:nonrec3D_N}).
}
\label{fig:3D_RHO_ELLIPSOID_DandN}
\end{figure}

Similarly to \S\ref{test:nonrec_N} (\textsc{Test2N}) for the 2D case, 
we solve a reaction-diffusion problem with Neumann boundary conditions on the whole boundary. The domain is the same ellipsoid of \S\ref{test:nonrec3D_D}.
In the right panel of Fig.~\ref{fig:3D_RHO_ELLIPSOID_DandN} we plot the convergence factors of the TGCS (with set-up from \S\ref{sect:setup}). Also in this case we observe an improvement of the convergence factor when the optimal
$\Delta \tau (\tilde{\vartheta})$ is chosen (defined in \eqref{DtauOptTheta} with \eqref{dtauoptLFA3D_Neumann}). 
The efficiency is considerably improved with respect to the case of constant $\Delta \tau$. The discrepancy of results between constant and optimal $\Delta \tau$ strategies is even more pronounced than the 2D case. Therefore, the implication of adopting an optimal $\Delta \tau$ that varies according to the distance of ghost points from the boundary remains a reasonable strategy to achieve an appropriate efficiency in view of future applications, such as AMR on larger 3D problems.

\section{Conclusions}\label{sec:concl}

We have presented a Boundary Local Fourier Analysis (BLFA) theory to design an efficient relaxation scheme that is able to smooth the residual of the boundary conditions along the tangential direction, without affecting the smoothing performance of the internal equations. 
The smoother is embedded in a multigrid framework and the convergence factor is not degraded from boundary effects.

The model problem is addressed with unfitted boundary methods and a ghost-point technique is adopted to discretize the boundary conditions.

It is known that equations for the boundary conditions pose problems to iterative solvers since their residuals behave differently than the residuals for internal points. In the context of smoothers for multigrid methods,
boundary conditions can be relaxed by using a fictitious time evolution, where the time step $\Delta \tau$ is a relaxation parameter. 

In~\cite{CocoRusso:Elliptic} the authors adopted a constant relaxation parameter for all ghost points and performed some over-relaxations on the boundary in order to mitigate the boundary effects. In this paper, we substantially improve the efficiency by adopting a variable $\Delta \tau$ that depends on the distance from the ghost point to the boundary and is able to 
smooths the residual along the tangential direction. The optimal value is found theoretically by a BLFA and confirmed from numerical experiments. The improvement of the efficiency is justified when the majority of grid points is located in the vicinity of the boundary, such as in the AMR framework. In this paper we perform the analysis on a uniform Cartesian grid and we leave the application to AMR approaches to future efforts.

For rectangular domains where a boundary does not pass through grid points, the prediction from the BLFA is excellent, both for Dirichlet and Neumann boundary conditions. For curved boundary, there might be a small degradation of the convergence factor and this becomes more evident for complex-shaped domains such as the flower-shaped one. We believe that this degradation will vanish when the grid is sufficiently refined. Also, the theoretical argument that is adopted for non-rectangular domain is a heuristic extension of the analysis performed for rectangular domains. Improvement of the theoretical parameters for the non-rectangular case can be performed by a more comprehensive analysis that accounts, for example, for the local curvature of the boundary. 

We finally point out that the BLFA method presented here does not rely on the internal smoother, namely on the smoother adopted for the internal equation (weighted-Jacobi, Gauss-Seidel Red-Black, Gauss-Seidel LEX, etc.). In other words, 
the scope of this paper is to provide a general method to find optimal relaxation parameters for the boundary conditions regardless of the other ingredients of the multigrid
and such that the convergence factor is not degraded from boundary effects. For a specific problem, a proper choice of the other multigrid components (internal smoother, transfer operators, etc.) must be performed carefully in order to optimize the efficiency. Also, the quantitative values shown in this paper, such as those of Remark \ref{remark:dtau}, strongly depend on these choices. Finally, the BLFA will only have the task of embedding the relaxation of boundary conditions into the multigrid method while maintaining the optimal convergence factor.

\section*{Acknowledgments}
All authors acknowledge support from GNCS--INDAM (National Group for Scientific Computing, Italy).


\bibliographystyle{amsalpha}
\bibliography{bibliography}
\end{document}